\def\ifundefined#1#2{\expandafter\ifx\csname#1\endcsname\relax\input #2\fi}
\def\Ma{{\cal\char'101}}
\def\Mb{{\cal\char'102}}
\def\Mc{{\cal\char'103}}
\def\Md{{\cal\char'104}}
\def\Me{{\cal\char'105}}
\def\Mf{{\cal\char'106}}
\def\Mg{{\cal\char'107}}
\def\Mh{{\cal\char'110}}
\def\Mi{{\cal\char'111}}
\def\Mj{{\cal\char'112}}
\def\Mk{{\cal\char'113}}
\def\Ml{{\cal\char'114}}
\def\Mm{{\cal\char'115}}
\def\Mn{{\cal\char'116}}
\def\Mo{{\cal\char'117}}
\def\Mp{{\cal\char'120}}
\def\Mq{{\cal\char'121}}
\def\Mr{{\cal\char'122}}
\def\Ms{{\cal\char'123}}
\def\Mt{{\cal\char'124}}
\def\Mu{{\cal\char'125}}
\def\Mv{{\cal\char'126}}
\def\Mw{{\cal\char'127}}
\def\Mx{{\cal\char'130}}
\def\My{{\cal\char'131}}
\def\Mz{{\cal\char'132}}

\def\pmb#1{\setbox0=\hbox{$#1$}       
     \kern-.025em\copy0\kern-\wd0
     \kern.05em\copy0\kern-\wd0
     \kern-.025em\box0}


\def\endproof{$\hfill \square$}




\def\Cross{\bigm| \kern-5.5pt \not \ \, }
\def\cross{\mid \kern-5.0pt \not \ \, }             
\def\notto{\hbox{$~\rightarrow~\kern-1.5em\hbox{/}\ \ $}}

\def\al{\alpha}

\def\sg{\sigma}

\hyphenation{math-ema-ticians}
\hyphenation{pa-ra-meters}
\hyphenation{pa-ra-meter}
\hyphenation{lem-ma}
\hyphenation{lem-mas}
\hyphenation{to-po-logy}
\hyphenation{to-po-logies}
\hyphenation{homo-logy}
\hyphenation{homo-mor-phy}

\def\nSigma{\Sigma \kern-8.3pt \bigm|\,}

\def\got#1{\hbox{\teneuler #1}}

\font\teneufm=eufm10
\font\eighteufm=eufm8
\font\fiveeufm=eufm5

\newfam\eufam
\textfont\eufam=\teneufm
\scriptfont\eufam=\eighteufm
\scriptscriptfont\eufam=\fiveeufm

\def\got{\fam=\eufam\teneufm}

\def\boxit#1{\vbox{\hrule\hbox{\vrule\kern2.0pt
       \vbox{\kern2.0pt#1\kern2.0pt}\kern2.0pt\vrule}\hrule}}

\def\vlra#1{\hbox{\kern-1pt
       \hbox{\raise2.38pt\hbox{\vbox{\hrule width#1 height0.26pt}}}
       \kern-4.0pt$\rightarrow$}}

\def\vlla#1{\hbox{$\leftarrow$\kern-1.0pt
       \hbox{\raise2.38pt\hbox{\vbox{\hrule width#1 height0.26pt}}}}}

\def\vlda#1{\hbox{$\leftarrow$\kern-1.0pt
       \hbox{\raise2.38pt\hbox{\vbox{\hrule width#1 height0.26pt}}}
       \kern-4.0pt$\rightarrow$}}

\def\longra#1#2#3{\,\lower3pt\hbox{${\buildrel\mathop{#2}
\over{{\vlra{#1}}\atop{#3}}}$}\,}

\def\longla#1#2#3{\,\lower3pt\hbox{${\buildrel\mathop{#2}
\over{{\vlla{#1}}\atop{#3}}}$}\,}

\def\longda#1#2#3{\,\lower3pt\hbox{${\buildrel\mathop{#2}
\over{{\vlda{#1}}\atop{#3}}}$}\,}

\def\overrightharpoonup#1{\vbox{\ialign{##\crcr
	$\rightharpoonup$\crcr\noalign{\kern-1pt\nointerlineskip}
	$\hfil\displaystyle{#1}\hfil$\crcr}}}

\catcode`@=11
\def\@@dalembert#1#2{\setbox0\hbox{$#1\rm I$}
  \vrule height.90\ht0 depth.1\ht0 width.04\ht0
  \rlap{\vrule height.90\ht0 depth-.86\ht0 width.8\ht0}
  \vrule height0\ht0 depth.1\ht0 width.8\ht0
  \vrule height.9\ht0 depth.1\ht0 width.1\ht0 }
\def\dalembert{\mathord{\mkern2mu\mathpalette\@@dalembert{}\mkern2mu}}

\def\@@varcirc#1#2{\mathord{\lower#1ex\hbox{\m@th${#2\mathchar\hex0017 }$}}}
\def\varcirc{\mathchoice
  {\@@varcirc{.91}\displaystyle}{\@@varcirc{.91}\textstyle}
{\@@varcirc{.45}\scriptscriptstyle}}
\catcode`@=12

\font\tensf=cmss10 \font\sevensf=cmss8 at 7pt
\newfam\sffam
\textfont\sffam=\tensf\scriptfont\sffam=\sevensf

\input amssym.def
\input amssym

\magnification=1200

\font\bigsll=cmsl10 scaled\magstep3
\tolerance=500
\overfullrule=0pt
\centerline{\bigsll Integral models in unramified mixed characteristic (0,2) of }
\medskip\medskip
\centerline{\bigsll hermitian orthogonal Shimura varieties of PEL type, Part I}
\bigskip\medskip\noindent
\centerline{Adrian Vasiu, Binghamton University, March 26, 2012}
\bigskip
\centerline{Final version to appear in J. Ramanujan Math. Soc.}
\bigskip\smallskip\noindent
{\bf Abstract.} Let $(G,\Mx)$ be a Shimura variety of PEL type such that $G_{{\bf Q}_2}$ is a split ${\bf GSO}_{2n}$ group with $n\ge 2$. We  prove the existence of the integral canonical models of ${\rm Sh}(G,\Mx)/H_2$ in unramified mixed characteristic $(0,2)$, where $H_2$ is a hyperspecial subgroup of $G({\bf Q}_2)$. 
\bigskip\smallskip\noindent
{\bigsll {\bf 1. Introduction}}
\bigskip
Let $k$ be an algebraically closed field of characteristic $2$. Let $W(k)$ be the ring of Witt vectors with coefficients in $k$. Let $\sg$ be the Frobenius automorphism of either $k$ or $W(k)$. Let $D$ be a {\it $2$-divisible group} over $k$. Let $(M,\Phi)$ be the (contravariant) {\it Dieudonn\'e module} of $D$. We recall that $M$ is a free $W(k)$-module of rank equal to the height $h_D$ of $D$ and $\Phi:M\to M$ is a $\sg$-linear endomorphism such that we have $2M\subseteq\Phi(M)\subseteq M$. Let $D^{\rm t}$ be the {\it Cartier dual} of $D$. 
\medskip\smallskip\noindent
{\bf 1.1. Standard deformation spaces.}
The simplest {\it formal deformation spaces} associated to $2$-divisible groups over $k$ are the following three:
\medskip 
{\bf (a)} the formal deformation space ${\got D}$ of $D$ over ${\rm Spf}\,W(k)$;
\smallskip
{\bf (b)} the formal deformation space ${\got D}_-$ of $(D,\lambda_D)$ over ${\rm Spf}\,W(k)$, where $\lambda_D:D\tilde\to D^{\rm t}$ is an (assumed to exist) isomorphism that is a {\it principal quasi-polarization} of $D$;
\smallskip
{\bf (c)} the formal deformation space ${\got D}_+$ of $(D,b_D)$ over ${\rm Spf}\,W(k)$, where $b_D:D\tilde\to D^{\rm t}$ is an (assumed to exist) isomorphism such that the perfect, bilinear form $b_M:M\times M\to W(k)$ induced naturally by $b_D$ is symmetric. \medskip
It is well known that both ${\got D}$ and ${\got D}_-$ are formally smooth over ${\rm Spf}\,W(k)$, cf. Grothendieck--Messing {\it deformation theory} (see [29, Chs. 4--5], [22, Cor. 4.8], etc.). 
\smallskip
In this paragraph we refer to (c). We assume that $h_D=2n$ with $n\in {\bf N}^\ast$. It is easy to check that if $(D,b_D)$ has a lift to ${\rm Spf}\,W(k)$ (equivalently to ${\rm Spec}\,W(k)$, cf. [29, Ch. II, Lem. 4.16]), then $b_M$ modulo $2W(k)$ is alternating. If $b_M$ modulo $2W(k)$ is alternating, then the subgroup scheme of ${\bf GL}_{M/2M}$ that fixes $b_M$ modulo $2W(k)$ is an ${\bf Sp}_{2n}$ group and thus the dimension of the tangent space of ${\got D}_+$ is ${{n(n+1)}\over 2}$. This dimension is greater than the dimension ${{n(n-1)}\over 2}$ predicted by geometric considerations in characteristic $0$ (on {\it hermitian symmetric domains} of type D III). Thus ${\got D}_+$ is not formally smooth over ${\rm Spf}\,W(k)$ and moreover the Grothendieck--Messing deformation theory does not provide a good understanding of ${\got D}_+$ (this is specific to characteristic $2$!). Thus, in the study of good formal subspaces of ${\got D}_+$, one encounters a serious difficulty. The difficulty splits naturally in three problems (parts) that can be described as follows. 
\medskip
{\bf (i)} Determine if $b_M$ modulo $2W(k)$ is or is not alternating.
\smallskip
{\bf (ii)} If $b_M$ modulo $2W(k)$ is alternating, then show that there exists a closed formal subscheme $\bar{\got D}_{++}$ of ${\got D}_+$ that has the following three properties: (ii.a) it is formally smooth over ${\rm Spf}\,k$ of dimension ${{n(n-1)}\over 2}$, (ii.b) a ${\rm Spf}\,k[[t]]$-valued point $y$ of ${\got D}_+$ factors through $\bar{\got D}_{++}$ if and only if the crystalline realization of the isomorphism of Barsotti--Tate groups of level $1$ over  ${\rm Spf}\,k[[t]]$ that is associated to $y$ and lifts $b_D[2]$, is alternating, and (ii.c) it lifts to a closed formal subscheme ${\got D}_{++}$ of ${\got D}_+$ which is formally smooth over ${\rm Spf}\,W(k)$. 
\smallskip
{\bf (iii)} If $\bar{\got D}_{++}$ exists, then show that one can choose ${\got D}_{++}$ in such a way that it has all the desired geometric interpretations. For instance, one required property is: if $V$ is a discrete valuation ring which is a finite extension of $W(k)$, then for each ${\rm Spf}\,V$-valued point of ${\got D}_{++}$, the corresponding Barsotti--Tate group of level $1$ over ${\rm Spec}\,V[{1\over 2}]$ is naturally endowed with a principal quasi-polarization (an alternating bilinear form). 
\medskip
If $\bar{\got D}_{++}$ exists, then properties (ii.a) and (ii.b) determine it uniquely. But if the {\it $2$-rank} of $D$ is positive, then the lift ${\got D}_{++}$ of $\bar{\got D}_{++}$ mentioned in (ii.c) is not unique. For instance, if $D$ is ordinary, then ${\got D}$ has a canonical structure of a formal torus and we can replace ${\got D}_{++}$ by its translate through a $2$-torsion ${\rm Spf}\,W(k)$-valued point of ${\got D}$ that factors through ${\got D}_+\setminus {\got D}_{++}$. This is why we also require in general the property (iii).
\smallskip
The formal deformation spaces ${\got D}$, ${\got D}_-$, and ${\got D}_+$  pertain naturally to the study of integral models in mixed characteristic $(0,2)$ of the simplest cases of {\it unitary}, {\it symplectic}, and {\it hermitian orthogonal} (respectively) {\it Shimura varieties} of {\it PEL type}. Shimura varieties of PEL type are moduli spaces of polarized abelian schemes endowed with suitable ${\bf Z}$-algebras of endomorphisms and with level structures and this explains the PEL type terminology (see [44], [27], [26], and [30]). 
\smallskip
Shimura varieties of PEL type are the simplest examples of Shimura varieties of abelian type (see [31]). The understanding of the zeta functions of Shimura varieties of abelian type depends on a good theory of their integral models. Such a theory was obtained in [39] and [41] (resp. in [42] and [24]) for cases of good reduction with respect to primes of characteristic at least $5$ (resp. at least $3$). But for refined applications to zeta functions one needs also a good theory in mixed characteristic $(0,2)$. Therefore we report here on recent progress towards such a theory.  
\medskip\smallskip\noindent
{\bf 1.2. Previous works.} 
Let ${\bf Z}_{(2)}$ be the localization of ${\bf Z}$ at the prime ideal $2{\bf Z}$. The previous status of the existence of smooth integral models of quotients of Shimura varieties of PEL (or even abelian) type in mixed characteristic $(0,2)$ can be summarized as follows.
\medskip\noindent
{\bf 1.2.1.} Mumford proved the existence of the moduli ${\bf Z}_{(2)}$-scheme $\Ma_{n,1,l}$ that parametrizes principally polarized abelian schemes over ${\bf Z}_{(2)}$-schemes which are of relative dimension $n$ and which are endowed with a level-$l$ symplectic similitude structure (see [34, Thms. 7.9 and 7.10]); here $n\in{\bf N}^{\ast}$ and $l\in 1+2{\bf N}^{\ast}$. The proof uses geometric invariant theory and standard deformations of abelian varieties. Artin's algebraization method can recover Mumford's result (see [1], [2], and [15, Ch. I, Subsect. 4.11]). 
\medskip\noindent
{\bf 1.2.2.} Drinfeld constructed good moduli spaces of $2$-divisible groups over $k$ of dimension $1$ (see [13]). See [32] and [20] for applications of them to compact, unitary Shimura varieties related to ${\bf SU}(1,n)$ groups over ${\bf R}$ (with $n\ge 1$): they provide proper, smooth integrals models of quotients of simple unitary Shimura varieties over localizations of rings of integers of number fields with respect to arbitrary primes of characteristic 2.
\medskip\noindent
{\bf 1.2.3.} Using Mumford's result, Serre--Tate deformation theory (see [29], [22], and [23]), and Grothendieck--Messing deformation theory, in [44] and [27] it is proved the existence of good integral models of quotients of unitary and symplectic Shimura varieties of PEL type in unramified mixed characteristic $(0,2)$. These integral models are finite schemes over $\Ma_{n,1,l}$, are smooth over ${\bf Z}_{(2)}$, and are moduli spaces of principally polarized abelian schemes which are of relative dimension $n$ and which are endowed with suitable ${\bf Z}$-algebras of endomorphisms and with level-$l$ symplectic similitude structures.
\medskip\noindent
{\bf 1.2.4.} Using the results of Subsubsection 1.2.3, in [41] it is proved the existence of good integral models of quotients of arbitrary unitary Shimura varieties in unramified mixed characteristic $(0,2)$.
\medskip\noindent
{\bf 1.2.5.} Recently, in [42] and [24] it is proved using new developments in the crystalline theory the existence of good integral models of quotients of Shimura varieties of Hodge type in unramified mixed characteristic $(0,2)$ which are finite schemes over $\Ma_{n,1,l}$ and whose special fibers are endowed with abelian schemes that have $2$-rank $0$ at all points.
\medskip
All the methods of Subsections 1.2.1 to 1.2.5 are {\it not} specific to the prime $2$.  
\medskip\smallskip\noindent
{\bf 1.3. Standard PEL situations.} The goal of this paper and of its subsequent Part II is to complete the proof started by Mumford of the existence of good integral models of Shimura varieties of PEL type in unramified mixed characteristic $(0,2)$. We now introduce the {\it standard PEL situations} used in [44], [27], and [26]. 
\smallskip
We start with a symplectic space $(W,\psi)$ over ${\bf Q}$ and with an injective map $$f:(G,\Mx)\hookrightarrow ({\bf GSp}(W,\psi),\Ms)$$
 of Shimura pairs (see [9], [10], [30, Ch. 1], and [39, Subsect. 2.4]). Here $({\bf GSp}(W,\psi),\Ms)$ defines a Siegel modular variety (see [30, Ex. 1.4]) and $(G,\Mx)$ is a Shimura pair of Hodge type. We identify $G$ with a {\it reductive} subgroup of ${\bf GSp}(W,\psi)$ via $f$. Let ${\bf S}:={\rm Res}_{{\bf C}/{\bf R}} {\bf G}_{m,{\bf C}}$ be the unique two dimensional torus over ${\bf R}$ with the property that ${\bf S}({\bf R})$ is the (multiplicative) group ${\bf G}_m({\bf C})$ of non-zero complex numbers. We recall that $\Mx$ (resp. $\Ms$) is a hermitian symmetric domain whose points are a $G({\bf R})$-conjugacy class of (resp. are the ${\bf GSp}(W,\psi)({\bf R})$-conjugacy class of all) monomorphisms ${\bf S}\hookrightarrow G_{{\bf R}}$ (resp. ${\bf S}\hookrightarrow {\bf GSp}(W\otimes_{\bf Q} {\bf R},\psi)$) over ${\bf R}$ that define Hodge ${\bf Q}$--structures on $W$ of type $\{(-1,0),(0,-1)\}$ and that have either $2\pi i\psi$ or $-2\pi i\psi$ as polarizations. Let $E(G,\Mx)$ be the number field that is the reflex field of $(G,\Mx)$ (see [10] and [30]). Let $v$ be a prime of $E(G,\Mx)$ of characteristic $2$. Let $O_{(v)}$ be the localization of the ring of integers of $E(G,\Mx)$ with respect to $v$ and let $k(v)$ be the residue field of $v$. Let ${\bf A}_f:=\widehat{\bf Z}\otimes_{\bf Z}{\bf Q}$ be the ring of finite ad\`eles. Let ${\bf A}_f^{(2)}$ be the ring of finite ad\`eles with the $2$-component omitted; we have ${\bf A}_f={\bf Q}_2\times{\bf A}_f^{(2)}$. Let $O(G)$ be the set of compact, open subgroups of $G({\bf A}_f)$ endowed with the inclusion relation. Let ${\rm Sh}(G,\Mx)$ be the {\it canonical model} over $E(G,\Mx)$ of the complex Shimura variety (see [9, Thm. 4.21 and Cor. 5.7]; see [10, Cor. 2.1.11] for the identity part which holds as [10, axioms 2.1.1.1 to 2.1.1.5] hold for $(G,\Mx)$)
$${\rm Sh}(G,\Mx)_{{\bf C}}:={\rm proj.}{\rm lim.}_{H\in O(G)} G({\bf Q})\backslash (\Mx\times G({\bf A}_f)/H)=G({\bf Q})\backslash (\Mx\times G({\bf A}_f)).\leqno (1)$$
\indent
Let $L$ be a ${\bf Z}$-lattice of $W$ such that $\psi$ induces a perfect form 
$\psi:L\times L\to{\bf Z}$ i.e., the induced injection $L\hookrightarrow L^\vee:={\rm Hom}(L,{\bf Z})$ is onto. Let $L_{(2)}:=L\otimes_{{\bf Z}} {\bf Z}_{(2)}$. Let $G_{{\bf Z}_{(2)}}$ be the schematic closure of $G$ in the reductive group scheme ${\bf GSp}(L_{(2)},\psi)$. Let $G_{{\bf Z}_2}:=G_{{\bf Z}_{(2)}}\times_{{\bf Z}_{(2)}} {\bf Z}_2$, $K_2:={\bf GSp}(L_{(2)},\psi)({\bf Z}_2)$, $H_2:=G({\bf Q}_2)\cap K_2=G_{{\bf Z}_{(2)}}({\bf Z}_2)$, and
$$\Mb:=\{b\in {\rm End}(L_{(2)})|b\;{\rm is}\; {\rm fixed}\; {\rm by}\; G_{{\bf Z}_{(2)}}\}.$$ Let $G_1$ be the subgroup of ${\bf GSp}(W,\psi)$ that fixes all elements of $\Mb[{1\over 2}]$. Let $\Mi$ be the {\it involution} of ${\rm End}(L_{(2)})$ defined by the identity $\psi(b(u),v)=\psi(u,\Mi(b)(v))$,
where $b\in {\rm End}(L_{(2)})$ and $u$, $v\in L_{(2)}$. As $\Mb=\Mb[{1\over 2}]\cap {\rm End}(L_{(2)})$, we have $\Mi(\Mb)=\Mb$. As the elements of $\Mx$ fix $\Mb\otimes_{\bf Z_{(2)}} {\bf R}$, the involution $\Mi$ of $\Mb$ is positive. Let ${\bf F}$ be an algebraic closure of ${\bf F}_2$. 
\smallskip
We will assume that the following four properties (axioms) hold:${}^1$ $\vfootnote{1}{One can check that the property (iv) implies the property (i).}$
\medskip
{\bf (i)} the $W({\bf F})$-algebra $\Mb\otimes_{{\bf Z}_{(2)}} W({\bf F})$ is a product of {\it matrix} $W({\bf F})$-algebras;
\smallskip
{\bf (ii)} the ${\bf Q}$--algebra $\Mb[{1\over 2}]$ is ${\bf Q}$--simple;
\smallskip
{\bf (iii)} the group $G$ is the identity component of $G_1$;
\smallskip
{\bf (iv)} the flat, affine group scheme $G_{{\bf Z}_{(2)}}$ over ${\bf Z}_{(2)}$ is reductive (i.e., it is smooth and its special fibre is connected and has a trivial unipotent radical).
\medskip
Assumption (iv) implies that $H_2$ is a hyperspecial subgroup of $G({\bf Q}_2)=G_{{\bf Q}_2}({\bf Q}_2)$ (cf. [38, Subsubsect. 3.8.1]) and that $v$ is unramified over $2$ (cf. [31, Cor. 4.7 (a)]). Let $G^{\rm der}$ be the derived group of $G$. As the case when $G$ is a torus is trivial from the point of view of good integral models, in all that follows we will assume that $G$ is not a torus (i.e., $G^{\rm der}$ is non-trivial). Assumption (ii) is not really required: it is inserted here only to ease the presentation. Due to properties (ii) and (iii) and the fact that $G$ is not a torus, one distinguishes the following three possible (and disjoint) cases (see [26, Sect. 7]):
\medskip
{\bf (A)} the group $G_{{\bf C}}^{\rm der}$ is a product of ${\bf SL}_n$ groups with $n\ge 2$ and, in the case $n=2$, the center of $G$ has dimension at least $2$;
\smallskip
{\bf (C)} the group $G_{{\bf C}}^{\rm der}$ is a product of ${\bf Sp}_{2n}$ groups with $n\ge 1$ and, in the case $n=1$, the center of $G$ has dimension $1$;
\smallskip
{\bf (D)} the group $G_{{\bf C}}^{\rm der}$ is a product of ${\bf SO}_{2n}$ groups with $n\ge 2$.
\medskip
We have $G\neq G_1$ if and only if we are in the case (D) i.e., if and only if $G^{\rm der}$ is not simply connected (cf. [26, Sect. 7]). In the case (A) (resp. (C) or (D)), one often says that ${\rm Sh}(G,\Mx)$ is a unitary (resp. a symplectic or a hermitian orthogonal) Shimura variety of PEL type (cf. the description of the intersection group $G_{\bf R}\cap {\bf Sp}(W\otimes_{\bf Q} {\bf R},\psi)$ in [37, Subsects. 2.6 and 2.7]). We are in the case (D) if and only if $\Mb\otimes_{{\bf Z}_{(2)}} {\bf R}$ is a product of matrix algebras over the quaternion ${\bf R}$-algebra ${\bf H}$ (see [37, Subsect. 2.1, (type III)]).
\smallskip
We refer to the quadruple $(f,L,v,\Mb)$ as a standard PEL situation in {\it mixed characteristic $(0,2)$}. Let $\Mm$ be the ${\bf Z}_{(2)}$-scheme which parametrizes isomorphism classes of principally polarized abelian schemes over ${\bf Z}_{(2)}$-schemes that are of relative dimension ${{\dim_{{\bf Q}}(W)}\over 2}$ and that are equipped with compatible level-$l$ symplectic similitude structures for all numbers $l\in 1+2{\bf N}^{\ast}$, cf. Subsubsection 1.2.1. We have a natural identification ${\rm Sh}({\bf GSp}(W,\psi),\Ms)/K_2=\Mm_{\bf Q}$ and an action of ${\bf GSp}(W,\psi)({\bf A}^{(2)}_f)$ on $\Mm$ (see [9, Thm. 4.21], [30], etc.). These symplectic similitude structures and this action are defined naturally via $(L,\psi)$ (see [39, Subsect. 4.1]). We identify ${\rm Sh}(G,\Mx)_{{\bf C}}/H_2=G_{{\bf Z}_{(2)}}({\bf Z}_{(2)})\backslash (\Mx\times G({\bf A}_f^{(2)}))$, cf. (1) and [31, Prop. 4.11]. From this identity and the analogous one for ${\rm Sh}({\bf GSp}(W,\psi),\Ms)_{{\bf C}}/K_2$, we get that the natural morphism ${\rm Sh}(G,\Mx)\to {\rm Sh}({\bf GSp}(W,\psi),\Ms)_{E(G,\Mx)}$ of $E(G,\Mx)$-schemes (see [9, Cor. 5.4]) induces a closed embedding (see also [39, Rm. 3.2.14]) 
$${\rm Sh}(G,\Mx)/H_2\hookrightarrow {\rm Sh}({\bf GSp}(W,\psi),\Ms)_{E(G,\Mx)}/K_2=\Mm_{E(G,\Mx)}.$$
\indent
Let  $\Mn$ be the schematic closure of ${\rm Sh}(G,\Mx)/H_2$ in $\Mm_{O_{(v)}}$. Let $(\Ma,\Lambda)$
be the pull back to $\Mn$ of the universal principally polarized abelian scheme over $\Mm$. 
\smallskip
We fix a monomorphism $O_{(v)}\hookrightarrow W(k)$ and in all that follows the pull backs of (morphisms between) $O_{(v)}$-schemes to ${\rm Spec}\,W(k)$ are with respect to it. For a morphism $y:{\rm Spec}\,k\to\Mn_{W(k)}$ of $W(k)$-schemes let
$$(A,\lambda_A):=y^*((\Ma,\Lambda)\times_{\Mn} \Mn_{W(k)}).$$ 
Let $(M_y,\Phi_y,\lambda_{M_y})$ be the principally quasi-polarized $F$-crystal over $k$ of $(A,\lambda_A)$. Thus $M_y$ is a free $W(k)$-module of rank $\dim_{{\bf Q}}(W)$, $\Phi_y:M_y\to M_y$ is a $\sigma$-linear endomorphism such that $\Phi_y(M_y)\supset pM_y$, and $\lambda_{M_y}$ is a perfect, alternating form on $M_y$ such that we have $\lambda_{M_y}(\Phi_y(u),\Phi_y(v))=2\sigma(\lambda_{M_y}(u,v))$ for all $u,v\in M_y$. Let $M_y^\vee:={\rm Hom}(M_y,W(k))$. We denote also by $\lambda_{M_y}$ the perfect, alternating form on $M_y^\vee$ induced naturally by $\lambda_{M_y}$. We have a natural action of $\Mb\otimes_{{\bf Z}_{(2)}} W(k)$ on $M_y^\vee$ (see Subsubsection 4.1.1). 
\smallskip
It is known that in the cases (A) and (C), $\Mn$ is regular and formally smooth over $O_{(v)}$ (cf. [44, Sect. 3]; see also [27, Sect. 6] and [26, Sect. 5]). Accordingly, in the whole paper we will assume that we are in the case (D). In this paper we study triples of the form $(M_y,\Phi_y,\lambda_{M_y})$ and the closed embedding $\Mn\hookrightarrow\Mm_{O_{(v)}}$. The goal of Part I is to prove the following results.
\medskip\smallskip\noindent
{\bf 1.4. Main Theorem.} {\it Let $(f,L,v,\Mb)$ be a standard PEL situation in mixed characteristic $(0,2)$ such that we are in the case (D). We use the notations of Subsection 1.3 and we assume that $G_{{\bf Z}_2}$ is a split ${\bf GSO}_{2n}$ group scheme with $n\ge 2$ (thus we have $k(v)={\bf F}_2$, cf. [31, Prop. 4.6 and Cor. 4.7 (b)]). Then the following two properties hold:
\medskip
{\bf (a)} there exist isomorphisms $L_{(2)}\otimes_{{\bf Z}_{(2)}} W(k)\tilde\to M_y^\vee$ of $\Mb\otimes_{{\bf Z}_{(2)}} W(k)$-modules that induce symplectic isomorphisms $(L_{(2)}\otimes_{{\bf Z}_{(2)}} W(k),\psi)\tilde\to (M_y^\vee,\lambda_{M_y})$;
\smallskip
{\bf (b)} the $O_{(v)}$-scheme $\Mn$ is regular and formally smooth.} 
\medskip
The locally compact, totally disconnected topological group $G({\bf A}_f^{(2)})$ acts on $\Mn$ continuously in the sense of [10, Subsubsect. 2.7.1]. Thus $\Mn$ is an {\it integral canonical model} of ${\rm Sh}(G,\Mx)/H_2$ over $O_{(v)}$ in the sense of [39, Def. 3.2.3 6)], cf. [39, Ex. 3.2.9 and Cor. 3.4.4]. Due to [43, Cor. 30], $\Mn$ is {\it the unique} integral canonical model of ${\rm Sh}(G,\Mx)/H_2$ over $O_{(v)}$ and thus it is {\it the final object} of the category of smooth integral models of ${\rm Sh}(G,\Mx)/H_2$ over $O_{(v)}$ (here the word smooth is used as in [30, Def. 2.2]). In simpler words, $\Mn$ is {\it the best} smooth integral model of ${\rm Sh}(G,\Mx)/H_2$ over $O_{(v)}$.
\smallskip
If $A$ is an ordinary abelian variety, then a result of Noot implies that the normalization $\Mn^{\rm n}_{W(k)}$ of $\Mn_{W(k)}$ is regular and formally smooth over $W(k)$ at all points in $\Mn^{\rm n}_{W(k)}(k)$ that map to $y\in\Mn_{W(k)}(k)$ (see [35, Cor. 3.8]). This and the next result proved in Subsection 7.2 play key roles in the proof of Theorem 1.4 (b).
\medskip\noindent
{\bf 1.4.1. Proposition.} {\it The following two properties hold:
\medskip
{\bf (a)} If $y\in\Mn_{W(k)}(k)$ is an ordinary point (i.e., if the abelian variety $A$ is ordinary), then there exists a unique lift $z_{\rm can}\in\Mn_{W(k)}(W(k))$ of $y$ such that the abelian scheme $z_{\rm can}^*(\Ma\times_{\Mn} \Mn_{W(k)})$ over $W(k)$ is the canonical lift of $A$.  
\smallskip
{\bf (b)} The ordinary locus of $\Mn_k$ is Zariski dense in $\Mn_k$.}
\medskip
Proposition 1.4.1 (b) points out that the cases covered by Theorem 1.4 (b) are in some sense at the opposite pole to the cases mentioned in Subsubsection 1.2.5. 
We would also like to emphasize that:
\medskip
{\bf (i)}  no particular case of the Main Theorem was known before, and 
\smallskip
{\bf (ii)} as a {\it new feature} compared with the results mentioned in Subsubsection 1.2.5, $\Mn$ itself is regular and formally smooth over $O_{(v)}$ (in [42] and [24] one works with the normalizations of the corresponding schematic closures in Mumford moduli schemes such as $\Mm_{O_{(v)}}$ and one can only prove that these normalizations are regular).   
\medskip\noindent
{\bf 1.4.2. Example.} We assume that $\Mb\otimes_{{\bf Z}_{(2)}} {\bf R}={\bf H}$, that $\Mb\otimes_{{\bf Z}_{(2)}} {\bf Z}_2=M_2({\bf Z}_2)$, and that we have  $\dim_{\bf Q}(W)=4n$ for some $n\ge 2$. Thus $\Mb[{1\over 2}]$ is a definite quaternion algebra over ${\bf Q}$ that splits at $2$ and that has a maximal order $\Mb$ which is a semisimple ${\bf Z}_{(2)}$-algebra; moreover we have a $\Mb[{1\over 2}]$-module $W$ of rank $4n$ such that the non-degenerate alternating form $\psi$ on $W$ defines a positive involution of $\Mb[{1\over 2}]$. As $\Mb\otimes_{{\bf Z}_{(2)}} {\bf R}={\bf H}$, we are in the case (D) and the group $G^{\rm der}_{\bf C}$ is an ${\bf SO}_{2n}$ group. It is well known that $G_{{\bf Q}_2}$ is a form of the split ${\bf GSO}_{2n}$ group over ${\bf Q}_2$. Thus if $G_{{\bf Z}_2}$ is a reductive group scheme, then it is either a split or a non-split ${\bf GSO}_{2n}$ group scheme. The Main Theorem applies only if $G_{{\bf Z}_2}$ is a split ${\bf GSO}_{2n}$ group scheme. 
\smallskip
In this paragraph we assume that $G_{{\bf Q}_2}$ is a split ${\bf GSO}_{2n}$ group. It is easy to see that there exists a ${\bf Z}_2$-lattice $L_2$ of $W\otimes_{\bf Q} {\bf Q}_2$ such that $\psi$ induces a perfect, alternating form $\psi:L_2\times L_2\to {\bf Z}_2$ and the schematic closure of $G_{{\bf Q}_2}$ in the reductive group scheme ${\bf GSp}(L_2,\psi)$ is a split ${\bf GSO}_{2n}$ (and thus also a split reductive) group scheme. We can choose the ${\bf Z}$-lattice $L$ of $W$ such that we have $L_2=L\otimes_{\bf Z} {\bf Z}_2=L_{(2)}\otimes_{{\bf Z}_{(2)}} {\bf Z}_2$. For such a choice of $L$, the group scheme $G_{{\bf Z}_{(2)}}$ is reductive and $G_{{\bf Z}_2}$ is a split ${\bf GSO}_{2n}$ group scheme; therefore the Main Theorem applies. 
\medskip\noindent
{\bf 1.4.3. On the new ideas and contents.} The deformation theories of Subsubsection 1.2.3 do not suffice to show that $\Mn_{W(k)}$ is formally smooth over $W(k)$ at points above $y\in\Mn_{W(k)}(k)$ (see Subsection 1.1 and Subsubsections 4.1.2 and 4.3.2). This explains why the proof of the Main Theorem uses also some of the techniques of [14] and [39], the main results of [35], and the following five new ideas.
\medskip\noindent
{\bf (i)} We get versions of the group theoretical results [26, Lem. 7.2, Cor. 7.3] that pertain to the case (D) in mixed characteristic $(0,2)$ (see Subsections 3.1, 3.3, 3.4, and 5.2).
\smallskip\noindent
{\bf (ii)} We use the crystalline Dieudonn\'e theory of [4] and [5] in order to prove Proposition 5.1 that surpasses (in the geometric context of the Main Theorem) the problem 1.1 (i) and that (together with the idea (i)) is the very essence of the proof of Theorem 1.4 (a).
\smallskip\noindent
{\bf (iii)} We use a modulo $2$ version of Faltings deformation theory [14, Sect. 7] as a key ingredient in the proofs of Propositions 6.7 and 1.4.1 (b) (see Subsections 6.3 to 6.7 and 7.2). Subsections 6.3 to 6.7 surpass (in the geometric context of the Main Theorem) the modulo $2$ version of the problems 1.1 (ii) and (iii). In particular, Proposition 6.7 proves that the reduced scheme $\Mn_{k(v),\rm red}$ of $\Mn_{k(v)}$ is regular and formally smooth over $k(v)$. 
\smallskip\noindent
{\bf (iv)} We use a counting argument of suitable ${\bf Z}_2$-lattices (see Lemma 3.1.4) to show based on [35] that the number of {\it quasi-canonical lifts} ${\rm Spec}\,W(k)\to \Mn_{W(k)}$ of a fixed ordinary point of $\Mn_{W(k)}(k)$ is exactly $2^{{n(n-1)}\over {2}}$ (i.e., it is exactly as predicted by Theorem 1.4 (b) and Proposition 1.4.1). This implies that $\Mn_{W(k)}$ is regular at all ordinary points of $\Mn_{W(k)}(k)$.
\smallskip\noindent
{\bf (v)} We combine the ideas of (iii) and (iv) with a classical {\it lemma of Hironaka} to conclude that $\Mn$ itself is regular (see Subsection 7.3). 
\medskip
We recall that if the abelian variety $A$ is ordinary, then a quasi-canonical lift ${\rm Spec}\,W(k)\to \Mn_{W(k)}$ of $y$ is a lift that produces an abelian scheme over $W(k)$ which lifts $A$ and whose Hodge filtration is defined by the maximal direct summand of $M_y$ normalized by ${1\over 2}\Phi_y$. Ideas (i) and (ii) are specific to the PEL context of Subsection 1.3. Ideas (iii) to (v) can be adapted to all characteristics and to all Shimura varieties of Hodge type. 
\smallskip
Sections 2 and 3 present tools from the crystalline theory and from the theories of group schemes and of involutions of matrix algebras over commutative ${\bf Z}$-algebras (respectively). Section 5 proves Theorem 1.4 (a). Section 6 checks that $\Mn_{k(v),\rm red}$ is regular and formally smooth over $k(v)$. Section 7 proves Theorem 1.4 (b) and Proposition 1.4.1. Section 4 lists notations, properties, and strategies that are used in Sections 5 to 7. 
\medskip\smallskip\noindent
{\bf 1.5. On Part II.} The hypothesis of the Main Theorem that $G_{{\bf Z}_2}$ is a split ${\bf GSO}_{2n}$ group scheme is inserted only to ease notations and to be able to apply the results of [35]. Part II of this paper will prove the Main Theorem without this hypothesis. The proof of the Main Theorem can be easily adapted to the general case. The main idea of Part II will be to use relative PEL situations (similar but different from the ones of [39, Subsubsect. 4.3.16 and Sect. 6]) to get that for proving that $\Mn$ is regular and formally smooth over $O_{(v)}$, we can  assume that $G_{{\bf Z}_2}$ is a split ${\bf GSO}_{2n}$ group scheme. 
\bigskip\smallskip\noindent
{\bigsll {\bf 2. Crystalline preliminaries}}
\bigskip
Let $B(k)$ be the field of fractions of $W(k)$. We denote also by $\sg$ the Frobenius automorphism of $B(k)$. If $Z$ is a $W(k)$-scheme annihilated by some power of $2$, we use Berthelot's crystalline site $CRIS(Z/{\rm Spec}\,W(k))$ of [3, Ch. III, Sect. 4]. See [4] and [5] for the crystalline Dieudonn\'e functor ${\bf D}$. If $\Me$ (resp. $E$) is a $2$-divisible group (resp. a finite, flat group scheme annihilated by $2$) over some $W(k$)-scheme, let $\Me^{\rm t}$ (resp. $E^{\rm t}$) be its Cartier dual. For a morphism of affine schemes ${\rm Spec}\,S_1\to {\rm Spec}\,S$ and for $Z$ (or $Z_S$ or $Z_*$) an $S$-scheme, let $Z_{S_1}$ (or $Z_{S_1}$ or $Z_{*,S_1}$) be $Z\times_S S_1$. A bilinear form $\lambda_M$ on a free $S$-module $M$ of finite rank is called perfect if it induces an $S$-linear isomorphism $M\tilde\to M^{\vee}$ where $M^{\vee}:={\rm Hom}(M,S)$; if $\lambda_M$ is alternating, $(M,\lambda_M)$ is called a symplectic space over $S$. Let $\pmb{\mu}_{2,S}$ be the $2$-torsion subgroup scheme of the rank $1$ split torus ${\bf G}_{m,S}$ over $S$. Let ${\bf G}_m(S)$ be the group of invertible elements of $S$. If $S$ is an ${\bf F}_2$-algebra, let $\Phi_S$ be the Frobenius endomorphism of either $S$ or ${\rm Spec}\,S$ (thus $\Phi_k=\sg$), let $M^{(2)}:=M\otimes_S {}_{\Phi_S} S$, let $Z^{(2)}:=Z\times_S {}_{\Phi_S} S$, and let $\pmb{\al}_{2,S}$ be the finite, flat, group scheme over $S$ of global functions of square $0$.
\medskip\smallskip\noindent
{\bf 2.1. Ramified data.} Let $V$ be a discrete valuation ring that is a finite extension of $W(k)$. Let $K$ be the field of fractions of $V$. Let $e:=[V:W(k)]$. Let $R:=W(k)[[t]]$, with $t$ as an independent variable. We fix a uniformizer $\pi$ of $V$. Let $f_{e}\in R$ be the Eisenstein polynomial of degree $e$ that has $\pi$ as a root. Let $R_e$ be the local $W(k)$-subalgebra of $B(k)[[t]]$ formed by formal power series $\sum_{i=0}^\infty a_it^i$ with the properties that for all $i\in{\bf N}$ we have $b_i:=a_i[{i\over {e}}]!\in W(k)$ and that the sequence $(b_i)_{i\in{\bf N}}$ of elements of $W(k)$ converges to $0$. Let $I_e(1):=\{\sum_{i=0}^\infty a_it^i\in R_e|a_0=0\}$. 
Let $\Phi_{R_e}$ be the Frobenius lift of $R_e$ that is compatible with $\sg$ and that takes $t$ to $t^2$. For $m\in{\bf N}^{\ast}$ let
$$U_m:=k[t]/(t^m)=k[[t]]/(t^m).$$
\indent
Let $S_e$ be the $W(k)$-subalgebra of $B(k)[t]$ generated by $W(k)[t]$ and by the ${{f_{e}^i}\over {i!}}$'s with $i\in{\bf N}$. As $f_{e}$ is an Eisenstein polynomial, $S_e$ is also $W(k)$-generated by $W(k)[t]$ and by the ${{t^{ei}}\over {i!}}$'s with $i\in{\bf N}$. Therefore the $W(k)$-algebra $S_e$ depends (as its notation suggests) only on $e$. The $2$-adic completion of $S_e$ is $R_e$. Thus we have a $W(k)$-epimorphism $R_e\twoheadrightarrow V$ that takes $t$ to $\pi$ and whose kernel has a natural divided power structure. We identify canonically $R_e/I_e(1)=W(k)$, $V=R/(f_e)$, and $V/2V=U_e$. 
\medskip\smallskip\noindent
{\bf 2.2. On finite, flat group schemes.} Let $E_0$ be a finite, flat group scheme annihilated by $2$ over a commutative $W(k)$-algebra $S$; this implies that $E_0$ is commutative. Let $\Delta_{E_0}:E_0\hookrightarrow E_0\times_S E_0$ be the diagonal monomorphism of $S$-schemes.
\medskip\noindent
{\bf 2.2.1. Definition.} A principal quasi-polarization of $E_0$ is an isomorphism $\lambda_{E_0}:E_0\tilde\to E_0^{\rm t}$ over $S$ such that the composite $co_{E_0}:E_0\to \pmb{\mu}_{2,S}$ of $(1_{E_0},\lambda_{E_0})\circ \Delta_{E_0}$ with the coupling morphism $cu_{E_0}:E_0\times_S E_0^{\rm t}\to \pmb{\mu}_{2,S}$ factors through the identity section of  $\pmb{\mu}_{2,S}$. 
\medskip
Let $r\in{\bf N}^{\ast}$. Until Section 3 we will take $S=V$, $E_0$ to be of rank (order) $2^{2r}$, and $\lambda_{E_0}$ to be a principal quasi-polarization. Thus the perfect bilinear form on $E_{0,K}$ induced by $\lambda_{E_0}$ is alternating. Let $(N_0,\phi_0,\upsilon_0,\nabla_0,b_{N_0})$ be the evaluation of ${\bf D}(E_0,\lambda_{E_0})$ (equivalently of ${\bf D}((E_0,\lambda_{E_0})_{U_e})$) at the trivial thickening of ${\rm Spec}\,U_e$. Thus $N_0$ is a free $U_e$-module of rank $2r$, the maps $\phi_0:N_0^{(2)}\to N_0$ and $\upsilon_0:N_0\to N_0^{(2)}$ are $U_e$-linear, $\nabla_0$ is a connection on $N_0$, and $b_{N_0}$ is a perfect bilinear form on $N_0$. We have $\phi_0\circ \upsilon_0=0$ and $\upsilon_0\circ \phi_0=0$ as well as the adjunction identity $b_{N_0}(\phi_0(x),y)=b_{N_0}(x,\upsilon_0(y))^2$ for all $x,y\in N_0$. As $\lambda_{E_0}$ is a principal quasi-polarization, the bilinear morphism $cu_{E_0}\circ (1_{E_0},\lambda_{E_0})$ over $S$ is symmetric. Thus $b_{N_0}$ is a symmetric form on $N_0$. We expect that (at least) under mild conditions the form $b_{N_0}$ on $N_0$ is in fact alternating. Here is an example over $k$ that implicitly points out that some conditions might be indeed required.
\medskip\noindent
{\bf 2.2.2. Example.} Let $E_1$ be the $2$-torsion group scheme of the unique supersingular $2$-divisible group over $k$ of height $2$. The evaluation $(N_1,\phi_1,\upsilon_1)$ of ${\bf D}(E_1)$ at the trivial thickening of ${\rm Spec}\,k$ is a $2$-dimensional $k$-vector space spanned by elements $e_1$ and $f_1$ that satisfy $\phi_1(e_1)=\upsilon_1(e_1)=f_1$ and $\phi_1(f_1)=\upsilon_1(f_1)=0$. Let $b_{N_1}$ be the perfect, non-alternating symmetric bilinear form on $N_1$ given by the rules: $b_{N_1}(f_1,e_1)=b_{N_1}(e_1,e_1)=b_{N_1}(e_1,f_1)=1$ and $b_{N_1}(f_1,f_1)=0$. For all $u$, $x\in N_1$ we have $b_{N_1}(\phi_1(x),u)=b_{N_1}(x,\upsilon_1(u))^2$ and thus also $b_{N_1}(\phi_1(x),\phi_1(x))=b_{N_1}(\upsilon_1(x),\upsilon_1(x))=0$. It is well known  that the perfect bilinear form $b_{N_1}$ induces an isomorphism $\lambda_{E_1}:E_1\tilde\to E_1^{\rm t}$ over $k$. As we are in characteristic $2$, the morphism $co_{E_1}:E_1\to \pmb{\mu}_{2,k}$ over $k$ is a homomorphism over $k$. But as $E_1$ has no toric part, $co_{E_1}$ is trivial. Thus $\lambda_{E_1}$ is a principal quasi-polarization.
\medskip\noindent
{\bf 2.2.3. Lemma.}  {\it If $V=W(k)$, then the perfect form $b_{N_0}$ on $N_0$ is alternating.} 
\medskip\noindent
{\it Proof:} 
We have $e=1$ and thus $U_e=k$. Let ${\bf D}(E_0)=(N_0,F_0,\phi_0,\phi_1,\nabla_0)$ be the Dieudonn\'e module of $E_0$ as used in [40, Construction 2.2]. Thus $F_0={\rm Ker}(\phi_0)$ and $\phi_1:F_0\to N_0$ is a $\sg$-linear map such that we have $N_0=\phi_0(N_0)+\phi_1(F_0)$. Let $x\in N_0$ and $u\in F_0$. We have $b_{N_0}(\phi_0(x),\phi_0(x))=2b_{N_0}(x,x)^2=0$, $b_{N_0}(\phi_1(u),\phi_1(u))=2b_{N_0}(u,u)^2=0$, and $b_{N_0}(\phi_0(x),\phi_1(u))=b_{N_0}(x,u)^2$. Let $v\in N_0$. We choose $x$ and $u$ such that we have an identity $v=\phi_0(x)+\phi_1(u)$. As $b_{N_0}$ is symmetric, we compute that $b_{N_0}(v,v)=b_{N_0}(\phi_0(x),\phi_0(x))+b_{N_0}(\phi_1(u),\phi_1(u))=0+0=0$. Thus $b_{N_0}$ is alternating.\endproof
\medskip
If $V\neq W(k)$, then the ring $V/4V$ has no Frobenius lift. Thus the proof of Lemma 2.2.3 can not be adapted if $V\neq W(k)$. Based on this, for $V\neq W(k)$ we will study pairs of the form $(N_0,b_{N_0})$ only in those cases that are related naturally to the geometric context of Theorem 1.4 (see Proposition 5.1 below).
\bigskip\smallskip\noindent
{\bigsll {\bf 3. Group schemes and involutions}}
\bigskip
Let $n\in{\bf N}^{\ast}$. Let ${\rm Spec}\,S$ be an affine scheme. We recall that a reductive group scheme $\Mk$ over $S$ is a smooth, affine group scheme over $S$ whose fibres are connected and have trivial unipotent radicals. Let $\Mk^{\rm ad}$  and $\Mk^{\rm der}$ be the adjoint and the derived (respectively) group schemes of $\Mk$, cf. [12, Vol. III, Exp. XXII, Def. 4.3.6 and Thm. 6.2.1]. Let ${\rm Lie}(\ast)$ be the {\it Lie algebra} of a smooth, closed subgroup scheme $\ast$ of $\Mk$. If $M$ is a free $S$-module of finite rank, let ${\bf GL}_M$ be the reductive group scheme over $S$ of linear automorphisms of $M$ and let $\Mt(M):=\oplus_{s,m\in{\bf N}} M^{\otimes s}\otimes_S M^{\vee\otimes m}$. Each $S$-linear isomorphism $i_S:M\tilde\to\tilde M$ of free $S$-modules of finite rank, extends naturally to an $S$-linear isomorphism (denoted also by) $i_S:\Mt(M)\tilde\to\Mt(\tilde M)$ and thus we will speak about $i_S$ taking some tensor of $\Mt(M)$ to some tensor of $\Mt(\tilde M)$. We identify ${\rm End}(M)=M\otimes_S M^\vee$. If $(M,\lambda_M)$ is a symplectic space over $S$, let ${\bf Sp}(M,\lambda_M):={\bf GSp}(M,\lambda_M)^{\rm der}$. We often use the same notation for two elements of some modules (like involutions, endomorphisms, bilinear forms, etc.) that are obtained one from another via extensions of scalars and restrictions. The reductive group schemes ${\bf GL}_{n,S}$, ${\bf Sp}_{2n,S}$, etc., are over $S$. 
\smallskip
Let ${\bf SO}_{2n,S}^+$, ${\bf O}_{2n,S}^+$, and ${\bf GSO}_{2n,S}^+$ be the split ${\bf SO}_{2n}$, ${\bf O}_{2n}$, and ${\bf GSO}_{2n}$ (respectively) group schemes over $S$. The constructions of ${\bf SO}_{2n,{\bf Z}_{(2)}}^+$ and ${\bf O}_{2n,{\bf Z}_{(2)}}^+$ (and implicitly of ${\bf SO}_{2n,S}^+$ and ${\bf O}_{2n,S}^+$ with $S$ a ${{\bf Z}_{(2)}}$-algebra) are recalled in Subsection 3.1. Let ${\bf GSO}_{2n,S}^+$ be the quotient of ${\bf SO}_{2n,S}^+\times_S {\bf G}_{m,S}$ by a $\pmb{\mu}_{2,S}$ subgroup scheme embedded diagonally. 
\smallskip
In Subsections 3.1 and 3.2 we present some general facts on ${\bf SO}_{2n}^+$ group schemes in mixed characteristic $(0,2)$; see [7, plates I and IV] for the weights used. In Subsections 3.3 and 3.4 we include complements on involutions and on non-alternating symmetric bilinear forms. In Subsection 3.5 we review some standard properties of the Shimura pair $(G,\Mx)$ we introduced in Subsection 1.3. \medskip\smallskip\noindent
{\bf 3.1. The $\Md_n$ group scheme.} We consider the quadratic form 
$${\got Q}_n(x):=x_1x_2+\cdots+x_{2n-1}x_{2n}\;\;\; {\rm defined}\;\;\;{\rm for}\;\;\; x=(x_1,\ldots,x_{2n})\in \Ml_n:={\bf Z}_{(2)}^{2n}.$$ 
For $\al\in{\bf Z}_{(2)}$ and $x\in\Ml_n$ we have ${\got Q}_n(\al x)=\al^2{\got Q}_n(x)$. Let $\Mq_n$ be the closed subgroup scheme of ${\bf GL}_{\Ml_n}$ that fixes ${\got Q}_n$. Let $\Md_n$ be the schematic closure of the identity component of $\Mq_{n,\bf Q}$ in $\Mq_n$. The fibres of $\Mq_n$ are smooth and have identity components that are split reductive groups (see [6, Subsect. 23.6] for the fibre of $\Mq_n$ over ${\bf F}_2$). We get that $\Md_n$ is a smooth, affine group scheme over ${\bf Z}_{(2)}$ whose generic fibre is connected and whose special fibre $\Md_{n,{\bf F}_2}$ has an identity component which is a split reductive group. If $\Md_n^0$ is the open subgroup scheme of $\Md_n$ which is the union of the identity components of the fibres of $\Md_n$, then $\Md_n^0$ is the complement in $\Md_n$ of a divisor and thus the homomorphism $\Md_n^0\to\Md_n$ is affine. Therefore $\Md_n^0$ is an affine, smooth group scheme whose fibres are reductive groups. We conclude that $\Md_n^0$ is a reductive group scheme. This implies that for each finite, discrete valuation ring extension $\flat$ of ${\bf Z}_2$, $\Md_n^0(\flat)$ is a maximal compact subgroup of $\Md_n^0(\flat[{1\over 2}])$ (see [38, Subsect. 3.2 and Subsubsect. 3.8.1] applied to $\Md_{n,\flat}(\flat)$) and thus there exists no $\flat$-valued point of $\Md_n$ which lifts a $\flat/2\flat$-valued point of $\Md_n$ that does not belong to $\Md_n^0$. Thus, as $\Md_n$ is smooth, we get that $\Md_n=\Md_n^0$ and therefore $\Md_{n,{\bf F}_2}$ is connected. Moreover $\Md_n$ is a reductive group scheme; it has a maximal split torus of rank $n$ and thus it is split. We conclude that ${\bf SO}_{2n,{\bf Z}_{(2)}}^+$ is isomorphic to $\Md_n$ and in this way we have constructed it up to isomorphism (cf. the uniqueness of a split reductive group scheme associated to a root datum; see [12, Vol. III, Exp. XXIII, Cor. 5.1]). We get that $\Md_1\tilde\to{\bf G}_{m,{\bf Z}_{(2)}}$, that $\Md_2$ is the quotient of a product of two ${\bf SL}_{2,{\bf Z}_{(2)}}$ group schemes by a $\pmb{\mu}_{2,{\bf Z}_{(2)}}$ subgroup scheme embedded diagonally, and that for $n\ge 3$ the group scheme $\Md_n$ is semisimple and $\Md_n^{\rm ad}$ is a split, absolutely simple, adjoint group scheme of $D_n$ Lie type.  
\smallskip
Next we study the rank $2n$ faithful representation 
$$\rho_n:\Md_n\hookrightarrow {\bf GL}_{\Ml_n}.$$
\noindent
If $n\ge 4$, then  $\rho_n$ is associated to the minimal weight $\varpi_1$. Moreover $\rho_3$ is associated to the weight $\varpi_2$ of the $A_3$ Lie type and $\rho_2$ is the tensor product of the standard rank $2$ representations of the mentioned two ${\bf SL}_{2,{\bf Z}_{(2)}}$ group schemes. Thus $\rho_n$ is isomorphic to its dual and, up to a ${\bf G}_m({\bf Z}_{(2)})$-multiple, there exists a unique perfect, symmetric bilinear form ${\got B}_n$ on $\Ml_n$ fixed by $\Md_n$ (the case $n=1$ is trivial). In fact we can define ${\got B}_n$ by
$${\got B}_n(u,v):={\got Q}_n(u+v)-{\got Q}_n(u)-{\got Q}_n(v)\;\;\;\; \forall u, v\in \Ml_n.$$ 
We can recover ${\got Q}_n$ from ${\got B}_n$ via the formula ${\got Q}_n(x)={{{\got B}_n(x,x)}\over 2}$, where $x\in \Ml_n$. This formula makes sense as $2$ is a non-zero divisor in ${\bf Z}_{(2)}$. We emphasize that we can not recover in a canonical way the reduction modulo $2{\bf Z}_{(2)}$ of ${\got Q}_n$ from the reduction modulo $2{\bf Z}_{(2)}$ of ${\got B}_n$. If $n\ge 2$, then the representation of a split maximal torus $\Mt_n$ of $\Md_n$ on $\Ml_n$ is a direct sum of rank $1$ subrepresentations and the $2n$ characters of $\Mt_n$ associated to these subrepresentations are distinct and moreover they are permuted transitively by the Weyl group of $\Md_n$ with respect to $\Mt_n$. This implies that for $n\ge 2$ the special fibre $\rho_{n,{\bf F}_2}$ of $\rho_n$ is an absolutely irreducible representation. Let $J(2n)$ be the matrix of ${\got B}_n$ with respect to the standard ${\bf Z}_{(2)}$-basis for $\Ml_n$; it is formed by $n$ diagonal blocks that are ${0\,1}\choose {1\,0}$. Therefore ${\got B}_n$ modulo $2{\bf Z}_{(2)}$ is alternating. 
\smallskip
As $\Mq_{n,{\bf F}_2}$ is smooth, we easily get that $\Md_n$ is the identity component of $\Mq_n$. We check that $\Mq_n$ is isomorphic to ${\bf O}^+_{2n,{\bf Z}_{(2)}}$ i.e., we have a non-trivial short exact sequence
$$1\to\Md_n\to\Mq_n\to {({\bf Z}/2{\bf Z})}_{{\bf Z}_{(2)}}\to 1\leqno (2)$$
that splits. We can  assume that $n\ge 2$ (as the case $n=1$ is easy). As $\rho_{n,{\bf F}_2}$ is absolutely irreducible for $n\ge 2$, the centralizer of  $\Md_{n,{\bf F}_2}$ in $\Mq_{n,{\bf F}_2}$ is a $\pmb{\mu}_{2,{\bf F}_2}$ group scheme and thus connected. Thus the group of connected components of $\Mq_{n,{\bf F}_2}$ acts on $\Md_{n,{\bf F}_2}$ via outer automorphisms that are as well automorphisms of $\rho_{n,{\bf F}_2}$. This implies that $\Mq_{n,{\bf F}_2}$ has at most two connected components (even if $2\le n\le 4$). But it is well known that $\Mq_{n,{\bf Q}}$ has two connected components and thus the open closed subscheme $\Mq_n\setminus\Md_n$ of $\Mq_n$ has a special fibre which is either empty or connected. The $2n\times 2n$ block matrix $\sharp:=$${J(2)\,O}\choose {O\,I_{2n-2}}$ is an element of $(\Mq_n\setminus\Md_n)({\bf Z}_{(2)})$ of order $2$ (as this is so over $\bf Q$); here $I_{2n-2}$ is the $(2n-2)\times (2n-2)$ identity matrix. From the last two sentences we get that $\Mq_n=\Md_n\sqcup \sharp\Md_n$ is a smooth group scheme over ${\bf Z}_{(2)}$ and that the short exact sequence (2) exists and splits (due to the existence of $\sharp$). The short exact sequence (2) is non-trivial as its generic fibre over $\bf Q$ is non-trivial. One could rephrase (2) by: $\Md_n$ is the kernel of the Dickson invariant epimorphism $\Mq_n\twoheadrightarrow {({\bf Z}/2{\bf Z})}_{{\bf Z}_{(2)}}$. We recall an obvious fact. 
\medskip\noindent
{\bf 3.1.1. Fact.} {\it Let $b_N$ be a symmetric bilinear form on a free module $N$ of finite rank over a commutative ${\bf F}_2$-algebra $S$. We have:
\medskip
{\bf (a)} The quadratic map $q_N:N\to S$ that maps $x\in N$ to $b_N(x,x)\in S$ is additive and its kernel ${\rm Ker}(q_N)$ is an $S$-submodule of $N$. We have $N={\rm Ker}(q_N)$ if and only if $b_N$ is alternating. 
\smallskip
{\bf (b)} If $S$ is a field, then $\dim_S(N/{\rm Ker}(q_N))\le [S:\Phi_S(S)]$.}
\medskip\noindent
{\bf 3.1.2. Fact.} {\it Let $\Mg_n$ be the flat, closed subgroup scheme of ${\bf GL}_{\Ml_n}$ generated by $\Md_n$ and by the center $\Mz_n$ of ${\bf GL}_{\Ml_n}$. Then $\Mg_n$ is a reductive group scheme isomorphic to ${\bf GSO}_{2n,{\bf Z}_{(2)}}^+$.}
\medskip\noindent
{\it Proof:} The kernel of the product homomorphism $\Pi_n:\Md_n\times_{{\bf  Z}_{(2)}} \Mz_n\to {\bf GL}_{\Ml_n}$ is isomorphic to ${\pmb{\mu}}_{2,{\bf  Z}_{(2)}}$. The quotient group scheme $\Md_n\times_{{\bf  Z}_{(2)}} \Mz_n/{\rm Ker}(\Pi_n)$ is reductive, cf. [12, Vol. III, Exp. XXII, Prop. 4.3.1]. The resulting homomorphism $\Md_n\times_{{\bf  Z}_{(2)}} \Mz_n/{\rm Ker}(\Pi_n)\to {\bf GL}_{\Ml_n}$ has fibres that are closed embeddings and thus (cf. [12, Vol. II, Exp. XVI, Cor. 1.5 a)]) it is a closed embedding. Thus we can identify $\Mg_n$ with $\Md_n\times_{{\bf  Z}_{(2)}} \Mz_n/{\rm Ker}(\Pi_n)$. Therefore $\Mg_n$ is a reductive group scheme isomorphic to ${\bf GSO}_{2n,{\bf Z}_{(2)}}^+$.\endproof
\medskip\noindent
{\bf 3.1.3. Lemma.} {\it Let $l$ be an algebraically closed field whose characteristic is either $0$ or $2$. Let $\Ml_n\otimes_{{\bf Z}_{(2)}} l=F_1\oplus F_0$ be a direct sum decomposition such that we have ${\got Q}_n(x)=0$ for all $x\in F_1\cup F_0$. Then the normalizer $\Mp_{n,l}$ of $F_1$ in $\Mg_{n,l}$ is a parabolic subgroup of $\Mg_{n,l}$ and we have $\dim(\Mg_{n,l}/\Mp_{n,l})={{n(n-1)}\over 2}$.}
\medskip\noindent
{\it Proof:} We can assume that $n\ge 2$. As we have ${\got Q}_n(x)=0$ for all $x\in F_1\cup F_0$, the cocharacter of ${\bf GL}_{\Ml_n\otimes_{{\bf Z}_{(2)}} l}$ that fixes $F_0$ and that acts as the inverse of the identity character of ${\bf G}_{m,l}$ on $F_1$, factors through $\Mg_{n,l}$. Let $\mu_l:{\bf G}_{m,l}\to\Mg_{n,l}$ be the resulting cocharacter. The fact that $\Mp_{n,l}$ is a parabolic subgroup of $\Mg_{n,l}$ is well known and it is a general consequence of the existence of $\mu_l$ (for instance, it is a particular case of [8, Prop. 2.1.8]). Let $T_l$ be a maximal torus of $\Mg_{n,l}$ through which $\mu_l$ factors; it is also a maximal torus of $\Mp_{n,l}$. Let ${\rm Lie}(\Mg_{n,l})={\rm Lie}(T_l)\bigoplus_{\al\in\Theta} {\got g}_{\al}$ be the root decomposition relative to $T_l$; thus $\Theta$ is a root system of $D_n$ Lie type. Let $\Theta_0$ be the subset of $\Theta$ formed by roots $\al$ with the property that $\mu_l$ acts on ${\got g}_{\al}$ either trivially or via the inverse of the identity character of ${\bf G}_{m,l}$. To check that $\dim(\Mg_{n,l}/\Mp_{n,l})={{n(n-1)}\over 2}$, it suffices to show that $\Theta\setminus\Theta_0$ has ${{n(n-1)}\over 2}$ roots. But it is well known that there exists a basis $\Delta=\{\alpha_1,\ldots,\alpha_n\}$ of $\Theta$ such that $\Theta\setminus\Theta_0$ is the set of roots of $\Theta$ that are linear combinations of elements of $-\Delta$ with non-negative integer coefficients, the coefficient of $\alpha_{n-1}$ being $-1$. There exist ${{n(n-1)}\over 2}$ such roots (cf. [7, Plate IV]) and thus $\Theta\setminus\Theta_0$ has precisely ${{n(n-1)}\over 2}$ roots.\endproof
\medskip\noindent
{\bf 3.1.4. Lemma.} {\it Let $\Ml_n\otimes_{{\bf Z}_{(2)}} {\bf Z_2}=L_0\oplus L_1$ be a direct sum decomposition such that we have ${\got Q}_n(x)=0$ for all $x\in L_0\cup L_1$. Let ${\got S}_n$ be the set of all ${\bf Z}_2$-lattices $\Ml$ of $\Ml_n\otimes_{{\bf Z}_{(2)}} {\bf Q_2}$ that have the following three properties:
\medskip
{\bf (i)} we have $L_1\subset\Ml\subset {1\over 2}L_1\oplus L_0$;
\smallskip
{\bf (ii)} it surjects onto $L_0$ producing a short exact sequence $0\to L_1\to\Ml\to L_0\to 0$ of ${\bf Z}_2$-modules;
\smallskip
{\bf (iii)} ${\got B}_n$ induces a perfect bilinear form on $\Ml$ whose reduction modulo $2{\bf Z}_2$ is alternating.
\medskip
Then the set ${\got S}_n$ has exactly $2^{{n(n-1)}\over 2}$ elements.} \medskip\noindent
{\it Proof:} Let $\Mp_n$ be the parabolic subgroup scheme of $\Mg_{n,{\bf Z}_2}$ that normalizes $L_1$, cf. Lemma 3.1.3. Let $\Mu_n$ be its unipotent radical: it is a smooth group scheme over ${\bf Z}_2$ that has connected fibres and that fixes both $L_1$ and $(\Ml_n\otimes_{{\bf Z}_{(2)}} {\bf Z_2})/L_1$. Therefore $\Mu_n$ is commutative and in fact it is isomorphic to the vector group scheme over ${\bf Z}_2$ defined by ${\rm Lie}(\Mu_n)$. The relative dimension of $\Mu_n$  is the relative dimension of $\Mg_{n,{\bf Z}_2}/\Mp_n$ and thus (cf. Lemma 3.1.3) it is ${{n(n-1)}\over 2}$. We conclude that $\Mu_n$ is isomorphic to ${\bf G}_{a,{\bf Z}_2}^{{{n(n-1)}\over 2}}$.
\smallskip
Let $\{u_1,\ldots,u_n\}$ be a ${\bf Z}_2$-basis for $L_1$. Let $\{v_1,\ldots,v_n\}$ be a ${\bf Z}_2$-basis for $L_0$ such that we have ${\got B}_n(u_i,v_j)=\delta_{ij}$. For $j\in\{1,\ldots,n\}$ let $\tilde v_j\in\Ml$ be such that it maps to $v_j$, cf. (ii). We have $w_j:=\tilde v_j-v_j\in {1\over 2}L_0$, cf. (i). By replacing each $\tilde v_j$ with $\tilde v_j+\tilde u_j$ for some $\tilde u_j\in L_1$ (cf. (i)), due to the property (iii) we can assume that ${\got B}_n(\tilde v_i,\tilde v_j)=0$ for all $i,j\in\{1,\ldots,n\}$. Let $(\tilde v_1',\ldots,\tilde v_n')$ be another $n$-tuple of $\Ml^n$ that maps to $(v_1,\ldots,v_n)$ and such that we have ${\got B}_n(\tilde v_i',\tilde v_j')=0$ for all $i,j\in\{1,\ldots,n\}$. Let $\tilde u_j:=\tilde v_j'-\tilde v_j\in L_1$ and let $g\in {\pmb GL}_{\Ml_n\otimes_{{\bf Z}_{(2)}} {\bf Z_2}}({\bf Z}_2)$ be such that it fixes each $u_i$ and maps each $v_j$ to $v_j+\tilde u_j$. As ${\got B}_n(\tilde v_i',\tilde v_j')=0={\got B}_n(\tilde v_i,\tilde v_j)$, we easily get that ${\got B}_n(v_i,\tilde u_j)+{\got B}_n(v_j,\tilde u_i)=0$ and thus that ${\got B}_n(g(v_i),g(v_j))={\got B}_n(v_i+\tilde u_i,v_j+\tilde u_j)=0$. As we also have ${\got B}_n(g(u_i),g(u_j))={\got B}_n(u_i,u_j)=0$ and ${\got B}_n(g(u_i),g(v_j))={\got B}_n(u_i,v_j+\tilde u_j)=\delta_{ij}$, we conclude that $g$ fixes ${\got B}_n$. As the determinant of $g$ is $1$, we have $g\in\Md_n({\bf Z}_2)$. From this and the very definition of $g$ we get that $g\in\Mu_n({\bf Z}_2)$.
\smallskip
Similarly we argue that the element $\tilde g\in {\pmb GL}_{\Ml_n\otimes_{{\bf Z}_{(2)}} {\bf Z_2}}({\bf Q}_2)$ that fixes each $u_i$ and that maps each $v_j$ to $\tilde v_j$ is in fact an element of $\Mu_n({\bf Q}_2)$. The replacement of $(\tilde v_1,\ldots,\tilde v_n)$ by $(\tilde v_1',\ldots,\tilde v_n')$ induces the replacement of $\tilde g$ by $g\tilde g=\tilde g g$. Thus the class $[\tilde g]\in \Mu_n({\bf Q}_2)/\Mu_n({\bf Z}_2)$ is well defined. As each $w_j\in {1\over 2}L_0$, in fact we have $[\tilde g]\in \Gamma_n$, where $\Gamma_n:=[\Mu_n({\bf Q}_2)/\Mu_n({\bf Z}_2)][2]$ is the $2$-torsion subgroup. As $\Mu_n$ is isomorphic to ${\bf G}_{a,{\bf Z}_2}^{{{n(n-1)}\over 2}}$, $\Gamma_n$ is an elementary abelian $2$-group of order $2^{{n(n-1)}\over 2}$.  
\smallskip
We have $\Ml=\tilde g(\Ml_n\otimes_{{\bf Z}_{(2)}} {\bf Z_2})$ and thus $\Ml$ is uniquely determined by $[\tilde g]\in \Gamma_n$. Thus the function ${\got S}_n\to\Gamma_n$ that maps $\Ml$ to $[\tilde g]$ is injective. As this function is obviously surjective, we conclude that it is a bijection. Therefore the set ${\got S}_n$ has $2^{{n(n-1)}\over 2}$ elements.\endproof
 \medskip\noindent
{\bf 3.2. Lemma.} {\it Let $S$ be a commutative, faithfully flat ${\bf Z}_{(2)}$-algebra. Let $\Ml$ be a free $S$-submodule of $\Ml_n\otimes_{{\bf Z}_{(2)}} S[{1\over 2}]$ of rank $2n$ such that $\Ml[{1\over 2}]=\Ml_n\otimes_{{\bf Z}_{(2)}} S[{1\over 2}]$ and we get a perfect bilinear form ${\got B}_n:\Ml\otimes_S \Ml\to S$. Suppose that there exists a reductive, closed subgroup scheme $\Md_n(\Ml)$ of ${\bf GL}_{\Ml}$ whose fibre over $S[{1\over 2}]$ is $\Md_{n,S[{1\over 2}]}$ (it is the schematic closure of $\Md_{n,S[{1\over 2}]}$ in ${\bf GL}_{\Ml}$). Then the perfect bilinear form on $\Ml/2\Ml$ induced by ${\got B}_n$ is alternating.}
\medskip\noindent
{\it Proof:} The statement of the lemma is local in the flat topology of ${\rm Spec}\,S$. Thus we can  assume that $S$ is local and $\Md_n(\Ml)$ is split, cf. [12, Vol. III, Exp. XIX, Prop. 6.1]. Let $T$ be a maximal split torus of $\Md_n(\Ml)$. We consider the $T$-invariant direct sum decomposition $\Ml=\oplus_{i=1}^{2n} \Mv_i$ into free $S$-modules of rank 1. Let $\Mw$ be an $S$-basis for $\Ml$ formed by elements of $\Mv_i$'s. As $T$-acts on each $\Mv_i$ via a non-trivial character and as $T$ fixes ${\got B}_n$, for each $x\in\Mw$ we have ${\got B}_n(x,x)=0$. Based on this and on the fact that ${\got B}_n$ is symmetric, from Fact 3.1.1 (a) we get that the perfect bilinear form on $\Ml/2\Ml$ induced by ${\got B}_n$ is alternating.\endproof
\medskip\smallskip\noindent
{\bf 3.3. Involutions.} See [25, Ch. I] for standard terminology and properties of involutions of (finite dimensional) semisimple algebras over fields. Let ${\rm Spec}\,S$ be a connected, reduced, affine scheme. Let $\Mf=\oplus_{i=1}^s \Mf_i$ be a product of matrix $S$-algebras. Let $I_{\Mf}$ be an involution of $\Mf$ that fixes $S$. Thus $I_{\Mf}$ is an $S$-linear automorphism of $\Mf$ of order $2$ such that for all $x$, $u\in \Mf$ we have $I_{\Mf}(xu)=I_{\Mf}(u)I_{\Mf}(x)$. 
\smallskip
Let $\eta$ be the permutation of $\{1,\ldots,s\}$ such that we have $I_{\Mf}(\Mf_i)=\Mf_{\eta(i)}$ for all $i\in\{1,\ldots,s\}$. The order of $\eta$ is either $1$ or $2$. Let $C_0$ be a subset of $\{1,\ldots,s\}$ whose elements are permuted transitively by $\eta$. Let $\Mf_0:=\oplus_{i\in C_0} \Mf_i$. Let $I_0$ be the restriction of $\Mi_{\Mf}$ to $\Mf_0$. We refer to $(\Mf_0,I_0)$ as a simple factor of $(\Mf,I_{\Mf})$. We say $(\Mf_0,I_0)$ is: 
\medskip
--  of second type, if $C_0$ has two elements; 

-- of first type, if $C_0$ has only one element and $I_0\neq 1_{\Mf_0}$;

-- trivial, if $C_0$ has only one element and $I_0=1_{\Mf_0}$.  
\medskip\noindent
{\bf 3.3.1. Lemma.} {\it Suppose that $S$ is also local and $(\Mf_0,I_0)$ is of first type (thus $\Mf_0$ is $\Mf_i$ for some $i\in\{1,\ldots,s\}$). We identify $\Mf_0={\rm End}(M_0)$, where $M_0$ is a free $S$-module of finite rank. We have:
\medskip
{\bf (a)} There exists a perfect bilinear form $b_0$ on $M_0$ such that the involution of $\Mf_0$ is defined by $b_0$ i.e., for all $u$, $v\in M_0$ and all $x\in \Mf_0$ we have $b_0(x(u),v)=b_0(u,I_0(x)(v))$. 
\smallskip
{\bf (b)}
The form $b_0$ (of (a)) is uniquely determined up to a ${\bf G}_m(S)$-multiple. 
\smallskip
{\bf (c)} If $S$ is also integral, then $b_0$ is either alternating or non-alternating symmetric.} 
\medskip\noindent
{\it Proof:} Let ${\bf Aut}({\bf GL}_{M_0})$ be the group scheme of automorphisms of ${\bf GL}_{M_0}$. To prove (a) we follow [25, Ch. I]. We identify the opposite $S$-algebra of $\Mf_0$ with ${\rm End}(M_0^\vee)$. Therefore $I_0$ defines naturally an $S$-isomorphism $I_0:{\rm End}(M_0)\tilde\to {\rm End}(M_0^\vee)$. We claim that each such $S$-isomorphism is defined naturally by an $S$-linear isomorphism $c_0:M_0\tilde\to M_0^\vee$. To check this claim, we fix an $S$-isomorphism $I_1:{\rm End}(M_0^\vee)\tilde\to {\rm End}(M_0)$ defined by an $S$-linear isomorphism $c_1:M_0^\vee\tilde\to M_0$ and it suffices to show that the $S$-automorphism $a_0:=I_1\circ I_0$ of ${\rm End}(M_0)$ is defined (via inner conjugation) by an element $d_0\in {\bf GL}_{M_0}(S)$. The element $a_0$ defines naturally an $S$-valued automorphism $\tilde a_0\in {\bf Aut}({\bf GL}_{M_0})(S)$ of ${\bf GL}_{M_0}$. 
\smallskip
Let ${\bf PGL}_{M_0}:={\bf GL}_{M_0}^{\rm ad}$. We have a short exact sequence $1\to {\bf PGL}_{M_0}\to {\bf Aut}({\bf GL}_{M_0})\to E_0\to 1$, where the group scheme $E_0$ over $S$ is either trivial or ${{\bf Z}/2{\bf Z}}_S$ (cf. [12, Vol. III, Exp. XXIV, Thm. 1.3]). We check that $\tilde a_0\in {\bf PGL}_{M_0}(S)$. We need to show that the image $\tilde a_2$ of $\tilde a_0$ in $E_0(S)$ is the identity element. As $S$ is reduced, to check this last thing we can assume that $S$ is a field and therefore the fact that $\tilde a_2$ is the identity element is implied by the Skolem--Noether theorem. Thus $\tilde a_0\in {\bf PGL}_{M_0}(S)$. 
\smallskip
As $S$ is local, all torsors of ${\bf G}_{m,S}$ are trivial. Thus there exists $d_0\in {\bf GL}_{M_0}(S)$ that maps to $\tilde a_0$. Therefore $c_0:=c_1^{-1}\circ d_0$ exists. Let $b_0$ be the perfect bilinear form on $M_0$ defined naturally by $c_0$. As $S$ is reduced, as in the case of fields one checks that for all $u$, $v\in M_0$ and all $x\in \Mf_0$ we have $b_0(x(u),v)=b_0(u,I_0(x)(v))$.  This proves (a). 
\smallskip
We check (b). Both $b_0$ and $c_0$ are uniquely determined by $d_0$. As $a_0$ and $\tilde a_0$ are uniquely determined by $I_0$, we get that $d_0$ is uniquely determined by $I_0$ up to a scalar multiplication with an element of ${\bf G}_m(S)$. From this (b) follows. For (c) we can  assume that $S$ is a field and this case is well known (for instance, see [25, Ch. I, Subsect. 2.1]).\endproof
\medskip\noindent
{\bf 3.3.2. Definition.} Suppose that $S$ is local and integral. Let $b_0$ be as in Lemma 3.3.1 (a). Let $J$ be an ideal of $S$. If the reduction of $b_0$ modulo $J$ is alternating (resp. is non-alternating symmetric), then we say that the reduction of the simple factor $(\Mf_0,I_0)$ modulo $J$ is of alternating (resp. of orthogonal) first type.
\medskip
We have the following converse form of Lemma 3.2.
\medskip\smallskip\noindent
{\bf 3.4. Proposition.} {\it Suppose that the commutative ${\bf Z}_{(2)}$-algebra $S$ is local, noetherian, $2$-adically complete, and strictly henselian. Let $b_M$ be a perfect symmetric bilinear form on a free $S$-module $M$ of finite rank $2n$ with the property that $b_M$ modulo $2S$ is alternating. Then the following three properties hold.
\medskip
{\bf (a)} There exists an $S$-basis $\Mw$ for $M$ with respect to which the matrix of $b_M$ is $J(2n)$. 
\smallskip
{\bf (b)} Let $q\in {\bf N}^{\ast}\setminus \{1\}$ and let $\Mw_q=\{u_{1,q},v_{1,q},\ldots,u_{n,q},v_{n,q}\}$ be an $S/2^qS$-basis for $M/2^qM$ such that for all $i,j\in\{1,\ldots,n\}$ we have $b_M(u_{i,q},u_{j,q})=0$, $b_M(u_{i,q},v_{j,q})=\delta_{ij}$, and $b_M(v_{i,q},v_{j,q})\in 2^{q-1}S/2^qS$. Then there exists an $S$-basis $\{u_1,v_1,\ldots,u_n,v_n\}$ for $M$ with respect to which the matrix of $b_M$ is $J(2n)$ and which has the property that for each $i\in\{1,\ldots,n\}$ the reduction of $u_i$ modulo $2^{q-1}S$ is $u_{i,q}$ modulo $2^{q-1}S/2^qS$ and the reduction of $v_i$ modulo $2^{q-2}S$ is $v_{i,q}$ modulo $2^{q-2}S/2^qS$.
\smallskip
{\bf (c)} Suppose $S$ is also a faithfully flat ${\bf Z}_{(2)}$-algebra. Let $q_M:M\to S$ be the quadratic form defined by $q_M(x):={{b_M(x,x)}\over 2}$, where $x\in M$. Let ${\bf O}(M,q_M)$ be the closed subgroup scheme of ${\bf GL}_M$ that fixes $q_M$. Let  ${\bf SO}(M,q_M)$ be the identity component of ${\bf O}(M,q_M)$. Let ${\bf GSO}(M,q_M)$ be the closed subgroup scheme of ${\bf GL}_M$ generated by ${\bf SO}(M,q_M)$ and by the center of ${\bf GL}_M$. Then ${\bf O}(M,q_M)$ (resp. ${\bf SO}(M,q_M)$ or ${\bf GSO}(M,q_M)$) is isomorphic to $\Mq_{n,S}$ (resp. $\Md_{n,S}$ or $\Mg_{n,S}$). If $S$ is also reduced and if $K_S$ denotes its ring of fractions, then ${\bf SO}(M,q_M)$ is the Zariski closure of ${\bf SO}(M,q_M)_{K_S}$ in ${\bf GL}_M$.}
\medskip\noindent
{\it Proof:} We prove (a). Our hypotheses imply that $S/2S$ is local, noetherian, and strictly henselian. As $b_M$ modulo $2S$ is alternating, there exists an $S/2S$-basis for $M/2M$ with respect to which the matrix of $b_M$ modulo $2S$ is $J(2n)$. By induction on $q\in{\bf N}^{\ast}$ we show that there exists an $S/2^qS$-basis for $M/2^qM$ with respect to which the matrix of $b_M$ modulo $2^{q}S$ is $J(2n)$. The passage from $q$ to $q+1$ goes as follows. We recall that we denote also by $b_M$ its reductions modulo $2^qS$ or $2^{q+1}S$.  Let $\Mw_q=\{u_{1,q},v_{1,q},\ldots,u_{n,q},v_{n,q}\}$ be an $S/2^qS$-basis for $M/2^qM$ such that we have $b_M(u_{i,q},v_{j,q})=\delta_{ij}$ and $b_M(u_{i,q},u_{j,q})=b_M(v_{i,q},v_{j,q})=0$ for all $i$, $j\in \{1,\ldots,n\}$. 
\smallskip
Let $\Mw_{q+1}=\{u_{1,q+1},v_{1,q+1},\ldots,u_{n,q+1},v_{n,q+1}\}$ be an $S/2^{q+1}S$-basis for $M/2^{q+1}M$ which modulo $2^{q-1}S/2^{q+1}S$ is $\Mw_q$ modulo $2^{q-1}S/2^qS$, which modulo $2^qS/2^{q+1}S$ is an $S/2^qS$-basis for $M/2^qM$ with respect to which the matrix of $b_M$ modulo $2^qS$ is $J(2n)$, and which is such that for all $i$, $j\in\{1,\ldots,n\}$ with $i\neq j$ we have $b_M(u_{i,q+1},v_{j,q+1})=b_M(u_{i,q+1},u_{j,q+1})=b_M(v_{i,q+1},v_{j,q+1})=0$. One constructs $\Mw_{q+1}$ by first considering an arbitrary $S/2^{q+1}S$-basis $\tilde{\Mw}_{q+1}=\{\tilde u_{1,q+1},\tilde v_{1,q+1},\ldots,\tilde u_{n,q+1},\tilde v_{n,q+1}\}$ for $M/2^{q+1}M$ which modulo $2^qS/2^{q+1}S$ is $\Mw_q$ and then by modifying inductively for $l\in\{2,\ldots,n\}$ the pair $(\tilde u_{l,q+1},\tilde v_{l,q+1})$ so that modulo $2^qS/2^{q+1}S$ it is not changed and its $S/2^{q+1}S$-span is perpendicular on the  $S/2^{q+1}S$-span of $\{\tilde u_{1,q+1},\tilde v_{1,q+1},\ldots,\tilde u_{l-1,q+1},\tilde v_{l-1,q+1}\}$ with respect to $b_M$.
\smallskip
Let $u$, $v\in M/2^{q+1}M$ be such that $b_M(u,u)=2^qa$, $b_M(v,v)=2^qc$, and $b_M(u,v)=1+2^qe$, where $a$, $c$, $e\in S/2^{q+1}S$. By replacing $v$ with $(1+2^qe)v$, we can  assume that $e=0$. We show that there exists $x\in S/2^{q+1}S$ such that 
$$b_M(u+2^{q-1}xv,u+2^{q-1}xv)=2^q(a+x+2^{2q-2}x^2c)$$ 
is $0$. If $q\ge 2$, then we can take $x=-a$. If $q=1$, then as $x$ we can take any element of $S/2^{q+1}S$ that modulo $2S/2^{q+1}S$ is a solution of the equation in $t$
$$\bar a+t+\bar ct^2=0,\leqno (3)$$ 
where $\bar a$ and $\bar c\in S/2S$ are the reductions modulo $2S/2^{q+1}S$ of $a$ and $c$ (respectively). The $S/2S$-scheme ${\rm Spec}\,S/2S[t]/(\bar a+t+\bar ct^2)$ is \'etale and has points over the maximal point of ${\rm Spec}\,S/2S$. As $S/2S$ is strictly henselian, the equation (3) has solutions in $S/2S$. Thus $x$ exists even if $q=1$. Therefore by replacing $(u,v)$ with $(u+2^{q-1}xv,(1+2^{2q-1}xc)v)$, we can also assume that $a=0$. By replacing $(u,v)$ with $(u,v+2^{q-1}c u)$ we can also assume that $c=0$. Thus $b_M(u,u)=b_M(v,v)=0$ and $b_M(u,v)=1$. Under all these replacements, the $S/2^{q+1}S$-span of $\{u,v\}$ does not change. 
\smallskip
Applying the previous paragraph to all pairs $(u,v)$ in $\{(u_{i,q+1},v_{i,q+1})|1\le i\le n\}$, we get that we can choose $\Mw_{q+1}$ such that the matrix of $b_M$ modulo $2^{q+1}S$ with respect to $\Mw_{q+1}$ is $J(2n)$. This completes the induction.
\smallskip
As $S$ is $2$-adically complete and as $\Mw_{q+1}$ modulo $2^{q-1}S/2^{q+1}S$ is $\Mw_q$ modulo $2^{q-1}S/2^{q}S$, there exists an $S$-basis $\Mw$ for $M$ that modulo $2^qS$ coincides with $\Mw_{q+1}$ modulo $2^qS/2^{q+1}S$, for all $q\in{\bf N}^{\ast}$. The matrix of $b_M$ with respect to $\Mw$ is $J(2n)$. Thus (a) holds. Part (b) follows from the proof of (a). 
\smallskip
The triple $(S,M,q_M)$ is isomorphic to the extension to $S$ of the triple $({\bf Z}_{(2)},\Ml_n,{\got Q}_n)$ of Subsection 3.1, cf. (a). Thus (b) follows from the definitions of $\Md_n$, $\Mq_n$, and $\Mg_n$ (see Subsection 3.1 and Fact 3.1.2). \endproof
\medskip\smallskip\noindent
{\bf 3.5. A review.} In this subsection we use the notations of Subsection 1.3 and review standard properties of the Shimura pair $(G,\Mx)$. Let $G^0:=G\cap {\bf Sp}(W,\psi)$; it is the normal subgroup of $G$ that fixes $\psi$. Let $h:{\bf S}\hookrightarrow G_{{\bf R}}$ be an element of $\Mx$. If $W\otimes_{{\bf Q}} {\bf C}=F_h^{-1,0}\oplus F^{0,-1}_h$ is the Hodge decomposition defined by $h$, let $\mu_h:{\bf G}_{m,{\bf C}}\to G_{{\bf C}}$ be the Hodge cocharacter that fixes $F^{0,-1}_h$ and that acts as the identity character of ${\bf G}_{m,{\bf C}}$ on $F_h^{-1,0}$. The cocharacter $\mu_h$ acts on the ${\bf C}$-span of $\psi$ via the identity character of ${\bf G}_{m,{\bf C}}$. The image through $h$ of the compact subtorus of ${\bf S}$ contains the center of ${\bf Sp}(W,\psi)_{\bf R}$ and it is contained in the extension to ${\bf R}$ of the identity component $G^{00}$ of $G^0$. Thus $G^{00}$ contains the center of ${\bf Sp}(W,\psi)$ and it is easy to see that this implies that the epimorphism $G/G^{00}\twoheadrightarrow G/G^0$ is in fact an isomorphism.  Thus $G^0=G^{00}$ is connected and therefore it is a reductive group. Let $G^0_{{\bf Z}_{(2)}}$ be the schematic closure of $G^0$ in $G_{{\bf Z}_{(2)}}$. 
\medskip\noindent
{\bf 3.5.1. Some group schemes.} Let $\Mb_1$ be the centralizer of $\Mb$ in ${\rm End}(L_{(2)})$. Let $G_{2,{\bf Z}_{(2)}}$ be the centralizer of $\Mb$ in ${\bf GL}_{L_{(2)}}$; it is the group scheme over ${\bf Z}_{(2)}$ of invertible elements of $\Mb_1$. Due to the property 1.3 (i), the $W({\bf F})$-algebra $\Mb_1\otimes_{{\bf Z}_{(2)}} W({\bf F})$ is also a product of matrix $W({\bf F})$-algebras. Therefore $G_{2,W({\bf F})}$ is a product of ${\bf GL}$ groups schemes over $W({\bf F})$ and thus $G_{2,{\bf Z}_{(2)}}$ is a reductive group scheme. As $\Mi(\Mb)=\Mb$ we have also $\Mi(\Mb_1)=\Mb_1$. The involution $\Mi$ of $\Mb_1$ defines an involution (denoted also by $\Mi$) of the group of points of $G_{2,{\bf Z}_{(2)}}$ with values in each fixed ${\bf Z}_{(2)}$-scheme. Let $G_1^0:=G_1\cap {\bf Sp}(W,\psi)$. Let $G^0_{1,{\bf Z}_{(2)}}$ be the schematic closure of $G_1^0$ in ${\bf Sp}(L_{(2)},\psi)$; it is the maximal flat, closed subgroup scheme of $G_{2,{\bf Z}_{(2)}}$ with the property that for each commutative ${\bf Z}_{(2)}$-algebra $\sharp$ and for every element $\dag\in G^0_{1,{\bf Z}_{(2)}}(\sharp)$ we have $\Mi(\dag)=\dag^{-1}$.  
\smallskip
Let $(\Mb_1\otimes_{{\bf Z}_{(2)}} {\bf Z}_2,\Mi)=\oplus_{j\in \kappa} (\Mb_j,\Mi)$ be the product decomposition of $(\Mb_1\otimes_{{\bf Z}_{(2)}} {\bf Z}_2,\Mi)$ into simple factors. Each $\Mb_j$ is a two sided ideal of $\Mb_1\otimes_{{\bf Z}_{(2)}} {\bf Z}_2$ that is a product of one or two simple ${\bf Z}_2$-algebras. It is easy to see that $(\Mb_1\otimes_{{\bf Z}_{(2)}} W({\bf F}),\Mi)$ has no trivial factor and that all simple factors of $(\Mb_j,\Mi)\otimes_{{\bf Z}_2} W({\bf F})$ have the same type. We have a product decomposition $G_{2,{\bf Z}_2}=\prod_{j\in \kappa} G_{2,j}$, where $G_{2,j}$ is defined by the property that ${\rm Lie}(G_{2,j})$ is the Lie algebra associated to $\Mb_j$. We have also product decompositions $G^0_{1,{\bf Z}_2}=\prod_{j\in \kappa} G^0_{1,j}$ and $G^0_{{\bf Z}_2}=\prod_{j\in \kappa} G^0_j$, where $G^0_{1,j}:=G^0_{1,{\bf Z}_2}\cap G_{2,j}$ and $G^0_{j}:=G^0_{{\bf Z}_2}\cap G_{2,j}$. 
\smallskip 
The double monomorphism $G^0_{j,W({\bf F})}\hookrightarrow G^0_{1,j,W({\bf F})}\hookrightarrow G_{2,j,W({\bf F})}$ over $W({\bf F})$ is a product of double monomorphisms over $W({\bf F})$ that are of one of the following three types (we recall that $G$ is not a torus): 
\medskip
-- ${\bf GL}_{n,W({\bf F})}={\bf GL}_{n,W({\bf F})}\hookrightarrow {\bf GL}_{n,W({\bf F})}\times {\bf GL}_{n,W({\bf F})}$ (with $n\ge 2$), if each simple factor of $(\Mb_j,\Mi)\otimes_{{\bf Z}_2} W({\bf F})$ is of second type; 

-- ${\bf Sp}_{2n,W({\bf F})}={\bf Sp}_{2n,W({\bf F})}\hookrightarrow {\bf GL}_{2n,W({\bf F})}$ (with $n\in{\bf N}^{\ast}$), if each simple factor of $(\Mb_j,\Mi)\otimes_{{\bf Z}_2} W({\bf F})$ is of symplectic first type;

--  ${\bf SO}_{2n,W({\bf F})}\hookrightarrow {\bf O}_{2n,W({\bf F})}\hookrightarrow {\bf GL}_{2n,W({\bf F})}$ (with $n\ge 2$), if each simple factor of $(\Mb_j,\Mi)\otimes_{{\bf Z}_2} W({\bf F})$ is of orthogonal first type. [The case $n=1$ is excluded as ${\bf SO}_{2,W({\bf F})}$ is a torus and as $\Mb$ is the maximal ${\bf Z}_{(2)}$-subalgebra of ${\rm End}(L_{(2)})$ fixed by $G_{{\bf Z}_{(2)}}$.]

\medskip
Thus, as we are in the case (D), each simple factor of $(\Mb_j,\Mi)\otimes_{{\bf Z}_2} W({\bf F})$ is of orthogonal first type and the quotient group scheme $G^0_{1,j,W({\bf F})}/G^0_{j,W({\bf F})}$ is a product of one or more $({{\bf Z}/2{\bf Z}})_{W({\bf F})}$ group schemes. 
\smallskip
We get that the quotient group $G_1/G$ is a finite, non-trivial $2$-torsion group. We get also that $G^0_{{\bf C}}$ is isomorphic to a product of ${\bf SO}_{2n,{\bf C}}$ groups ($n\ge 2$) and thus that $G^{\rm der}_{{\bf Z}_{(2)}}=G^0_{{\bf Z}_{(2)}}$. From this and [37, Subsects. 2.6 and 2.7] we get that $G^0_{{\bf R}}$ is a product of ${\bf SO}^*_{2n}$ groups. We recall that ${\bf SO}^*_{2n}$ is the semisimple group over ${\bf R}$ whose ${\bf R}$-valued points are those elements of ${\bf SL}_{2n}({\bf C})$ that fix both the quadratic form $z_1^2+\cdots+z_{2n}^2$ and the skew hermitian form $-z_1\bar z_{n+1}+z_{n+1}\bar z_1-\cdots-z_n\bar z_{2n}+z_{2n}\bar z_n$.
\smallskip
The $G({\bf R})$-conjugacy class $\Mx$ of $h$ is a disjoint union of connected hermitian symmetric domains (cf. [10, Cor. 1.1.17]) which (cf. the structure of $G_{{\bf R}}^0$) are products of isomorphic irreducible hermitian symmetric domains of D III type i.e., of the form ${\bf SO}^*_{2n}({\bf R})/{\bf U}_n({\bf R})$ with $n\ge 2$. But ${\bf SO}^*_{2n}({\bf R})/{\bf U}_n({\bf R})$ has real dimension $n(n-1)$, cf. [21, Ch. X, Sect. 6, Table V]. If $G_{{\bf Z}_2}$ is isomorphic to ${\bf GSO}_{2n,{\bf Z}_2}^+$, then each connected component of $\Mx$ is isomorphic to ${\bf SO}^*_{2n}({\bf R})/{\bf U}_n({\bf R})$ and has (complex) dimension ${{n(n-1)}\over 2}$. 
\medskip\noindent
{\bf 3.5.2. Extra tensors.} 
For $g\in G({\bf A}_f^{(2)})\subset G({\bf A}_f)$, let $L_g$ be the ${\bf Z}$-lattice of $W$ such that we have $L_g\otimes_{{\bf Z}} \widehat{{\bf Z}}=g(L\otimes_{{\bf Z}} \widehat{{\bf Z}})$. We have $L_{(2)}=L_g\otimes_{{\bf Z}} {\bf Z}_{(2)}$. Let $(v_\al)_{\al\in\Mj}$ be a family of tensors of $\Mt(W)$ such that $G$ is the subgroup of ${\bf GL}_W$ that fixes $v_{\al}$ for all $\al\in\Mj$, cf. [11, Prop. 3.1 (c)]. We choose the set $\Mj$ such that $\Mb\subseteq\Mj$ and for $b\in\Mb$ we have $v_b=b\in {\rm End}(W)=W\otimes_{{\bf Q}} W^\vee$. We recall that if $h\in\Mx$ is as in the beginning of Subsection 3.5 and if 
$$w=[h,g_w]\in G_{{\bf Z}_{(2)}}({\bf Z}_{(2)})\backslash (\Mx\times G({\bf A}_f^{(2)}))=({\rm Sh}(G,\Mx)/H_2)({\bf C})=\Mn({\bf C}),$$ 
then the complex Lie group associated to $w^*(\Ma)$ is $L_{g_w}\backslash L\otimes_{{\bf Z}} {\bf C}/F^{0,-1}_h$, $v_b$ is the Betti realization of the ${\bf Z}_{(2)}$-endomorphism of $w^*(\Ma)$ defined by $b\in\Mb$, and the non-degenerate alternating form on $L_{(2)}=H_1(w^*(\Ma),{\bf Z}_{(2)})$ defined by $w^*(\Lambda)$ is a ${\bf G}_m({\bf Z}_{(2)})$-multiple of $\psi$ (see [31, Ch. 3] and [39, Subsect. 4.1]). We get also a natural identification 
$$L_{(2)}\otimes_{{\bf Z}_{(2)}} {\bf Z}_2=(H^1_{\acute et}(w^*(\Ma),{\bf Z}_2))^\vee$$ 
under which $v_b$ and $\psi$ get identified with the $2$-adic \'etale realization of the ${\bf Z}_{(2)}$-endomorphism of $w^*(\Ma)$ defined by $b\in\Mb$ and respectively with a ${\bf G}_m({\bf Z}_{2})$-multiple of the perfect form on $(H^1_{\acute et}(w^*(\Ma),{\bf Z}_2))^\vee$ defined by $w^*(\Lambda)$. 
\medskip\noindent
{\bf 3.5.3. On the case of ${\bf GSO}_{2n,{\bf Z}_2}^+$.} In this subsubsection we  assume $G_{{\bf Z}_2}$ is isomorphic to ${\bf GSO}_{2n,{\bf Z}_2}^+$. From Subsubsection 3.5.1 we get that $n\ge 2$ and that $G_{2,{\bf Z}_2}$ is isomorphic to ${\bf GL}_{2n,{\bf Z}_2}$. Thus we can identify $\Mb_1\otimes_{{\bf Z}_{(2)}} {\bf Z}_2$ with ${\rm End}(\Mv)$, where $\Mv$ is a free ${\bf Z}_2$-module of rank $2n$. Let $s\in{\bf N}^{\ast}\setminus\{1\}$ be such that as $\Mb_1\otimes_{{\bf Z}_{(2)}} {\bf Z}_2$-modules we can identify
$$L_{(2)}\otimes_{{\bf Z}_{(2)}} {\bf Z}_2=\Mv^s.\leqno (4)$$
Let $b_{\Mv}$ be a perfect bilinear form on $\Mv$ that defines the involution $\Mi$ of $\Mb_1\otimes_{{\bf Z}_{(2)}} {\bf Z}_2$, cf. Lemma 3.3.1 (a). It is unique up to a ${\bf G}_m({\bf Z}_2)$-multiple (cf. Lemma 3.3.1 (b)), it is fixed by $G^{\rm der}_{{\bf Z}_2}=G^0_{{\bf Z}_2}$, and it is symmetric (as $(\Mb_1\otimes_{{\bf Z}_{(2)}} {\bf Z}_2,\Mi)$ is of orthogonal first type). Moreover the abelianization $G_{{\bf Z}_2}/G^{\rm der}_{{\bf Z}_2}=G_{{\bf Z}_2}/G^0_{{\bf Z}_2}$ of $G_{{\bf Z}_2}$ is a split torus of rank $1$ which acts on the ${\bf Z}_2$-spans of $b_{\Mv}$ and $\psi$ via the same character.  
\medskip\noindent
{\bf 3.5.4. Lemma.} {\it Each simple factor of $(\Mb_1\otimes_{{\bf Z}_{(2)}} {\bf F},\Mi)$ is of alternating first type.} 
\medskip\noindent 
{\it Proof:} Let $j\in\kappa$ and let $({\rm End}(W({\bf F})^{2n}),\Mi)$ be a simple factor of $(\Mb_j\otimes_{{\bf Z}_{2}} W({\bf F}),\Mi)$; it is of orthogonal first type (see Subsubsection 3.5.1). Let ${\got B}_n^\prime$ be a perfect, symmetric bilinear form on $W({\bf F})^{2n}$ that defines the involution $\Mi$ of ${\rm End}(W({\bf F})^{2n})$, cf. Lemma 3.3.1 (a) and (c). Let $V$ and $K$ be as in  Subsection 2.1 and such that there exists a $K$-linear isomorphism $(W({\bf F})^{2n}\otimes_{W({\bf F})} K,{\got B}_n^\prime)\tilde\to (\Ml_n\otimes_{{\bf Z}_{(2)}} K,{\got B}_n)$, to be viewed as an identification. This allows us to identify $\Md_{n,K}$ with the generic fibre of a direct factor of $G^0_{j,V}$ and thus the schematic closure of $\Md_{n,K}$ in ${\bf GL}_{W({\bf F})^{2n}\otimes_{W({\bf F})} V}$ is a reductive group scheme isomorphic to the mentioned factor. From this and Lemma 3.2 we get that the reduction modulo $2V$ of the perfect, symmetric bilinear form ${\got B}_n^\prime$ on $W({\bf F})^{2n}\otimes_{W({\bf F})} V$ is alternating. Thus ${\got B}_n^\prime$ modulo $2W({\bf F})$ is alternating.\endproof 
\medskip\noindent
{\bf 3.5.5. Remark.} Let $G^{\Mi}_{2,j,{\bf F}}$ be the subgroup scheme of $G_{2,j,{\bf F}}$ defined by the property that for each commutative ${\bf F}$-algebra $\sharp$ and for every element $\dag\in G^{\Mi}_{2,j,{\bf F}}(\sharp)$ we have $\Mi(\dag)=\dag^{-1}$. From Lemma 3.5.4 we get that the double monomorphism $G^0_{j,{\bf F}}\hookrightarrow G^0_{1,j,{\bf F}}\hookrightarrow G^{\Mi}_{2,j,{\bf F}}$ over ${\bf F}$ is a product of double monomorphisms of the form ${\bf SO}_{2n,{\bf F}}\hookrightarrow {\bf O}_{2n,{\bf F}}\hookrightarrow {\bf Sp}_{2n,{\bf F}}$.
\bigskip\smallskip\noindent
{\bigsll {\bf 4. Crystalline notations and basic properties}}
\bigskip
Until the end we will use the notations of Subsections 1.3, 2.1, and 3.5 and of Subsubsection 3.5.2; moreover we will take the algebraically closed field $k$ to have a countable transcendental degree over ${\bf F}_2$. Each continuous action of a totally discontinuous locally compact group on a scheme will be in the sense of [10, Subsubsect. 2.7.1] and it will be a right action. In this section we mainly introduce the basic crystalline setting, properties, and strategy needed in the proof of the Main Theorem. Let $d:=\dim_{{\bf C}}(\Mx)$. Let $\Mb^{\rm opp}$ be the opposite ${\bf Z}_{(2)}$-algebra of $\Mb$.
\medskip\smallskip\noindent
{\bf 4.1. Basic setting.} The continuous actions of $G({\bf A}_f^{(2)})$ on ${\rm Sh}(G,\Mx)/H_2$ and $\Mm$ induce an action of $G({\bf A}_f^{(2)})$ on $\Mn$ (which can be easily checked to be continuous). Let $H_0$ be a compact, open subgroup of $G({\bf A}_f^{(2)})$ such that we have $H_0\subseteq {\rm Ker}({\bf GSp}(L,\psi)(\widehat{\bf Z})\to {\bf GSp}(L,\psi)({\bf Z}/m{\bf Z}))$ for some $m\in 1+2{\bf N}^{\ast}$. The group $H_0$ acts freely on $\Mm$ (see Serre's result [33, Sect. 21, Thm. 5]) and the quotient scheme $\Mm/H_0$ exists and is a pro-\'etale cover of the Mumford scheme $\Ma_{{\dim_{\bf Q}(W)\over 2},1,m}$; the natural morphism $\Mm\to\Mm/H_0$ is also a pro-\'etale cover. We get that $H_0$ acts freely on the closed subscheme $\Mn$ of $\Mm_{O_{(v)}}$. This allows us to consider the quotient scheme $\Mn/H_0$. To be short, here we {\it define} $\Mn/H_0$ to be the schematic closure of ${\rm Sh}(G,\Mx)/(H_2\times H_0)$ in $\Mm_{O_{(v)}}/H_0$. 
\smallskip
The inverse image $\Mn/H_0\times_{\Mm_{O_{(v)}}/H_0} \Mm_{O_{(v)}}$ of $\Mn/H_0$ in $\Mm_{O_{(v)}}$ is a pro-\'etale cover of $\Mn/H_0$ and thus it is a flat, closed subscheme of $\Mm_{O_{(v)}}$; therefore it is the schematic closure in $\Mm_{O_{(v)}}$ of its generic fibre which is ${\rm Sh}(G,\Mx)/H_2$ and thus it is $\Mn$ itself. Thus $\Mn$ is a pro-\'etale cover of $\Mn/H_0$. 
\smallskip
We check that the flat $O_{(v)}$-scheme $\Mn/H_0$ is of finite type (to be compared with [39, proof of Prop. 3.4.1]). We consider a compact, open subgroup $K_0$ of ${\bf GSp}(W,\psi)({\bf A}_f^{(2)})$ which contains $H_0$, which is contained in ${\rm Ker}({\bf GSp}(L,\psi)(\widehat{\bf Z})\to {\bf GSp}(L,\psi)({\bf Z}/m{\bf Z}))$ (with $m$ as above), and which has the property that ${\rm Sh}(G,\Mx)/(H_2\times H_0)$ is a closed subscheme of $\Mm_{E(G,\Mx)}/K_0$. Let $\Mn/H_0(K_0)$ be the schematic closure of  ${\rm Sh}(G,\Mx)/(H_2\times H_0)$ in $\Mm_{O_{(v)}}/K_0$. As $\Mm_{O_{(v)}}/K_0$ is an \'etale cover of $\Ma_{{\dim_{\bf Q}(W)\over 2},1,m}$, $\Mn/H_0(K_0)$ is a flat $O_{(v)}$-scheme of finite type. As $\Mm/H_0\to\Mm/K_0$ is a pro-\'etale cover, we have a natural pro-finite morphism $\Mn/H_0\to \Mn/H_0(K_0)$ whose generic fibre is an isomorphism; thus this morphism is birational and induces an isomorphism at the level of normalizations. As the finitely generated $O_{(v)}$-algebras are excellent (see [28, Subsects. (34.A) and (34.B)]), the scheme $\Mn/H_0(K_0)$ is excellent and thus its normalization is a finite $\Mn/H_0(K_0)$-scheme. From the last two sentences we get that the morphism $\Mn/H_0\to \Mn/H_0(K_0)$ is finite and thus the $O_{(v)}$-scheme $\Mn/H_0$ is of finite type. The relative dimension of $\Mn/H_0$ over $O_{(v)}$ is $d$. We will denote $(\Mn/H_0)_k$ and $(\Mn/H_0)_{W(k)}$ by $\Mn_k/H_0$ and $\Mn_{W(k)}/H_0$ (respectively).  
\medskip\noindent
{\bf 4.1.1. On ${\bf Z}_{(2)}$-endomorphisms.} We consider the closed subscheme $\Mn^{\rm mod}$ of $\Mm_{O_{(v)}}$ which is the moduli scheme that parametrizes principally polarized abelian schemes of relative dimension ${\dim_{\bf Q}(W)\over 2}$ over $O_{(v)}$-schemes which are endowed with a ${\bf Z}_{(2)}$-algebra $\Mb$ of ${\bf Z}_{(2)}$-endomorphisms, which have compatible level-$m$ symplectic similitude structures for all $m\in {\bf N}\setminus 2{\bf N}$, which are subject to extra axioms that are modeled on the injective map $f:(G,\Mx)\hookrightarrow ({\bf GSp}(W,\psi),\Ms)$, and which was considered first in [26, Sect. 5]. We recall here the following two key axioms for the case of a connected, affine $O_{(v)}$-scheme ${\rm Spec}\,S$ and of a quadruple $(A_S,\lambda_{A_S},\Mb,(\kappa_m)_{m\in {\bf N}\setminus 2{\bf N}})$ over it of the type mentioned.
\medskip
{\bf (i) (The determinant axiom)}: If $\{\alpha_1,\ldots,\alpha_s\}$ is a ${\bf Z}_{(2)}$-basis for $\Mb$ and if $X_1,\ldots,X_s$ are independent variables, then the determinant of the linear endomorphism $\sum_{j=1}^s X_j\alpha_j$ of the Lie algebra of the abelian scheme $A_S$ over ${\rm Spec}\,S$ is the extension to $S$ of a universal determinant over $O_{(v)}$ that is of a similar nature and that it is associated naturally to the faithful representation $\Mb\hookrightarrow {\rm End}(L_{(2)})$ and the involution $\Mi$ of ${\rm End}(L_{(2)})$ (it does not depend on the given $\Mx$, cf. [26, Lemmas 4.2 and 4.3 and Sect. 5]). 
\smallskip
{\bf (ii)} For each odd prime $l\in {\bf N}$, the symplectic similitude isomorphism $\kappa_{l^\infty}:(W\otimes_{{\bf Q}} {\bf Q}_l,\psi)\tilde\to (T_l(A_S)\otimes_{{\bf Z}_l} {\bf Q}_l,\lambda_{A_S})$ induced by the level-$l^t$ symplectic similitude structures $\kappa_{l^t}$ of $(A_S,\lambda_{A_S})$ with $t\in {\bf N}^{\ast}$, is also a $\Mb\otimes_{{\bf Z}_{(2)}} {\bf Q}$-linear isomorphism (here $\lambda_{A_S}$ denotes also the perfect alternating form on the $l$-adic Tate-module $T_l(A_S)$ of $A_S$ induced by $\lambda_{A_S}$).
\medskip
The axiom (ii) implies that the morphism $\Mn^{\rm mod}\to\Mm_{O_{(v)}}$ is injective on geometric points (i.e., a ${\bf Q}$--endomorphism of $A_S$ is uniquely determined by its $l$-adic \'etale realization with $l$ an odd prime). The Serre--Tate deformation theory implies that  $\Mn^{\rm mod}\to\Mm_{O_{(v)}}$ induces injections at the level of tangent spaces associated to geometric points. Based on the last two sentences, one easily gets that the morphism $\Mn^{\rm mod}\to\Mm_{O_{(v)}}$ which is projective (see loc. cit.) is in fact a closed embedding as stated above.
\smallskip
The generic fibre $\Mn^{\rm mod}_{E(G,\Mx)}$ is a finite disjoint union of Shimura varieties of the form ${\rm Sh}(G,\Mx)/H_2$ and this forms another way to get that ${\rm Sh}(G,\Mx)/H_2$ is a closed subscheme of $\Mm_{E(G,\Mx)}$. Therefore $\Mn$ is a closed subscheme of $\Mn^{\rm mod}$ and this implies that for each $b\in\Mb\subseteq\Mj$ we can speak about the ${\bf Z}_{(2)}$-endomorphism of $\Ma$ that corresponds naturally to $b$. Accordingly, in all that follows we denote also by $\Mb$ the ${\bf Z}_{(2)}$-algebra of ${\bf Z}_{(2)}$-endomorphisms with which each abelian scheme obtained from $\Ma$ by pull back is endowed. For $b\in\Mb$, we denote also by $b$ different de Rham (crystalline) realizations of such ${\bf Z}_{(2)}$-endomorphisms that correspond to $b$. In particular, we will speak about the ${\bf Z}_{(2)}$-monomorphism $\Mb^{\rm opp}\hookrightarrow {\rm End}(M_y)$ that makes $M_y$ to be a $\Mb^{\rm opp}\otimes_{{\bf Z}_{(2)}} W(k)$-module and makes $M_y^\vee$ to be a $\Mb\otimes_{{\bf Z}_{(2)}} W(k)$-module. 
\smallskip
In what follows it is convenient to consider as well the finite level version of $(\Ma,\Lambda,\Mb)$. Thus we will choose the compact, open subgroup $H_0$ of $G({\bf A}_f^{(2)})$ such that the triple $(\Ma,\Lambda,\Mb)$ is the pull back of a similar triple $(\Ma_{H_0},\Lambda_{H_0},\Mb)$ over $\Mn/H_0$  (this holds if $H_0\subseteq {\rm Ker}({\bf GSp}(L,\psi)(\widehat{\bf Z})\to {\bf GSp}(L,\psi)({\bf Z}/l{\bf Z}))$ for an odd prime $l>>0$).
\medskip\noindent
{\bf 4.1.2. The moduli difficulty.} The difficulty we face in the case (D) is that the moduli scheme  $\Mn^{\rm mod}_{W(k)}$ is never formally smooth over $W(k)$ at the point $y\in\Mn_{W(k)}(k)\subseteq \Mn^{\rm mod}_{W(k)}(k)$. We detail this fact in our present context. 
\smallskip
We would like to show that $\Mn_{W(k)}$ is regular and formally smooth over $W(k)$ at each point $y\in\Mn_{W(k)}(k)$. As we are in the case (D), the trouble is that the formal deformation space ${\got D}_y$ over ${\rm Spf}\,k$ of the triple $(A,\lambda_A,\Mb)$ is not formally smooth over ${\rm Spf}\,k$ of dimension $d$. One can check this starting from the fact (see Remark 3.5.5) that the dimension of the subgroup of ${\bf GSp}(L_{(2)}\otimes_{{\bf Z}_{(2)}} {\bf F}_2,\psi)$ that centralizes $\Mb\otimes_{{\bf Z}_{(2)}} {\bf F}_2$ is greater than $\dim(G)$. 
\smallskip
We explain this in the case when $G_{{\bf Z}_2}$ is isomorphic to ${\bf GSO}^+_{2n,{\bf Z}_2}$. The dimension of either $\Mn_{E(G,\Mx)}$ or $\Mn^{\rm mod}_{E(G,\Mx)}$ is $d={{n(n-1)}\over 2}$ (cf. end of Subsubsection 3.5.1 for the last equality). One always expects that Theorem 1.4 (a) holds and thus implicitly (cf. Remark 3.5.5) that the centralizer of $\Mb^{\rm opp}\otimes_{{\bf Z}_{(2)}} W(k)$ in ${\bf Sp}(M_y,\lambda_{M_y})$ has a special fibre $C_k$ which is an ${\bf Sp}_{2n,k}$ group. If there exists a $W(k)$-valued point of $\Mn^{\rm mod}_{W(k)}$ which lifts $y\in \Mn^{\rm mod}_{W(k)}(k)$ (for instance, this would hold if $\Mn^{\rm mod}_{W(k)}$ would be formally smooth over $W(k)$ at $y\in\Mn^{\rm mod}_{W(k)}(k)$), then one can easily check based on Lemma 2.2.3 that $C_k$ is an ${\bf Sp}_{2n,k}$ group. If $C_k$ is an ${\bf Sp}_{2n,k}$ group, then the deformation theories of Subsubsection 1.2.3 imply that the tangent space of ${\got D}_y$ (equivalently of the point $y\in\Mn^{\rm mod}_{W(k)}(k)$) has dimension ${{n(n+1)}\over 2}$ which is greater than $d$ and therefore $\Mn^{\rm mod}_{W(k)}$ is not formally smooth over $W(k)$ at $y\in\Mn^{\rm mod}_{W(k)}(k)$. We conclude that $\Mn^{\rm mod}_{W(k)}$ is not formally smooth over $W(k)$ at $y\in\Mn^{\rm mod}_{W(k)}(k)$.
\smallskip
Thus the deformation theories of Subsubsection 1.2.3 do not suffice to show that $\Mn_{W(k)}$ is formally smooth over $W(k)$ at $y$. This explains why below we will use heavily Section 3 and a few other crystalline theories and new ideas. 
\medskip\noindent
{\bf 4.1.3. Symmetric principal quasi-polarization.}  In this subsubsection we  assume $G_{{\bf Z}_2}$ is isomorphic to ${\bf GSO}_{2n,{\bf Z}_2}^+$. Let $s\in{\bf N}^{\ast}\setminus\{1\}$ and $(\Mv,b_{\Mv})$ be as in Subsubsection 3.5.3. The projection of $L_{(2)}\otimes_{{\bf Z}_{(2)}} {\bf Z}_2$ on each factor $\Mv$ associated naturally to (4) belongs to $\Mb\otimes_{{\bf Z}_{(2)}} {\bf Z}_2$; thus it can be identified with an idempotent ${\bf Z}_2$-endomorphism of $\Ma_{\Mn_{W(k)}}$ (cf. Subsubsection 4.1.1) and therefore also with an idempotent endomorphism of the $2$-divisible group of $\Ma_{\Mn_{W(k)}}$. Thus the principally quasi-polarized $2$-divisible group of $(\Ma,\Lambda)_{\Mn_{W(k)}}$ is of the form
$$(\Me^s,\Lambda_{\Me^s}).$$
\indent
We show the existence of an isomorphism $b_{\Me}:\Me\tilde\to\Me^{\rm t}$ that corresponds naturally to $b_{\Mv}$. We denote also by $b_{\Mv}$ the ${\bf Z}_2$-linear isomorphism $\Mv\tilde\to\Mv^{\vee}$ induced by $b_{\Mv}$. Let 
$$b_{\Mv}^s:L_{(2)}\otimes_{{\bf Z}_{(2)}} {\bf Z}_2=\Mv^s\tilde\to (L_{(2)}\otimes_{{\bf Z}_{(2)}} {\bf Z}_2)^{\vee}=(\Mv^s)^{\vee}=(\Mv^{\vee})^s$$
be the ${\bf Z}_2$-linear isomorphism which is the direct sum of $s$ copies of $b_{\Mv}$. We denote by $b_{\psi}:L_{(2)}\otimes_{{\bf Z}_{(2)}} {\bf Z}_2\tilde\to (L_{(2)}\otimes_{{\bf Z}_{(2)}} {\bf Z}_2)^{\vee}$ the ${\bf Z}_2$-linear isomorphism induced by $\psi$. The composite $b_{\psi}^{-1}\circ b_{\Mv}^s\in {\rm End}(L_{(2)}\otimes_{{\bf Z}_{(2)}} {\bf Z}_2)$ is fixed by $G^{\rm der}_{{\bf Z}_2}=G^0_{{\bf Z}_2}$ and $G_{{\bf Z}_2}/G^{\rm der}_{{\bf Z}_2}=G_{{\bf Z}_2}/G^0_{{\bf Z}_2}$ acts trivially on it (cf. the last sentence of Subsubsection 3.5.3). Thus $b_{\psi}^{-1}\circ b_{\Mv}^s$ is fixed by $G_{{\bf Z}_2}$ and therefore it belongs to $\Mb\otimes_{{\bf Z}_{(2)}} {\bf Z}_2$. Thus we can speak about the ${\bf Z}_2$-isomorphism of $\Ma$ that corresponds to $b_{\psi}^{-1}\circ b_{\Mv}^s$. By composing it with the ${\bf Z}_2$-isomorphism (in fact it is an isomorphism) from $\Ma$ to its dual $\Ma^{\rm t}$  that corresponds to $b_{\psi}$ (equivalently to the principal polarization $\Lambda$ of $\Ma$), we get a ${\bf Z}_2$-isomorphism from $\Ma$ to $\Ma^{\rm t}$ that corresponds naturally to $b_{\Mv}^s$. At the level of $2$-divisible groups we get an isomorphism $b_{\Me}^s:\Me^s\tilde\to (\Me^s)^{\rm t}=(\Me^{\rm t})^s$. The notations match i.e., $b_{\Me}^s$ is the direct sum of $s$ copies of an isomorphism $b_{\Me}:\Me\tilde\to\Me^{\rm t}$ (as one can easily check this over $\Mn_{B(k)}$ using first \'etale cohomology groups with coefficients in ${\bf Z}_2$ as in Subsubsection 3.5.2). We call $b_{\Me}$ a symmetric principal quasi-polarization of $\Me$. 
\smallskip
As $b_{\Mv}$ modulo $2{\bf Z}_2$ is alternating (cf. Lemma 3.5.4), the generic fibre of $b_{\Me}[2]:\Me[2]\tilde\to\Me[2]^{\rm t}$ is a principal quasi-polarization in the sense of Definition 2.2.1. Thus $\lambda_{\Me[2]}:=b_{\Me}[2]$ itself is a principal quasi-polarization.
\medskip\noindent
{\bf 4.1.4. Quasi-sections.} 
As $\Mn_{W(k)}$ is a pro-\'etale cover of the flat $W(k)$-scheme $\Mn_{W(k)}/H_0$ of finite type, the completion of the local ring of $y\in\Mn_{W(k)}(k)$ is a local, complete, noetherian, flat $W(k)$-algebra of residue field $k$. Thus from the Fact 4.1.5 below applied to this completion we get that there exists a lift 
$$z:{\rm Spec}\,V\to\Mn_{W(k)}$$
of $y:{\rm Spec}\,k\to\Mn_{W(k)}$, where $V$ is a discrete valuation ring as in Subsection 2.1. The composite of $z$ with the pro-\'etale cover $\Mn_{W(k)}\to \Mn_{W(k)}/H_0$ is a quasi-section of the $W(k)$-scheme $\Mn_{W(k)}/H_0$ of finite type. We emphasize that $z$ is not necessarily a closed embedding. Let 
$$(A_V,\lambda_{A_V}):=z^*((\Ma,\Lambda)_{\Mn_{W(k)}}).$$ 
\noindent
We denote also by $\lambda_{A_V}$ the perfect form on $(H^1_{dR}(A_V/V))^\vee$ that is the de Rham realization of $\lambda_{A_V}$. Let $F^1_V$ be the Hodge filtration of $H^1_{dR}(A_V/V)$ defined by $A_V$. Under the canonical identification $H^1_{dR}(A_V/V)\otimes_V k=M_y/2M_y$, the $k$-vector space $F^1_V\otimes_V k$ gets identified with the kernel of the reduction of $\Phi_y$ modulo $2W(k)$.
\smallskip
We fix an $O_{(v)}$-monomorphism $i_V:V\hookrightarrow{\bf C}$. Let $w=[h,g_w]\in\Mn({\bf C})=\Mn_{W(k)}({\bf C})$ be the composite of ${\rm Spec}\,{\bf C}\to {\rm Spec}\,V$ with $z$; thus $A_{{\bf C}}$ is $w^*(\Ma)$ of Subsubsection 3.5.2 and we have an identification $L_{(2)}\otimes_{{\bf Z}_{(2)}} {\bf Z}_2=(H^1_{\acute et}(A_{{\bf C}},{\bf Z}_{2}))^\vee$ of $\Mb$-modules that takes $\psi$ to a unit of ${\bf Z}_2$ times the $2$-adic \'etale realization of (the extension to ${\bf C}$ of) $\lambda_{A_V}$. 
\medskip\noindent
{\bf 4.1.5. Fact.} {\it Let $\Mo$ be a flat $W(k)$-algebra which is a local, complete, noetherian ring of residue field $k$ (thus $\Mo$ is a quotient of $W(k)[[t_1,\ldots,t_s]]$ for some $s\in {\bf N}^{\ast}$). Then there exists a $W(k)$-homomorphism $\xi:\Mo\to V$, where $V$ is a discrete valuation ring as in Subsection 2.1.}
\medskip\noindent
{\it Proof:} By replacing $\Mo$ with a quotient of it by a minimal prime ideal of it of characteristic $0$, to prove the existence of $\xi$ we can assume that $\Mo$ is integral. We proceed by induction on $v=\dim(\Mo)\in {\bf N}^{\ast}$. If $v=1$, then the normalization of $\Mo$ is isomorphic to a $W(k)$-algebra $V$ as in Subsection 2.1 and we can take $\xi$ to be the resulting inclusion $\Mo\hookrightarrow V$. For $v\ge 2$, the passage from $v-1$ to $v$ goes as follows. Let $\Lambda_1(\Mo)$ be the set of prime ideals of $\Mo$ of height $1$. For each non-zero element $a$ of the maximal ideal ${\got m}_\Mo$ of $\Mo$, the number of prime ideals in $\Lambda_1(\Mo)$ containing $a$ is a positive integer. As $v\ge 2$, ${\got m}_\Mo$ is not a finite union of elements of $\Lambda_1(\Mo)$ and thus the set $\Lambda_1(\Mo)$ is infinite. From the last two sentences we get that there exists ${\got p}\in \Lambda_1(\Mo)$ such that $\Mo/{\got p}$ has characteristic $0$. We have $\dim(\Mo/{\got p})\le v-1$ and therefore from the inductive assumption we get there exists a $W(k)$-homomorphism $\xi_{\got p}:\Mo/{\got p}\to V$; thus we can take $\xi$ to be the composite of the $W(k)$-epimorphism $\Mo\to\Mo/{\got p}$ with $\xi_{\got p}$. This ends the induction and thus the proof.\endproof
\medskip\smallskip\noindent
{\bf 4.2. Fontaine theory.} We fix an algebraic closure $\bar K$ of $K$. We identify $L_{(2)}\otimes_{{\bf Z}_{(2)}} {\bf Z}_2=(H^1_{\acute et}(A_{{\bf C}},{\bf Z}_{2}))^\vee=(H^1_{\acute et}(A_V\times_V \bar K,{\bf Z}_2))^\vee$ with the Tate-module of the $2$-divisible group of $A_{K}$. We consider the $2$-adic Galois representation  
$$\varrho:{\rm Gal}(\bar K/K)\to {\bf GL}_{L_{(2)}\otimes_{{\bf Z}_{(2)}} {\bf Q}_2}({\bf Q}_2)={\bf GL}_{(H^1_{\acute et}(A_V\times_V \bar K,{\bf Q}_2))^\vee}({\bf Q}_2)$$ 
associated to $A_K$; it factors through the group of ${\bf Q}_2$-valued points of the subgroup $G_{1,{\bf Q}_2}$ of ${\bf GL}_{L_{(2)}\otimes_{{\bf Z}_{(2)}} {\bf Q}_2}$ that fixes $\psi$ and each element of $\Mb$. As $G$ is the identity component of $G_1$ (cf. axiom 1.3 (iii)), the group $(G_1/G)({\bf Q}_2)$ is finite and thus there exists a finite field extension $\tilde K$ of $K$ which is contained in $\bar K$ and such that $\varrho({\rm Gal}(\bar K/\tilde K))$ is a subgroup of $G({\bf Q}_2)$. [One can check using the symplectic similitude level structures of $A_K$ and standard properties of Hodge cycles as in [11] that in fact the image of $\varrho$ is contained in $G({\bf Q}_2)$ and thus that we can take $\tilde K=K$; but this is irrelevant for what follows.] We get that ${\rm Gal}(\bar K/\tilde K)$ fixes (via $\varrho$) each tensor $v_{\alpha}\in\Mt(L_{(2)}\otimes_{{\bf Z}_{(2)}} {\bf Q}_2)$ with $\alpha\in\Mj$. 
\smallskip
Let $\tilde V$ be the ring of integers of $\tilde K$. Let $\tilde B$ be the Fontaine ring of $\tilde K$ as used in [17]; it is a commutative, integral $B(k)$-algebra equipped with a Frobenius endomorphism and a ${\rm Gal}(\bar K/\tilde K)$-action. Under the identification $L_{(2)}\otimes_{{\bf Z}_{(2)}} {\bf Z}_2=(H^1_{\acute et}(A_V\times_V \bar K,{\bf Z}_2))^\vee$, [17, Thm. 6.2] provides us with a functorial $\tilde B$-linear isomorphism 
$$c_y:L_{(2)}\otimes_{{\bf Z}_{(2)}} \tilde B\tilde\to M_y^\vee\otimes_{W(k)} \tilde B$$ which preserves all structures and in particular which is also a $\Mb\otimes_{{\bf Z}_{(2)}} \tilde B$-isomorphism. From Fontaine comparison theory we get that $c_y$ takes each $v_{\al}\in\Mt(W)$ to a tensor $t_{\al}\in\Mt(M_y^\vee[{1\over 2}])$ and takes $\psi$ to a unit $u_{\psi}$ of $\tilde B$ times the alternating form $\lambda_{M_y}$ on $M_y^\vee\otimes_{W(k)} \tilde B$. If $b\in\Mb\subseteq\Mj$, then $t_{b}$ is the de Rham (crystalline) realization $b$ of the ${\bf Z}_{(2)}$-endomorphism $b$ of $A_V$ (of $A$). From the property 1.3 (iv) we get that the group $G_{B(k)}$ is split. From this and the existence of the cocharacter $\mu_h$ of Subsection 3.5, we get that $G_{B(k)}$ has cocharacters that act on the $B(k)$-span of $\psi$ via the identity character of ${\bf G}_{m,B(k)}$. Thus by composing $c_y$ with a $\tilde B$-valued point of the image of such a cocharacter, we get the existence of a symplectic isomorphism 
$$c_y^\prime:(L_{(2)}\otimes_{{\bf Z}_{(2)}} \tilde B,\psi)\tilde\to (M_y^\vee\otimes_{W(k)} \tilde B,\lambda_{M_y})$$
that takes $v_{\al}$ to $t_{\al}$ for all $\al\in\Mj$. As the group $G^0=G\cap {\bf Sp}(W,\psi)$ is connected (see Subsection 3.5), the only class in $H^1({\rm Gal}(\bar K/B(k)),G^0_{B(k)})$ is the trivial one (see [36, Ch. II, Subsect. 3.3 and Ch. III, Subsect. 2.3]). From this and the existence of $c^\prime_y$, we get the existence of a $B(k)$-linear isomorphism
$$j_y:L_{(2)}\otimes_{{\bf Z}_{(2)}} B(k)\tilde\to M_y^\vee[{1\over 2}]\leqno (5)$$
that takes $\psi$ to $\lambda_{M_y}$ and takes $v_{\al}$ to $t_{\al}$ for all $\al\in\Mj$. 
\medskip\smallskip\noindent
{\bf 4.3. Crystalline decompositions and complements.} Until the end we  assume that $G_{{\bf Z}_2}$ is isomorphic to ${\bf GSO}_{2n,{\bf Z}_2}^+$. Thus $d={{n(n-1)}\over 2}$, cf. end of Subsubsection 3.5.1. Let $s\in{\bf N}^{\ast}\setminus\{1\}$ and $(\Mv,b_{\Mv})$ be as in Subsubsection 3.5.3. We consider the pull backs 
$$(\Me_z^s,\lambda_{\Me_z^s})\;\;\;{\rm and}\;\;\; b_z:\Me_z\tilde\to \Me_z^{\rm t} \leqno (6)$$ 
of $(\Me^s,\Lambda_{\Me^s})$ and $b_{\Me}$ (respectively) of Subsubsection 4.1.3 through $z$; thus $(\Me_z^s,\lambda_{\Me_z^s})$ is the principally quasi-polarized $2$-divisible group of $(A_V,\lambda_{A_V})$. But (6) induces a direct sum decomposition of $F$-crystals over $k$
$$(M_y,\Phi_y)=(N_y,\Phi_y)^s\leqno (7)$$
(here $(N_y,\Phi_y)$ is the $F$-crystal of the special fibre $\Me_{z,k}$ of $\Me_z$).
Let $b_{N_y}$ be the perfect bilinear form on $N_y$ that correspond to the special fibre $b_{z,k}$ of $b_z$ via the Dieudonn\'e functor ${\bf D}$. The functorial isomorphism $c_y$ takes $b_{\Mv}$ to the unit $u_{\psi}$ of $\tilde B$ times $b_{N_y}$ (see Subsection 4.2; we recall that in Subsubsection 4.1.3 we have constructed $b_{\Me}$ in terms of $b_{\psi}$ and of ${\bf Z}_2$-endomorphisms of $\Ma$). This implies that $c_y^\prime$ takes the symmetric bilinear form $b_{\Mv}$ on $\Mv$ to the bilinear form $b_{N_y}$ and therefore also $j_y$ takes $b_{\Mv}$ to $b_{N_y}$. Thus the perfect bilinear form $b_{N_y}$ on $N_y$ is symmetric. Let $(\bar N_y,b_{\bar N_y}):=(N_y,b_{N_y})\otimes_{W(k)} k$. From the end of Subsubsection 4.1.3 we get that:
\medskip\noindent
{\bf 4.3.1. Fact.} {\it The isomorphism $b_z[2]:\Me_z[2]\tilde\to\Me_z[2]^{\rm t}$ is a principal quasi-polarization.}
\medskip\noindent
{\bf 4.3.2. The crystalline difficulty.} Based on the Fact 4.3.1, one expects that the perfect bilinear form $b_{\bar N_y}$ on $\bar N_y$ is alternating. In fact, if by chance this would not be true, then it is easy to see that there exist no $W(k)$-valued points of $\Mn_{W(k)}$ that lift $y$. Unfortunately, due to the fact that the natural divided power structure on the ideal $2W(k)$ of $W(k)$ is not nilpotent modulo $4W(k)$, all the present classification results available on commutative finite flat group schemes of $2$-power rank (order) over $V$, either assume some connectivity property or do not relate directly to the \'etale cohomology of their generic fibres. Thus these classification results do not apply directly in the case when the $2$-rank of $A$ is positive. This forms the essence of the crystalline difficulty we encounter; in different aspects, it shows up in connection to all the bullets below.${}^1$ $\vfootnote{1} {The recent manuscripts on $2$-divisible groups over $V$ of E. Lau and T. Liu available at http://arxiv.org/abs/1006.2720 and http://www.math.purdue.edu/\~{}tongliu/research.html (respectively) could be used to overcome the crystalline difficulty in many situations. However, as the two manuscripts were still under review at the time this paper got accepted for publication and are a lot more involved than the simple, self-contained arguments we present below, we will not appeal here to them.}$ 
\medskip\noindent
{\bf 4.3.3. Strategy.} Our strategy to overcome the moduli and crystalline difficulties mentioned above can be summarized in seven steps as follows.
\medskip\noindent
$\bullet$ We use an elementary global versal deformation argument involving $\Mn_k/H_0$ to show that the perfect bilinear form $b_{\bar N_y}$ on $\bar N_y$ is alternating (see Proposition 5.1).
\smallskip\noindent
$\bullet$ We use Proposition 5.1 and Section 3 to prove Theorem 1.4 (a) (see Subsection 5.2).
\smallskip\noindent
$\bullet$ Based on Theorem 1.4 (a) and on Faltings deformation theory, we construct explicitly the supposed to be formal deformation space of $y$ in the reduced scheme $\Mn_{k,{\rm red}}$ of $\Mn_k$ (see Subsections 6.1 to 6.4). More precisely we construct a morphism ${\rm Spec}\,k[[t_1,\ldots,t_d]]\to \Mm_k$ which is supposed to define the mentioned formal deformation space. 
\smallskip\noindent
$\bullet$ Based on Subsections 6.1 to 6.4 and on Section 3, we prove that the formal deformation space of $y$ in $\Mn_{k,{\rm red}}$ is indeed the one supposed to be and thus it is formally smooth over ${\rm Spf}\,k$ (see Subsections 6.5 to 6.7). 
\smallskip\noindent
$\bullet$ By combining Lemma 3.1.4 with the main results of [35], we show that $\Mn_{W(k)}$ is smooth at all ordinary $k$-valued points of it (see Proposition 7.1). 
\smallskip\noindent
$\bullet$ Based on the above mentioned constructions, we show that the ordinary locus of $\Mn_k$ is Zariski dense in $\Mn_k$ (see Subsection 7.2). 
\smallskip\noindent
$\bullet$ Based on the last three bullets, a direct application of a lemma of Hironaka allows to conclude that $\Mn$ is regular and formally smooth over $O_{(v)}$ (see Subsection 7.3). 
\bigskip\smallskip\noindent
{\bigsll {\bf 5. Proof of Theorem 1.4 (a)}}
\bigskip
We use the notations of Subsections 1.3, 2.1, and 3.5 and of Subsubsections 3.5.2, 3.5.3, and 4.1.3. This section proves Theorem 1.4 (a). The main part of the proof is to show that Lemma 3.5.4 and Fact 4.3.1 transfer naturally to the crystalline contexts that pertain to $y$; the below proposition surpasses (in the geometric context of the Main Theorem) the problem 1.1 (i).
\medskip\smallskip\noindent 
{\bf 5.1. Proposition.} {\it We recall that $G_{{\bf Z}_2}\tilde\to {\bf GSO}_{2n,{\bf Z}_2}^+$. Then for each point $y\in\Mn_{W(k)}(k)$, the perfect bilinear form $b_{\bar N_y}$ on $\bar N_y$ is alternating (see Subsection 4.3 for notations).}
\medskip\noindent
{\it Proof:} We recall that $b_{N_y}$ is symmetric. If $V=W(k)$, then based on Fact 4.3.1 the proposition follows from Lemma 2.2.3 applied to $(\Me_z[2],b_z[2])$. The proof of the proposition in the general case (of an arbitrary index $[V:W(k)]\in{\bf N}^{\ast}$) is lengthy and thus below we will number its main parts. In Subsubsections 5.1.1 to 5.1.5 we will use only properties 1.3 (i) and (iii); but in Subsubsection 5.1.6 we will use also the property 1.3 (iv) and thus implicitly Lemma 3.5.4. We fix an irreducible component ${\rm Irr}_y$ of the reduced scheme of $\Mn_k/H_0$ that passes through the  $k$-valued point of $\Mn_k/H_0$ defined by $y$. Let ${\rm Irr}^0_y$ be an affine, open subscheme of ${\rm Irr}_y$ that contains the $k$-valued point of ${\rm Irr}_y$ defined by $y$. Let ${\rm Spec}\,S$ be a connected \'etale cover of a non-empty, affine, smooth, open subscheme of ${\rm Irr}^0_y$ and let $y_S$ be an arbitrary closed (thus $k$-valued) point of ${\rm Spec}\,S$. 
\medskip\noindent 
{\bf 5.1.1. Part I: notations.}  Due to the choice of $H_0$ before Subsubsection 4.1.2 and to the definition of $b_{\Me}$ in terms of $\Lambda$ and of ${\bf Z}_2$-endomorphisms of $\Ma$ that correspond to elements of $\Mb\otimes_{{\bf Z}_{(2)}} {\bf Z}_2$, we have a direct sum decomposition $(\Ma_{H_0}\times_{\Mn/H_0} \Mn_{W(k)}/H_0)[2^{\infty}]=\Me_0^s$ and an isomorphism $b_{\Me_0}:\Me_0\tilde\to \Me_0^{\rm t}$ of $2$-divisible groups whose pull backs to $\Mn_k\times_{\Mn_k/H_0} {\rm Spec}\, S$ are natural pull backs of the decomposition $\Ma_{\Mn_{W(k)}}[2^{\infty}]=\Me^s$ and of the isomorphism $b_{\Me}:\Me\tilde\to\Me^{\rm t}$ (respectively) of the Subsubsection 4.1.3. Let $E:=\Me_{0,S}[2]$ and $\lambda_E:E\tilde\to E^{\rm t}$ be the pull backs of $\Me_0[2]$ and $\lambda_{\Me_0[2]}:=b_{\Me_0}[2]$ (respectively) to ${\rm Spec}\,S$. From Fact 4.3.1 applied not to $z$ but to a morphism $\tilde z:{\rm Spec}\,\tilde V\to\Mn_{W(k)}$ whose image in $\Mn_k/H_0$ is the image of $y_S$ in $\Mn_k/H_0$ (here $\tilde V$ is a finite, discrete valuation ring extension of $W(k)$), we get that $\lambda_E$ is a principal quasi-polarization modulo the maximal ideal of $S$ that defines the closed point $y_s$ of ${\rm Spec}\,S$. By varying the closed point $y_S$ of ${\rm Spec}\,S$, we get that $\lambda_E$ is a principal quasi-polarization. We consider the evaluation
$$(N,\phi,\upsilon,\nabla,b_N)$$ 
of ${\bf D}(E,\lambda_E)$ at the trivial thickening of ${\rm Spec}\,S$. Thus $N$ is a projective $S$-module, $\phi:N^{(2)}\to N$ and $\upsilon:N\to N^{(2)}$ are $S$-linear maps, $\nabla$ is a connection on $N$, and $b_N$ is a symmetric bilinear form on $N$. Let $F:={\rm Ker}(\phi)$; it is the Hodge filtration of $N$ defined by $E$. We consider the $S$-submodule
$$X:=\{x\in N|b_N(x,x)=0\}\subseteq N.$$
\indent
Let $\tilde {\rm Irr}_y$ be the normalization of ${\rm Irr}_y$ in the field of fractions of $S$; ${\rm Spec}\,S$ is an open, Zariski dense subscheme of it. Let $\tilde y$ be a $k$-valued point of $\tilde {\rm Irr}_y$ whose image in ${\rm Irr}_y$ is the $k$-valued point of ${\rm Irr}_y$  defined by $y$. As ${\rm Spec}\, S$ maps onto a Zariski dense, open subscheme of the affine, integral scheme ${\rm Irr}_y$ and as $\lambda_E$ is induced by the isomorphism $b_{\Me_0,{\rm Irr}_y}$ that lifts $b_{z,k}$, we get that $b_N$ specializes in a natural way to $b_{\bar N_y}$ via the crystalline realization of the isomorphism $\lambda_{\Me_0[2],\tilde {\rm Irr}_y}$ and via the pull back through $\tilde y$ of this crystalline realization. 
\smallskip
Thus to prove the proposition, using specializations we get that it suffices to show that $b_N$ is alternating. To check this thing, we will often either shrink ${\rm Spec}\,S$ (i.e., we will replace ${\rm Spec}\,S$ by an open, dense subscheme of it) or replace ${\rm Spec}\,S$ by a connected \'etale cover of it. We will show that the assumption that $b_N$ is not alternating, leads to a contradiction (the argument will end up before Subsection 5.2). This assumption implies that $X\neq N$. 
\smallskip
As $b_N(F,F)=0$, we have $F\subseteq X$. As $S$ is a finitely generated, integral $k$-algebra, by shrinking ${\rm Spec}\,S$, we can assume that $N$ and $F$ are free $S$-modules, that the reduction of $b_N$ modulo each maximal ideal of $S$ is not alternating, and that we have a short exact sequence of nontrivial free $S$-modules
$$0\to X\to N\to N/X\to 0.\leqno (8)$$ 
We have ${\rm rk}_S(N)=2n$ and (cf. the existence of $\lambda_E$) ${\rm rk}_S(F)=n$. Let $r:={\rm rk}_S(N/X)$; thus ${\rm rk}_S(X)=2n-r$. As $F\subseteq X$, we have $r\in\{1,\ldots,n\}$. As the field of fractions of $S$ is not perfect, we can not use Fact 3.1.1 (a) and (b) to get directly that $r$ is $1$. 
\medskip\noindent
{\bf 5.1.2. Part II: Dieudonn\'e theory.} The short exact sequence (8)
is preserved by $\phi$ and $\upsilon$ and we have $\phi(N^{(2)})\subseteq X$ and $\upsilon(N)\subseteq X^{(2)}$. Thus the resulting Frobenius and Verschiebung maps of $N/X$ are both $0$ maps. As $b_N$ is parallel with respect to the connection $\nabla$ and as $X$ is uniquely determined by $b_N$, $\nabla$ restricts to a connection on $X$. By shrinking ${\rm Spec}\,S$, we can  assume that $S$ has a finite $2$-basis in the sense of [5, Def. 1.1.1]. From the fully faithfulness part of [5, Thm. 4.1.1 or 4.2.1] and the last three sentences, we get that to (8) corresponds a short exact sequence 
$$0\to \pmb{\al}_{2,S}^{r}\to E\to E_1\to 0$$
of commutative, finite, flat group schemes over $S$. Let 
$$E_2:=\lambda_E^{-1}(E_1^{\rm t})\subset \lambda_E^{-1}(E^{\rm t})=E.$$ 
\indent
The quotient group scheme $E/E_2$ is isomorphic to $(\pmb{\al}_{2,S}^r)^{\rm t}$ and therefore to $\pmb{\al}_{2,S}^r$. Let $Y$ be the direct summand of $N$ that defines the evaluation of ${\rm Im}({\bf D}(E/E_2)\to{\bf D}(E))$ at the trivial thickening of ${\rm Spec}\,S$. As $E/E_2\tilde\to \pmb{\al}_{2,S}^r$, we have $\phi(Y\otimes 1)=0$. Thus $Y\subseteq F$ and therefore the $S$-linear map $Y\to N/X$ is $0$. This implies that the homomorphism $\pmb{\al}_{2,S}^r\to E/E_2$ over $S$ is trivial, cf. [5, Thm. 4.1.1 or 4.2.1]. Thus the subgroup scheme $\pmb{\al}_{2,S}^r$ of $E$ is in fact a subgroup scheme of $E_2$. 
\smallskip
The perpendicular of $Y$ with respect to $b_N$ is $X$, cf. constructions. Thus $b_N$ induces naturally a perfect form $b_{\tilde N}$ on $\tilde N:=X/Y$ that (cf. definition of $X$) is alternating.  Let 
$$\tilde E:=E_2/\pmb{\al}_{2,S}^r.$$ 
Let $\lambda_{\tilde E}:\tilde E\tilde\to\tilde E^{\rm t}$ be the isomorphism of finite, flat group schemes over $S$ induced by $\lambda_E$. The evaluation of ${\bf D}(\tilde E,\lambda_{\tilde E})$ at the trivial thickening of ${\rm Spec}\,S$ is $(\tilde N,\tilde\phi,\tilde\upsilon,\tilde\nabla,b_{\tilde N})$, where $(\tilde\phi,\tilde\upsilon,\tilde\nabla)$ is obtained from $(\phi,\upsilon,\nabla)$ via a natural passage to quotients. 
\medskip\noindent
{\bf 5.1.3. Part III: computing $r$.} By shrinking ${\rm Spec}\,S$ we can assume that the $S$-module $\Omega_{S/k}$ of relative differentials is free of rank $d={{n(n-1)}\over 2}$. As $\Mm_k$ is a moduli space of principally polarized abelian schemes endowed with extra symplectic similitude structures, the abelian scheme $\Ma_k$ and thus also its $2$-divisible group is versal at each $k$-valued point of $\Mm_k$. From this and the fact that ${\rm Spec}\,S$ is an \'etale cover of a locally closed subscheme of $\Mm_k/H_0$ which is smooth, integral, and of dimension $d$, we get that the Kodaira--Spencer map $\Mh$ of $\nabla$ (i.e., of $E$) is injective modulo each maximal ideal of $S$ and thus its image ${\rm Im}(\Mh)$ is a free $S$-module of rank $d$ which is a direct summand of ${\rm Hom}(F,N/F)$. The connection $\nabla$ restricts to connections on $Y$ and $X$. From this and the inclusions $Y\subseteq F\subseteq X$, we get that ${\rm Im}(\Mh)$ is in fact a direct summand of the direct summand ${\rm Hom}(F/Y,X/F)$ of ${\rm Hom}(F,N/F)$ and thus it is canonically identified with the image of the Kodaira--Spencer map $\tilde\Mh$ of $\tilde E$. As $b_{\tilde N}$ is alternating and $r\ge 1$, it is well known that we have ${\rm rk}_S({\rm Im}(\tilde\Mh))\le {{(n-r)(n-r+1)}\over 2}\le d$. [The argument for the inequality ${\rm rk}_S({\rm Im}(\tilde\Mh))\le {{(n-r)(n-r+1)}\over 2}$ is the same as the one checking that the formal deformation space of a principally quasi-polarized $2$-divsible group over $k$ of height $2(n-r)$ has (to be compared with 1.1 (b)) a tangent space of dimension ${{(n-r)(n-r+1)}\over 2}$]. 
Therefore the three numbers  $d={{n(n-1)}\over 2}$, ${{(n-r)(n-r+1)}\over 2}$, and ${\rm rk}_S({\rm Im}(\tilde\Mh))={\rm rk}_S({\rm Im}(\Mh))$ must be equal. Thus $r=1$ and the pair $(\tilde E,\lambda_{\tilde E})$ is a versal deformation at each $k$-valued point of ${\rm Spec}\,S$ (i.e., and the Kodaira--Spencer map of $\tilde E$ is injective modulo each maximal ideal of $S$). As $r=1$ and $n\ge 2$, we have $\tilde N\neq 0$.
\medskip\noindent 
{\bf 5.1.4. Part IV: reduction to an ordinary context.} We check that by shrinking ${\rm Spec}\,S$ and by passing to an \'etale cover of ${\rm Spec}\,S$, we can  assume there exists a short exact sequence 
$$0\to \pmb{\mu}_{2,S}^{n-1}\to \tilde E\to ({\bf Z}/2{\bf Z})_S^{n-1}\to 0.\leqno (9)$$ 
This well known property is a consequence of the versality part of the end of Subsubsection 5.1.3. For the sake of completeness, we will include a self contained argument for (9). 
\smallskip
Let $\tilde F:=F/Y\subseteq \tilde N$. Let $J$ be a fixed maximal ideal of $S$. Let $(\tilde N_k,\tilde\phi_k,\tilde\upsilon_k,\tilde F_k,b_{\tilde N_k})$ be the reduction of $(\tilde N,\tilde\phi,\tilde\upsilon,\tilde F,b_{\tilde N})$ modulo $J$. By shrinking ${\rm Spec}\,S$, we can assume that $\tilde N$ and $\tilde F$ are free $S$-modules. We fix an isomorphism 
$$I_S:(\tilde N_k,\tilde F_k,b_{\tilde N_k})\otimes_k S\tilde\to (\tilde N,\tilde F,b_{\tilde N})$$ 
of filtered symplectic spaces over $S$. Let $\tilde N_k=\tilde F_k\oplus \tilde Q_k$ be a direct sum decomposition such that $\tilde Q_k$ is a maximal isotropic subspace of $\tilde N_k$ with respect to $b_{\tilde N_k}$. We define $\tilde H:={\bf Sp}(\tilde N_k,b_{\tilde N_k})$. Any two standard symplectic $S$-bases of $\tilde N_k\otimes_k S$ with respect to $b_{\tilde N_k}$ are $\tilde H(S)$-conjugate. This implies that there exist elements $h_S$, $\tilde h_S\in \tilde H(S)$ such that under $I_S$, $\tilde\phi$ and $\tilde\upsilon$ become isomorphic to $h_S\circ (\tilde\phi_k\otimes 1_S)$ and $(\tilde\upsilon_k\otimes 1_S)\circ\tilde h_S^{-1}h_S^{-1}$ (respectively). As $\tilde\upsilon\circ \tilde\phi=0$, we get that $\tilde h_S$ normalizes ${\rm Im}(\tilde\phi_k)\otimes_k S={\rm Ker}(\tilde\upsilon_k)\otimes_k S$. As for all $x\in\tilde N_k^{(2)}$ and $u\in\tilde N_k$ we have identities $b_{\tilde N_k}(h_S(\tilde\phi_k(x\otimes 1)),u\otimes 1)=b_{\tilde N_k}(x\otimes 1,(\tilde\upsilon_k\otimes 1_S)\tilde h_S^{-1}h_S^{-1}(u\otimes 1))$ and $b_{\tilde N_k}(\tilde\phi_k(x\otimes 1)),u\otimes 1)=b_{\tilde N_k}(x\otimes 1,\tilde\upsilon_k(u\otimes 1))$, we get $b_{\tilde N_k}(\tilde\phi_k(x\otimes 1),\tilde h_S(u\otimes 1)-u\otimes 1)=0$; thus $\tilde h_S(u\otimes 1)-u\otimes 1\in {\rm Im}(\tilde\phi_k)\otimes_k S={\rm Ker}(\tilde\upsilon_k)\otimes_k S$. From this and the fact that $\tilde h_S$ normalizes ${\rm Ker}(\tilde\upsilon_k)\otimes_k S$ we get that $(\tilde\upsilon_k\otimes 1_S)\circ\tilde h_S^{-1}=\tilde\upsilon_k\otimes 1_S$. Thus we can  assume that $\tilde h_S$ is the identity element of $\tilde H(S)$. 
\smallskip 
Let $\tilde h_o\in \tilde H(k)$ be such that $\tilde h_o\tilde\phi_k((\tilde Q_k)^{(2)})=\tilde Q_k$; this implies that $(\tilde N_k,\tilde h_o\tilde\phi_k,\tilde\upsilon_k \tilde h_o^{-1})$ is an ordinary Dieudonn\'e module over $k$. Thus there exists a non-empty open subscheme $\tilde\Mo$ of $\tilde H$ such that for each element $\tilde g\in\tilde\Mo(k)$, the triple $(\tilde N_k,\tilde g\circ\tilde\phi_k,\tilde\upsilon_k\circ\tilde g^{-1})$ is an ordinary Dieudonn\'e module over $k$. Let $\tilde P$ be the parabolic subgroup of $\tilde H$ that normalizes $\tilde F_k$. Let $\tilde C$ be the Levi subgroup of $\tilde P$ that normalizes $\tilde Q_k$. For $\tilde h\in \tilde P(k)$, we write $\tilde h=\tilde h_u\tilde c$, where $\tilde h_u$ is a $k$-valued point of the unipotent radical of $\tilde P$ and where $\tilde c\in\tilde C(k)$. Let
$$j_S:\tilde P\times_k {\rm Spec}\,S\to \tilde H$$ 
be the morphism that takes $(\tilde h,\tilde m)\in \tilde P(k)\times {\rm Spec}\,S(k)$ to $\tilde u:=\tilde h(h_S\circ \tilde m)((\tilde c)^{-1}\circ \sg)\in \tilde H(k)$. As ${\rm Im}(\tilde\Mh)$ has rank $d={{n(n-1)}\over 2}=\dim(S)$, the composite of $h_S:{\rm Spec}\,S\to \tilde H$ with the quotient epimorphism $\tilde H\twoheadrightarrow\tilde H/\tilde P$ is \'etale. Thus $j_S$ induces $k$-linear isomorphisms at the level of tangent spaces of $k$-valued points. Therefore $j_S$ is \'etale and thus $j_S^*(\tilde\Mo)$ is a non-empty open subscheme of $\tilde P\times_k {\rm Spec}\,S$. Therefore there exists a pair $(\tilde h,\tilde m)\in \tilde P(k)\times {\rm Spec}\,S(k)$ such that we have $\tilde u\in\tilde\Mo(k)$. The reduction of $(\tilde N,\tilde\phi,\tilde\upsilon)$ modulo the maximal ideal $\tilde J$ of $S$ that defines $\tilde m$, is isomorphic to $(\tilde N_k,(h_S\circ\tilde m)\tilde\phi_k,\tilde\upsilon_k(h_S\circ\tilde m)^{-1})$ and thus also to $(\tilde N_k,\tilde h(h_S\circ\tilde m)\tilde\phi_k\tilde h^{-1},\tilde h\tilde\upsilon_k(h_S\circ\tilde m)^{-1}\tilde h^{-1})=(\tilde N_k,\tilde u\tilde\phi_k,\tilde\upsilon_k\tilde u^{-1})$. Therefore the reduction of $\tilde E$ modulo $\tilde J$ is an ordinary truncated Barsotti--Tate group of level $1$ over $k$. Thus, by shrinking ${\rm Spec}\,S$, we can  assume that $\tilde E$ is an ordinary truncated Barsotti--Tate group of level $1$ over $k$. Therefore the short exact sequence (9) exists after a passage to an \'etale cover of ${\rm Spec}\,S$. 
\medskip\noindent
{\bf 5.1.5. Part V: filtrations.} Due to (9) and the existence of the short exact sequence $0\to \tilde E\to E_1\to \pmb{\al}_{2,S}\to 0$, we have naturally another short exact sequence $0\to \pmb{\mu}_{2,S}^{n-1}\to E_1\to \pmb{\al}_{2,S}\times_S ({\bf Z}/2{\bf Z})_S^{n-1}\to 0$. Due to (9) and the existence of a short exact sequence $0\to \pmb{\al}_{2,S}\to E_2\to\tilde E\to 0$, we get that we have naturally another short exact sequence $0\to \pmb{\al}_{2,S}\times_S \pmb{\mu}_{2,S}^{n-1}\to E_2\to ({\bf Z}/2{\bf Z})_S^{n-1}\to 0$. Thus $E$ has a filtration 
$$0\subset E_{(1)}\subset  E_{(2)}\subset  E,\leqno (10)$$
where $E_{(1)}$ is $\pmb{\mu}_{2,S}^{n-1}$, $E_{(2)}/E_{(1)}$ is the extension of $\pmb{\al}_{2,S}$ by $\pmb{\al}_{2,S}$, and $E/E_{(2)}$ is $({\bf Z}/2{\bf Z})_S^{n-1}$.
\smallskip
In order to benefit from the previous notations, until Subsection 5.2 we choose $y:{\rm Spec}\,k\to\Mn_{W(k)}$ such that it defines a morphism ${\rm Spec}\,k\to\Mn_{W(k)}/H_0$ that factors through ${\rm Spec}\,S$. We recall that $\Me_{z,k}$ is the special fibre of $\Me_z$. The $2$-ranks of $\Me_{z,k}$ and $\Me_{z,k}^{\rm t}$ are $n-1$ (see (10)). It is well known that this implies that $\Me_z$ has a filtration by $2$-divisible groups
$$0\subset\Me_{(1),z}\subset\Me_{(2),z}\subset\Me_z,\leqno (11)$$
where $\Me_{(1),z}=\pmb{\mu}_{2^{\infty},V}^{n-1}$, $\Me_{z}/\Me_{(2),z}=({\bf Q}_2/{\bf Z}_2)_V^{n-1}$, and the $2$-divisible group $\Me_{z}^\prime:=\Me_{(2),z}/\Me_{(1),z}$ over $V$ is connected, has a connected Cartier dual, and its height is equal to $2=2n-2(n-1)$. 
\medskip\noindent
{\bf 5.1.6. Part VI: the new pair $(\Me_z^\prime,b_z^\prime)$.} The filtration (11) is compatible in the natural way with the isomorphism $b_z:\Me_z\tilde\to\Me_z^{\rm t}$ and thus it induces naturally an isomorphism $b_z^\prime:\Me^\prime_z\tilde\to \Me_z^{\prime t}$. Let $T_2(\Me^\prime_{z,K})$ be the Tate-module of the generic fibre $\Me^\prime_{z,K}$ of $\Me^\prime_z$. As $\Me_z^\prime$ has height $2$, $T_2(\Me^\prime_{z,K})$ is a free ${\bf Z}_2$-module of rank $2$ on which ${\rm Gal}(\bar K/K)$-acts. To $b_z^\prime$ corresponds a perfect, symmetric bilinear form (denoted also by) $b_z^\prime$ on $T_2(\Me^\prime_{z,K})$ that modulo $2{\bf Z}_2$ is alternating (cf. Fact 4.3.1).  Let ${\bf GSO}^\prime$ be the schematic closure in ${\bf GL}_{T_2(\Me^\prime_{z,K})}$ of the identity component of the subgroup of ${\bf GL}_{T_2(\Me^\prime_{z,K})[{1\over 2}]}$ that normalizes the ${\bf Q}_2$-span of $b_z^\prime$. As ${\bf GSO}^\prime_{W({\bf F})}$ is isomorphic to ${\bf GSO}_{2,W({\bf F})}$ (cf. Proposition 3.4 (c)), ${\bf GSO}^\prime$ is a torus of rank $2$. The Galois representation ${\rm Gal}(\bar K/K)\to {\bf GL}_{T_2(\Me^\prime_{z,K})}({\bf Z}_2)$ induced by $\varrho$ of Subsection 4.2 factors through ${\bf GSO}^\prime({\bf Z}_2)$. Let $\Mb^\prime$ be the semisimple, commutative ${\bf Z}_2$-algebra of endomorphisms of $\Me^\prime_{z}$ whose Lie algebra is ${\rm Lie}({\bf GSO}^\prime)$. 
\smallskip
Let $e$ and $R_e$ be as in  Subsection 2.1. Let $(N_z^\prime,\Phi_z^\prime,\Upsilon_z^\prime,\nabla_{z}^\prime,b_{N_z^\prime})$ be the projective limit indexed by $m\in{\bf N}^{\ast}$ of the evaluations of ${\bf D}((\Me_z^\prime,b_{\Me_z}^\prime)_{U_e})$ at the thickenings attached naturally to the closed embeddings ${\rm Spec}\,U_e\hookrightarrow {\rm Spec}\,R_e/2^mR_e$. Thus $N_z^\prime$ is a free $R_e$-module of rank $2$, we have $R_e$-linear maps $\Phi_z^\prime:N_z^\prime\otimes_{R_e} {}_{\Phi_{R_e}} R_e\to N_z^\prime$ and $\Upsilon_z^\prime:N_z^\prime\to N_z^\prime\otimes_{R_e} {}_{\Phi_{R_e}} R_e$, $\nabla_z^\prime:N_z^\prime\to N_z^\prime\otimes_{R_e} R_e\, dt$ is an integrable and topologically nilpotent connection, and $b_{N_z^\prime}$ is a perfect, symmetric bilinear form on $N_z^\prime$. The ${\bf Z}_2$-algebra $\Mb^\prime$ acts on $N_z^\prime$ and the $R_e\otimes_{{\bf Z}_2} \Mb^\prime$-module $N_z^\prime$ is free of rank $1$. 
\smallskip
To the natural decomposition $\Mb^\prime\otimes_{{\bf Z}_2} W({\bf F})\tilde\to W({\bf F})\oplus W({\bf F})$ of $W({\bf F})$-algebras corresponds a direct sum decomposition $N_z^\prime=N_z^{\prime (1)}\oplus N_z^{\prime (2)}$ into free $R_e$-modules of rank $1$. As $N_z^{\prime (1)}$ and $N_z^{\prime (2)}$ are isotropic with respect to $b_{N_z^\prime}$, $b_{N_z^\prime}$ modulo $2R_e$ is alternating (cf. Fact 3.1.1 (a)). Thus $(N_z^\prime,b_{N_z^\prime})\otimes_{R_e} k$ is a symplectic space over $k$ of rank $2$. 
\smallskip
The fibre of the filtration (11) over $k$ splits. Thus we have a direct sum decomposition 
$$(\bar N_y,b_{\bar N_y})=[(N_z^\prime,b_{N_z^\prime})\otimes_{R_e} k]\oplus (N_o,b_{N_o})\leqno (12)$$ 
over $k$, where $(N_o,b_{N_o})$ is part of a quadruple $(N_o,\phi_{N_o},\upsilon_{N_o},b_{N_o})$ that is the evaluation at the trivial thickening of ${\rm Spec}\,k$ of the Dieudonn\'e functor ${\bf D}$ applied to an ordinary truncated Barsotti--Tate group $E_o$ of level $1$ over $k$ equipped with an isomorphism $\lambda_{E_o}:E_o\tilde\to E_o^{\rm t}$. As $E_o\tilde\to({\bf Z}/p{\bf Z})_k^{n-1}\oplus\pmb{\mu}_{2,k}^{n-1}$, we get that $(N_o,b_{N_o})$ is a symplectic space over $k$. 
\smallskip
From the symplectic space properties of the last two paragraphs and from (12) we get that $(\bar N_y,b_{\bar N_y})$ is a symplectic space over $k$. Thus the reduction of $b_N$ modulo the maximal ideal of $S$ associated naturally to $y$ is alternating. In other words, we reached the desired contradiction. Thus $b_N$ is an alternating form on $N$. This proves the proposition.\endproof
\medskip\smallskip\noindent
{\bf 5.2. End of the proof of Theorem 1.4 (a).} Due to Proposition 5.1, the formula $q_{N_y}(x):={{b_{N_y}(x,x)}\over 2}$ defines a quadratic form on $N_y$. Let $\Mg_y$ be the schematic closure in ${\bf GSp}(M_y,\lambda_{M_y})$ of the identity component of the subgroup of ${\bf GSp}(M_y[{1\over 2}],\lambda_{M_y})$ that fixes $t_{b}$ for all $b\in\Mb\subseteq\Mj$. Based on (7), we can redefine $\Mg_y$ as ${\bf GSO}(N_y,q_{N_y})$ and thus $\Mg_y$ is a reductive group scheme (cf. Proposition 3.4 (c)). 
\smallskip
Let $j_y$ be as in Subsection 4.2. Let $L_y:=j_y^{-1}(M_y^\vee)$; it is a $W(k)$-lattice of $L_{(2)}\otimes_{{\bf Z}_{(2)}} B(k)$. We recall that $j_y$ takes $\psi$ to $\lambda_{M_y}$ and takes $b=v_{b}$ to $b=t_b$ for all $b\in\Mb$ (cf. (5); see Subsubsection 3.5.2 and Subsubsection 4.1.1 for these identities). Thus the $W(k)$-lattice $L_y$ of $L_{(2)}\otimes_{{\bf Z}_{(2)}} B(k)$ has the following three properties: 
\medskip\noindent
{\bf (i)} for all $b\in\Mb\otimes_{{\bf Z}_{(2)}} W(k)$ we have $b(L_y)\subseteq L_y$, 
\smallskip\noindent
{\bf (ii)} the schematic closure of $G_{B(k)}$ in ${\bf GL}_{L_y}$ is the reductive group scheme $j_y^{-1}\Mg_y j_y$,  
\smallskip\noindent
{\bf (iii)} and we get a perfect, alternating form $\psi:L_y\otimes_{W(k)} L_y\to W(k)$. 
\medskip
In the next paragraph we check that properties (i) to (iii) imply that there exists $g_y\in G^0(B(k))$ such that we have $g_y(L_{(2)}\otimes_{{\bf Z}_{(2)}} W(k))=L_y$. For this we consider a maximal (split) torus $T_{W(k)}$ of $G_{W(k)}$ which (cf. property (ii)) is also a maximal torus of $j_y^{-1}\Mg_y j_y$. The existence of $T_{W(k)}$ follows from the following two facts:
\medskip
$\bullet$ each two points of the building of $G_{B(k)}=(j_y^{-1}\Mg_y j_y)_{B(k)}$, belong to an apartment of the building (see [38, pp. 43--44] for this fact and for the language used) and each apartment  of the building is fixed by the $W(k)$-valued points of a uniquely determined torus over $W(k)$ whose generic fibre is a maximal torus of $G_{B(k)}$ (see [38, Sect. 1 and Subsect. 2.1]);
\smallskip
$\bullet$ there exists a unique hyperspecial point of the building of $G_{B(k)}$ which is fixed by $G_{W(k)}$ (resp. by $j_y^{-1}\Mg_y j_y$) (see [38, Subsubsect. 3.8.1]).
\medskip
 Let $L_{(2)}\otimes_{{\bf Z}_{(2)}} W(k)=\oplus_{i\in \chi} W_i$ and $L_y=\oplus_{i\in \chi} W_i^\prime$ be direct sum decompositions indexed by a set $\chi$ of characters of $T_{W(k)}$, such that for all $i\in \chi$ the actions of $T_{W(k)}$ on $W_i$ and $W_i^\prime$ are via the character $i$ of $T_{W(k)}$. The representation of $\Mb\otimes_{{\bf Z}_{(2)}} B(k)$ on $W_i[{1\over 2}]=W_i^\prime[{1\over 2}]$ is absolutely irreducible. Thus the representations of $\Mb\otimes_{{\bf Z}_{(2)}} W(k)$ on $W_i$ and $W_i^\prime$ are isomorphic and their fibres are absolutely irreducible. Therefore there exists $n_i\in{\bf Z}$ such that we have $p^{n_i}W_i=W_i^\prime$. Let $\mu_0:{\bf G}_{m,B(k)}\to {\bf GL}_{L_{(2)}\otimes_{{\bf Z}_{(2)}} B(k)}$ be a cocharacter such that it acts on $W_i[{1\over 2}]$ via the $n_i$-th power of the identity character of ${\bf G}_{m,B(k)}$. Let $g_y:=\mu_0(p)$. We have $g_y(L_{(2)}\otimes_{{\bf Z}_{(2)}} W(k))=L_y$. As $T_{W(k)}$ is a torus of ${\bf GSp}(L_{(2)}\otimes_{{\bf Z}_{(2)}} W(k),\psi)$, for each $i\in \chi$ there exists a unique $\tilde i\in \chi$ such that for every $i^\prime\in \chi\setminus\{\tilde i\}$ we have $\psi(W_i,W_{i^\prime})=0$. The map $\chi\to \chi$ that takes $i$ to $\tilde i$ is a bijection of order at most $2$. But as $L_y$ and $L_{(2)}\otimes_{{\bf Z}_{(2)}} W(k)$ are both self-dual $W(k)$-lattices of $W\otimes_{{\bf Q}} B(k)$ with respect to $\psi$ (cf. property (iii)), for all $i\in \chi$ we have $n_i+n_{\tilde i}=0$. Therefore $\mu_0$ fixes $\psi$. As for each $i\in \chi$ the $W(k)$-module $W_i$ is left invariant by $\Mb\otimes_{{\bf Z}_{(2)}} W(k)$, $\mu_0$ fixes all elements $b\in\Mb$. Thus $\mu_0$ is a cocharacter of $G^0_{1,B(k)}$ and therefore also of the identity component $G^0_{B(k)}$ of $G^0_{1,B(k)}$. Thus $g_y\in G^0(B(k))$ i.e., $g_y$ exists. 
\smallskip
By replacing $j_y$ with $j_yg_y$, we can assume that $j_y(L_{(2)}\otimes_{{\bf Z}_{(2)}} W(k))=j_y(L_y)=M_y^\vee$. Thus $j_y:L_{(2)}\otimes_{{\bf Z}_{(2)}} W(k)\tilde\to M_y^\vee$ is an isomorphism of $\Mb\otimes_{{\bf Z}_{(2)}} W(k)$-modules that induces a symplectic isomorphism $(L_{(2)}\otimes_{{\bf Z}_{(2)}} W(k),\psi)\tilde\to (M_y^\vee,\lambda_{M_y})$. Thus Theorem 1.4 (a) holds.\endproof 
\bigskip\smallskip\noindent
{\bigsll {\bf 6. The study of $\Mn_{k(v),\rm red}$}}
\bigskip
We recall that until the end we  assume that $G_{{\bf Z}_2}$ is isomorphic to ${\bf GSO}_{2n,{\bf Z}_2}^+$. Subsections 6.1 to 6.6 are intermediary steps towards the verification (see Proposition 6.7) that the reduced scheme $\Mn_{k(v),\rm red}$ of $\Mn_{k(v)}$ is regular and formally smooth over $k(v)$. We use the notations of Subsections 4.1 to 4.3. Let $\Mg_y$ be the reductive, closed subgroup scheme of ${\bf GL}_{M_y}$ introduced in the beginning of Subsection 5.2. 
\smallskip
In Subsections 6.1 and 6.2 we show the existence of a cocharacter $\mu_y:{\bf G}_{m,W(k)}\to\Mg_y$ that produces a direct sum decomposition $M_y=F^1\oplus F^0$ such that $F^1$ is the Hodge filtration defined by a lift $(\tilde A,\lambda_{\tilde A},\Mb)$ of $(A,\lambda_A,\Mb)$ to $W(k)$. Unfortunately, we can not check directly that we can  assume that $(\tilde A,\lambda_{\tilde A},\Mb)$ is the pull back of $(\Ma,\Lambda,\Mb)_{\Mn_{W(k)}}$ via a $W(k)$-valued point of $\Mn_{W(k)}$ that lifts $y$ (see Remark 6.2.1; this explains the length of this section). Subsection 6.3 associates to $\mu_y$ a smooth, closed subgroup scheme $\Mu_y$ of $\Mg_y$ of relative dimension $d$ over $W(k)$. Subsection 6.4 uses $\Mu_y$ and [14, Sect. 7] to construct a versal deformation of $(A,\lambda_A,\Mb)$ that defines naturally a morphism $m_{\Mr/2\Mr}:{\rm Spec}\,\Mr/2\Mr\to\Mm_k$ of $k$-schemes, where $\Mr$ is a $W(k)$-algebra of formal power series in $d$ variables. Subsections 6.5 to 6.7 use a modulo $2$ version of [14, Sect. 7] to check that $m_{\Mr/2\Mr}$ factors through a formally \'etale morphism $n_{\Mr/2\Mr}:{\rm Spec}\,\Mr/2\Mr\to\Mn_{k,\rm red}$ of $k$-schemes; the existence of $n_{\Mr/2\Mr}$ will surpass (in the geometric context of the Main Theorem) the modulo $2W(k)$ version of problems 1.1 (ii) and (iii). 
\medskip\smallskip\noindent
{\bf 6.1. Proposition.} {\it There exists a cocharacter $\mu_y:{\bf G}_{m,W(k)}\to\Mg_y$ that induces a direct sum decomposition $M_y=F^1\oplus F^0$ with the properties that $F^1$ lifts the Hodge filtration $\bar F^1_y:=F^1_V\otimes_V k$ of $M_y/2M_y$ defined by $A$ and that $\beta\in {\bf G}_m(W(k))$ acts on $F^i$ through $\mu_y$ as the multiplication with $\beta^{-i}$ ($i\in\{0,1\}$).}
\medskip\noindent
{\it Proof:}
Let $\tilde F^1_0$ be a direct summand of the $W(k)$-module $N_y$ that lifts the Hodge filtration of $\bar N_y=N_y/2N_y$. For $u\in\tilde F^1_0$ let $x:={{\Phi_y(u)}\over 2}\in N_y$. We have $b_{N_y}(x,x)={1\over 2}\sg(b_{N_y}(u,u))$ and (cf. Proposition 5.1) $b_{N_y}(x,x)\in 2W(k)$. Thus for $u\in\tilde F^1_0$ we have $b_{N_y}(u,u)\in 4W(k)$. Let $\{\tilde u_1,\ldots,\tilde u_n\}$ be a $W(k)$-basis for $\tilde F^1_0$. As $b_{N_y}$ is perfect, inductively for $l\in\{2,\ldots,n\}$ we can replace $\tilde u_l$ by an element in $\tilde u_l+2N_y$ such that we have $b_{N_y}(\tilde u_i,\tilde u_l)\in 4W(k)$ for all $i\in\{1,\ldots,l-1\}$. Thus we can choose $\tilde F^1_0$ such that we have $b_{N_y}(\tilde F^1_0,\tilde F^1_0)\subseteq 4W(k)$. 
\smallskip
We extend $\{\tilde u_1,\ldots,\tilde u_n\}$ to a $W(k)$-basis $\{\tilde u_1,\tilde v_1,\ldots,\tilde u_n,\tilde v_n\}$ for $N_y$ such that the following two things hold: for each $i\in\{1,\ldots,n\}$ we have $b_{N_y}(\tilde u_i,\tilde v_i)=1$ and for each $i$, $j\in\{1,\ldots,n\}$ with $i\neq j$ and for every pair $(u,v)\in\{(\tilde u_i,\tilde v_j),(\tilde v_i,\tilde v_j)\}$ we have $b_{N_y}(u,v)\in 4W(k)$ (one constructs such a $W(k)$-basis via the same inductive processes used in the previous paragraph and in the second paragraph of the proof of Proposition 3.4). 
\smallskip
As $b_{\bar N_y}$ is alternating (cf. Proposition 5.1) and as $b_{N_y}(\tilde u_i,\tilde u_j)\in 4W(k)$ for all $i$, $j\in\{1,\ldots,n\}$, from Proposition 3.4 (b) (applied with $S=W(k)$ and $q=2$) we get that there exists a $W(k)$-basis $\{u_1,v_1,\ldots,u_n,v_n\}$ for $N_y$ with respect to which the matrix of $b_{N_y}$ is $J(2n)$ and moreover for $i\in\{1,\ldots,n\}$ we have $u_i-\tilde u_i\in 2N_y$. Let $F^1_0$ and $F^0_0$ be the $W(k)$-spans of $\{u_1,\ldots,u_n\}$ and $\{v_1,\ldots,v_n\}$ (respectively). The direct sum decomposition $N_y=F^1_0\oplus F^0_0$ of $W(k)$-modules is such that $F^1_0$ and $\tilde F^1_0$ are congruent modulo $2W(k)$ and we have $b_{N_y}(F^1_0,F^1_0)=b_{N_y}(F^0_0,F^0_0)=0$. 
\smallskip
We refer to (7). As $M_y=N_y^s$, we can take the decomposition $M_y=F^1\oplus F^0$ such that we have $F^1:=(F^1_0)^s$ and $F^0:=(F^0_0)^s$. Thus we have:
\medskip
{\bf (i)} {\it the image of the cocharacter $\mu_y:{\bf G}_{m,W(k)}\to {\bf GL}_{M_y}$ that acts on $F^1$ and $F^0$ as desired, fixes all endomorphisms of $M_y$ defined by elements of $\Mb^{\rm opp}$.}
\medskip
As $b_{N_y}(F^1_0,F^1_0)=b_{N_y}(F^0_0,F^0_0)=0$ and as in Subsubsection 4.1.3 we have constructed $b_{\Me}$ in terms of $b_{\psi}$ and of ${\bf Z}_2$-endomorphisms of $\Ma$, we get also that: 
\medskip
{\bf (ii)} {\it the group scheme ${\bf G}_{m,W(k)}$ acts via $\mu_y$ on the $W(k)$-span of $b_{N_y}$ and thus also on the $W(k)$-span of $\lambda_{M_y}$, through the inverse of the identity character of ${\bf G}_{m,W(k)}$.}
\medskip
From properties (i) and (ii) we get that the cocharacter $\mu_y$ factors through the subgroup scheme of ${\bf GSp}(M_y,\lambda_{M_y})$ that fixes the ${\bf Z}_{(2)}$-algebra $\Mb^{\rm opp}$ of ${\rm End}(M_y)$. Therefore the cocharacter $\mu_y$ factors through $\Mg_y$.\endproof 
\medskip\noindent
{\bf 6.1.1. Corollary.} {\it The normalizer $\Mp_{y,k}$ of $\bar F^1_y=F^1\otimes_{W(k)} k$ in $\Mg_{y,k}:=\Mg_y\times_{W(k)} k$ is a parabolic subgroup of $\Mg_{y,k}$ such that we have $\dim(\Mg_{y,k}/\Mp_{y,k})={{n(n-1)}\over 2}=d$.}
\medskip\noindent
{\it Proof:} Let $q_{N_y}$ be as in Subsection 5.2 and let $q_{\bar N_y}$ be its reduction modulo $2W(k)$. With the notations of the previous proof, for $x\in F^1_0\cup F^0_0$ we have $q_{N_y}(x)=0$. As we can identify $\Mg_y$ with ${\bf GSO}(N_y,q_{N_y})$ and as $F^1=(F^1_0)^s$, we can also identify $\Mp_{y,k}$ with the normalizer of $F^1_0/2F^1_0$ in ${\bf GSO}(\bar N_y,q_{\bar N_y})$. Thus the corollary follows from Lemma 3.1.3.\endproof 
\medskip\smallskip\noindent
{\bf 6.2. Choice of $\mu_y$.} Let $(D_y,\lambda_{D_y},\Mb)$ be the principally quasi-polarized $2$-divisible group over $k$ endowed with endomorphisms of $y^*((\Ma,\Lambda,\Mb)_{\Mn_{W(k)}})$. We check that we can choose the cocharacter $\mu_y$ of Proposition 6.1 such that there exists a lift $(\tilde D,\lambda_{\tilde D},\Mb)$ of $(D_y,\lambda_{D_y},\Mb)$ to $W(k)$ with the property that $F^1$ is the Hodge filtration of $M_y$ defined by $\tilde D$. We consider the direct sum decomposition $M_y=M_{y(0)}\oplus M_{y(1)}\oplus M_{y(2)}$ left invariant by $\Phi_y$ and such that all slopes of $(M_{y(l)},\Phi_y)$ are $0$ if $l=0$, are $1$ if $l=2$, and belong to the interval $(0,1)$ if $l=1$. To it corresponds a product decomposition $D_y=\prod_{l=0}^2 D_{y(l)}$. 
\smallskip
We consider the Newton type of cocharacter $\nu_y:{\bf G}_{m,W(k)}\to {\bf GL}_{M_y}$ that acts on $M_{y(l)}$ via the $l$-th power of the identity character of ${\bf G}_{m,W(k)}$. The endomorphisms of $M_y$ defined by elements of $\Mb^{\rm opp}$ are fixed by (i.e., commute with) $\Phi_y$ and thus belong to $\oplus_{l=1}^2 {\rm End}(M_{y(l)})$. Therefore the cocharacter $\nu_y$ fixes all these endomorphisms. As we have $\lambda_{M_y}(\Phi_y(u),\Phi_y(v))=2\sg(\lambda_{M_y}(u,v))$ for all $u$, $v\in M_y$, we get that $\lambda_{M_y}(M_{y(0)},M_{y(0)})=\lambda_{M_y}(M_{y(2)},M_{y(2)})=\lambda_{M_y}(M_{y(0)}\oplus M_{y(2)},M_{y(1)})=0$ and thus the cocharacter $\nu_y$ also normalizes the $W(k)$-span of $\lambda_{M_y}$. We conclude that $\nu_y$ factors through $\Mg_y$.
\smallskip
The special fibre of $\nu_y$ normalizes $\bar F^1_y=F^1/2F^1$ (i.e., the kernel of $\Phi_y$ modulo $2W(k)$) and therefore it factors through $\Mp_{y,k}$. We can replace the role of $\mu_y$ by the one of an inner conjugate of it through an arbitrary element $h\in\Mg_y(W(k))$ that lifts an element of $\Mp_{y,k}(k)$. But there exist such elements $h$ with the property that $h\mu_yh^{-1}$ and $\nu_y$ commute. Thus not to introduce extra notations, we will assume that $\mu_y$ and $\nu_y$ commute. Therefore we have a direct sum decomposition $F^1=\oplus_{l=0}^2 F^1_{(l)}$, where $F^1_{(l)}:=F^1\cap M_{y(l)}$. There exists a unique $2$-divisible group $\tilde D_{(l)}$ over $W(k)$ that lifts $D_{y(l)}$ and whose Hodge filtration is $F^1_{(l)}$, cf. [16, 1.6 (ii) of p. 186] applied to the Honda triple $(M_{y(l)},{1\over 2}\Phi_y(F^1_{(l)}),\Phi_y)$. Let $\tilde D:=\prod_{l=0}^2 \tilde D_{(l)}$; the fact that there exists a lift $(\tilde D,\lambda_{\tilde D},\Mb)$ of $(D_y,\lambda_{D_y},\Mb)$ to $W(k)$ is also implied by loc. cit. 
\medskip\noindent
{\bf 6.2.1. Remark.} Let $(\tilde A,\lambda_{\tilde A},\Mb)$ be the principally polarized abelian scheme over $W(k)$ endowed with endomorphisms that lifts $(A,\lambda_A,\Mb)$ and whose principally quasi-polarized $2$-divisible group endowed with endomorphisms is $(\tilde D,\lambda_{\tilde D},\Mb)$, cf. Serre--Tate deformation theory and Grothendieck's algebraization theory (see [18, Ch. III, Thms. 5.4.1 and 5.4.5]). Let ${\rm Spec}\,W(k)\to\Mm_{W(k)}$ be the morphism of $W(k)$-schemes associated to $(\tilde A,\lambda_{\tilde A})$ and its symplectic similitude structures that lift those of $(A,\lambda_A)$. It is easy to check that it factors through the closed subscheme $\Mn^{\rm mod}_{W(k)}$ of $\Mm_{W(k)}$ (i.e., the determinant axiom mentioned in Subsubsection 4.1.1 holds for $\tilde A$ as it holds for $A_V$). But, as we are in the case (D), we have $G\neq G_1$ and the Hasse principle fails for $G$. Thus we can not prove directly that the last morphism of $W(k)$-schemes factors through the closed subscheme $\Mn_{W(k)}$ of $\Mn^{\rm mod}_{W(k)}$ (to be compared with [26, Sect. 8]). This explains why in Subsections 6.5 and 6.6 below we will work mainly modulo $2$ and not directly in mixed characteristic $(0,2)$. 
\medskip\smallskip\noindent
{\bf 6.3. Defining $\Mu_y$.} We have a natural direct sum decomposition of $W(k)$-modules ${\rm End}(M_y)={\rm End}(F^1)\oplus {\rm End}(F^0)\oplus {\rm Hom}(F^1,F^0)\oplus {\rm Hom}(F^0,F^1)$. Let $\Mu_y$ be the flat, closed subgroup scheme of ${\bf GL}_{M_y}$ defined by the rule: if $\tilde R$ is a commutative $W(k)$-algebra, then 
$$\Mu_y(\tilde R)=1_{M_y\otimes_{W(k)} \tilde R}+({\rm Lie}(\Mg_y)\cap {\rm Hom}(F^1,F^0))\otimes_{W(k)} \tilde R$$ 
(the last intersection being taken inside ${\rm End}(M_y)$). The group scheme $\Mu_y$ is smooth, connected, commutative, and its Lie algebra is the direct summand ${\rm Lie}(\Mg_y)\cap {\rm Hom}(F^1,F^0)$ of ${\rm Hom}(F^1,F^0)={\rm Hom}(F^1,M_y/F^1)$. Thus ${\rm Lie}(\Mu_y)\subseteq {\rm Lie}(\Mg_y)$. This implies $\Mu_{y,B(k)}\subset \Mg_{y,B(k)}$, cf. [6, Ch. II, Subsect. 7.1]. Thus $\Mu_y\subset \Mg_y$. The group scheme $\Mu_y$ acts trivially on both $F^0$ and $M_y/F^0$. Due to the existence of the cocharacter $\mu_y:{\bf G}_{m,W(k)}\to\Mg_y$, we have a direct sum decomposition ${\rm Lie}(\Mg_{y,k})={\rm Lie}(\Mu_{y,k})\oplus {\rm Lie}(\Mp_{y,k})$ of $k$-vector spaces. From the last two sentences we get that the group scheme $\Mu_{y,k}\cap\Mp_{y,k}$ is trivial and that the natural morphism $\Mu_{y,k}\to\Mg_{y,k}/\Mp_{y,k}$ of $k$-schemes is an open embedding. Thus the relative dimension of $\Mu_y$ over $W(k)$ is $d=\dim(\Mg_{y,k}/\Mp_{y,k})$, cf. Corollary 6.1.1. 
\smallskip
Let ${\rm Spf}\,\Mr$ be the completion of $\Mu_y$ along its identity section. We can identify $\Mr=W(k)[[t_1,\ldots,t_d]]$ in such a way that the ideal $(t_1,\ldots,t_d)$ defines the identity section of $\Mu_y$. 
Let $\Phi_{\Mr}$ be the Frobenius lift of $\Mr$ that is compatible with $\sg$ and takes $t_i$ to $t_i^2$ for all $i\in\{1,\ldots,d\}$. Let $d\Phi_{\Mr}:\oplus_{i=1}^d\Mr dt_i\to\oplus_{i=1}^d\Mr dt_i$ be the $(2,t_1,\ldots,t_d)$-completion of the differential map of $\Phi_{\Mr}$. We consider the universal element
$$u_{\rm univ}\in \Mu_y(\Mr)$$ 
defined by the natural morphism ${\rm Spec}\,\Mr\to \Mu_y$ of $W(k)$-schemes. 
\smallskip
The existence of $\tilde D$ allows us to apply [14, Sect. 7, Thm. 10]: there exists a unique connection $\nabla_y:M_y\otimes_{W(k)} \Mr\to M_y\otimes_{W(k)} \oplus_{i=1}^d \Mr dt_i$ such that we have an identity
$$\nabla_y\circ u_{\rm univ}(\Phi_y\otimes\Phi_{\Mr})=(u_{\rm univ}(\Phi_y\otimes\Phi_{\Mr})\otimes d\Phi_{\Mr})\circ \nabla_y\leqno (13)$$ 
of maps from $M_y\otimes_{W(k)} \Mr$ to $M_y\otimes_{W(k)} \oplus_{i=1}^d \Mr dt_i$; the connection $\nabla_y$ is integrable and topologically nilpotent. We emphasize that loc. cit. pertains to all primes including $2$ (primes at least $3$ are used in [14] only until [14, Sect. 7], cf. [14, Sect. 1]). Let $\nabla_{\rm triv}$ be the flat connection on $M_y\otimes_{W(k)} \Mr$ that annihilates $M_y\otimes 1$. Due to Formula (13), it is easy to check that the connections $\nabla_{\rm triv}+u_{\rm univ}^{-1}du_{\rm univ}$ and $\nabla_y$ are congruent modulo $(t_1,\ldots,t_d)$. From this, the fact that ${\rm Lie}(\Mu_y)$ is a direct summand of ${\rm Hom}(F^1,M_y/F^1)$, and the very definition of $u_{\rm univ}\in \Mu_y(\Mr)$ we get that:
\medskip\noindent
{\bf 6.3.1. Fact.} {\it  The Kodaira--Spencer map of $\nabla_y$ is injective modulo $(2,t_1,\ldots,t_d)$ and thus its image is a direct summand ${\got I}_y$ of ${\rm Hom}(F^1,M_y/F^1)\otimes_{W(k)} \Mr$ of rank $d$.} 
\medskip\smallskip\noindent
{\bf 6.4. A deformation.} We recall that the categories of $2$-divisible groups over ${\rm Spf}\,\Mr/2\Mr$ and respectively over ${\rm Spec}\,\Mr/2\Mr$, are canonically isomorphic (cf. [29, Ch. II, Lem. 4.16]). From this and [14, Sect. 7, Thm. 10] we get the existence of a $2$-divisible group $\tilde D_{\Mr/2\Mr}$ over $\Mr/2\Mr$ whose $F$-crystal over $\Mr/2\Mr$ is $(M_y\otimes_{W(k)} \Mr,u_{\rm univ}(\Phi_y\otimes\Phi_{\Mr}),\nabla_y)$. We have 
$$u_{\rm univ}(\Phi_y\otimes\Phi_{\Mr})\circ b=b\circ u_{\rm univ}(\Phi_y\otimes\Phi_{\Mr})\;\; \forall b\in\Mb^{\rm opp}.\leqno (14)$$ 
As $\nabla_y$ is uniquely determined by $u_{\rm univ}(\Phi_y\otimes\Phi_{\Mr})$ and due to Formula (14), $\nabla_y$ is invariant under the natural action of the group of invertible elements of $\Mb^{\rm opp}$ on $M_y\otimes_{W(k)} \Mr$. Each element of $\Mb^{\rm opp}$ is a sum of two invertible elements of $\Mb^{\rm opp}$. Thus as the ring $\Mr/2\Mr$ has a finite $2$-basis, from the fully faithfulness part of [5, Thm. 4.1.1 or 4.2.1] we get first that the $2$-divisible group $\tilde D_{\Mr/2\Mr}$ is uniquely determined (by its $F$-crystal over $\Mr/2\Mr$) and second that each $b\in \Mb^{\rm opp}$ is the crystalline realization of a unique endomorphism $b$ of $\tilde D_{\Mr/2\Mr}$. Thus we have a natural ${\bf Z}_{(2)}$-monomorphism $\Mb\hookrightarrow {\rm End}(\tilde D_{\Mr/2\Mr})$. A similar argument shows that there exists a unique principal quasi-polarization $\lambda_{\tilde D_{\Mr/2\Mr}}$ of $\tilde D_{\Mr/2\Mr}$ such that ${\bf D}(\lambda_{\tilde D_{\Mr/2\Mr}})$ is defined by the alternating form $\lambda_{M_y}$ on $M_y\otimes_{W(k)} \Mr$. As ${\rm Lie}(\Mg_y)$ is the direct summand of ${\rm Lie}({\bf GSp}(M_y,\lambda_{M_y}))$ that annihilates ${\rm Im}(\Mb^{\rm opp}\to {\rm End}(M_y))$ and as $\nabla_y$ annihilates ${\rm Im}(\Mb^{\rm opp}\to {\rm End}(M_y)\otimes_{W(k)} \Mr)$, we have 
$$\nabla_y-\nabla_{\rm triv}\in {\rm Lie}(\Mg_y)\otimes_{W(k)} \oplus_{i=1}^d \Mr dt_i\subseteq {\rm End}(M_y)\otimes_{W(k)} \oplus_{i=1}^d \Mr dt_i.$$
From this, as inside ${\rm End}(M_y)$ we have ${\rm Lie}(\Mu_y)={\rm Lie}(\Mg_y)\cap {\rm Hom}(F^1,M_y/F^1)$ (cf. definition of $\Mu_y$), we get that the direct summand ${\got I}_y$ of ${\rm Hom}(F^1,M_y/F^1)\otimes_{W(k)} \Mr$ introduced by Fact 6.3.1 is contained in the direct summand ${\rm Lie}(\Mu_y)\otimes_{W(k)} \Mr$ of ${\rm Hom}(F^1,M_y/F^1)\otimes_{W(k)} \Mr$ and thus by reasons of ranks it is  exactly ${\rm Lie}(\Mu_y)\otimes_{W(k)} \Mr$. Let
$$\Mc_y:=(M_y\otimes_{W(k)} \Mr,F^1\otimes_{W(k)} \Mr,u_{\rm univ}(\Phi_y\otimes\Phi_{\Mr}),\nabla_y,\lambda_{M_y},\Mb^{\rm opp}).$$ 
\indent
From Serre--Tate deformation theory we get the existence of a unique triple 
$$(A^\prime_{\Mr/2\Mr},\lambda_{A^\prime_{\Mr/2\Mr}},\Mb)$$ 
over $\Mr/2\Mr$ which lifts $(A,\lambda_A,\Mb)$ and whose principally quasi-polarized $2$-divisible group endowed with endomorphisms is $(\tilde D_{\Mr/2\Mr},\lambda_{\tilde D_{\Mr/2\Mr}},\Mb)$. 
To $(A^\prime_{\Mr/2\Mr},\lambda_{A^\prime_{\Mr/2\Mr}})$ and its symplectic similitude structures lifting those of $(A,\lambda_A)$, corresponds a morphism of $k$-schemes
$$m_{\Mr/2\Mr}:{\rm Spec}\,\Mr/2\Mr\to\Mm_k.$$
\noindent
{\bf 6.5. Specializing to $y$.} Let $R=W(k)[[t]]$ be as in  Subsection 2.1 and let ${\got R}$ be its field of fractions. We consider an arbitrary morphism 
$$n_{R/2R}:{\rm Spec}\,R/2R\to\Mn_k$$
of $k$-schemes that lifts the point $y\in\Mn_{W(k)}(k)=\Mn_k(k)$. Let $m_{R/2R}:{\rm Spec}\,R/2R\to\Mm_k$ be the composite of $n_{R/2R}$ with the morphism $\Mn_k\to\Mm_k$. Let $k_1$ be an algebraic closure of the field of  fractions $k((t))$ of $R/2R=k[[t]]$. Let the morphism $y_1:{\rm Spec}\,k_1\to\Mn_{W(k_1)}$ be defined naturally by $n_{R/2R}$. We recall from  Subsection 2.1 that $U_m=k[[t]]/(t^m)$, where $m\in{\bf N}^{\ast}$. Let 
$$n_{U_m}:{\rm Spec}\,U_m\to\Mn_k\;\; {\rm and}\;\; m_{U_m}:{\rm Spec}\,U_m\to\Mm_k$$
be the composites of the closed embedding ${\rm Spec}\,U_m\hookrightarrow {\rm Spec}\,R/2R$ with $n_{R/2R}$ and $m_{R/2R}$ (respectively). Let $(M_R,\Phi_{M_R},\nabla_{M_R},\lambda_{M_R})$ be the principally quasi-polarized $F$-crystal over $R/2R$ of the pull back $(\tilde A_{R/2R},\lambda_{\tilde A_{R/2R}})$ through $m_{R/2R}$ of the universal principally polarized abelian scheme over $\Mm_k$; thus $(M_R,\lambda_{M_R})$ is a symplectic space over $R$ of rank $\dim_{{\bf Q}}(W)$, etc. We have a natural ${\bf Z}_2$-monomorphism $\Mb^{\rm opp}\hookrightarrow {\rm End}(M_R,\Phi_{M_R},\nabla_{M_R})$. Let $F^1_{R/2R}$ be the Hodge filtration of $M_R/2M_R$ defined by $\tilde A_{R/2R}$ (i.e., the kernel of the reduction of $\Phi_{M_R}$ modulo $2R$). Let $\Mg_R^\prime$ be the schematic closure in ${\bf GL}_{M_R}$ of the identity component of the subgroup of ${\bf GSp}(M_R,\lambda_{M_R})_{{\got R}}$ that fixes the ${\got R}$-subalgebra $\Mb^{\rm opp}\otimes_{{\bf Z}_{(2)}} {\got R}$ of ${\rm End}(M_R\otimes_{R} {\got R})$. Let $M_R=N_{R}^s$ be the decomposition into $R$-modules that corresponds naturally to (4). Let $b_{N_R}$ be the perfect bilinear form on $N_{R}$ that corresponds naturally to the (pull back to ${\rm Spec}\, R/2R$ of the) isomorphism $b_{\Me}$ of Subsubsection 4.1.3; it is symmetric (cf. Subsection 4.3 applied to $y_1$). As $b_{N_R}$ modulo $2R$ is alternating (cf. Proposition 5.1 applied to $y_1$), the formula $q_{N_R}(x):={{b_{N_R}(x,x)}\over 2}$ defines a quadratic form on $N_R$. We can redefine $\Mg_R^\prime={\bf GSO}(N_R,q_{N_R})$. Thus $\Mg^\prime_R$ is a reductive, closed subgroup scheme of ${\bf GL}_{M_R}$, cf. Proposition 3.4 (c). 
\smallskip
Let $\Mp^\prime_{k((t))}$ be the parabolic subgroup of $\Mg^\prime_{k((t))}$ that is the normalizer of $F^1_{R/2R}\otimes_{R/2R} k((t))$ in $\Mg^\prime_{k((t))}$, cf. Corollary 6.1.1 applied to $y_1$. The $R$-scheme of parabolic subgroup schemes of $\Mg^\prime_{R}$ is projective, cf. [12, Vol. III, Exp. XXVI, Cor. 3.5]. Therefore the schematic closure $\Mp^\prime_{R/2R}$ of $\Mp^\prime_{k((t))}$ in $\Mg^\prime_{R/2R}$ is a parabolic subgroup scheme  of $\Mg^\prime_{R/2R}$ that  normalizes $F^1_{R/2R}$ and that (cf. Corollary 6.1.1 applied to $y_1$) has the same relative dimension as $\Mp_{y,k}$. Thus we have a monomorphism $\Mp^\prime_k\hookrightarrow\Mp_{y,k}$ over $k$ which by reasons of dimensions is an isomorphism over $k$, to be viewed as an identification. The special fibre $\bar\mu_y$ of $\mu_y$ factors through $\Mp_{y,k}=\Mp^\prime_k$. From [12, Vol. II, Exp. IX, Thms. 3.6 and 7.1] we get that $\bar\mu_y$ lifts to a cocharacter $\mu_{R/2R}:{\bf G}_{m,R/2R}\to\Mp^\prime_{R/2R}$ which at its turn lifts to a cocharacter $\mu_R:{\bf G}_{m,R}\to \Mg^\prime_R$. Let $M_R=F^{1\prime}_{R}\oplus F^{0\prime}_{R}$ be the direct sum decomposition such that $\beta\in {\bf G}_m(R)$ acts on $F^{i\prime}_R$ through $\mu_y$ as the multiplication with $\beta^{-i}$ ($i\in\{0,1\}$). From constructions we get that 
$$F^{1\prime}_R/2F^{1\prime}_R=F^1_{R/2R}\;\; {\rm and}\;\;F^{0\prime}_R\otimes_R k=F^0/2F^0.$$ 
\noindent
{\bf 6.6. Theorem.} {\it For $m\in{\bf N}^{\ast}$ let $\Mt_m:=n_{U_m}^*((\Ma,\Lambda,\Mb)_{\Mn_k})$. There exists a morphism $u_m:{\rm Spec}\,U_m\to {\rm Spec}\,\Mr/2\Mr$ of $k$-schemes such that $u_m^*(A^\prime_{\Mr/2\Mr},\lambda_{A^\prime_{\Mr/2\Mr}},\Mb)$ is isomorphic to $\Mt_m$, under an isomorphism that lifts the identity automorphism of $(A,\lambda_A,\Mb)$.}  
\medskip\noindent
{\it Proof:} We consider the reduction
$$(M_R(m),F^{1\prime}_R(m),\Phi_{M_R}(m),\nabla_{M_R}(m),\lambda_{M_R}(m),\Mg_R^\prime(m))$$ 
of $(M_R,F^{1\prime}_R,\Phi_{M_R},\nabla_{M_R},\lambda_{M_R},\Mg^\prime_R)$ modulo the ideal $(t^m)$ of $R$. This theorem is a geometric variant of a slight modification of [14, Sect. 7, Thm. 10 and Rm. iii)]. Following loc. cit., we show by induction on $m\in{\bf N}^{\ast}$ that the morphism $u_m$ exists and that there exists a $W(k)$-homomorphism $l(m):\Mr\to R/(t^m)$ which has the following two properties:
\medskip
{\bf (i)} it maps $(t_1,\ldots,t_d)$ to $(t)/(t^m)$ and modulo $2W(k)$ it defines the morphism $u_m$; 
\smallskip
{\bf (ii)} the extension of $\Mc_y$ via $l(m)$ is isomorphic to $(M_R(m),F^{1\prime}_R(m),\Phi_{M_R}(m),\nabla_{M_R}(m),\break \lambda_{M_R}(m),\Mb^{\rm opp})$ under an isomorphism $\My(m)$ that modulo $(t)/(t^m)$ is $1_{M_y}$. 
\medskip
The case $m=1$ follows from constructions. For $m\ge 2$ the passage from $m-1$ to $m$ goes as follows. We endow the ideal $J_m:=(t^{m-1})/(t^m)$ of $R/(t^m)$ with the natural divided power structure; thus $J_m^{[2]}=0$. The ideal $(2,J_m)$ of $R/(t^m)$ has a natural divided power structure that extends the one of $J_m$. Let $\tilde l(m):\Mr\to R/(t^m)$ be a $W(k)$-homomorphism that lifts $l(m-1)$. Let $\tilde u_m:{\rm Spec}\,U_m\to {\rm Spec}\,\Mr/2\Mr$ be defined by $\tilde l(m)$ modulo $2W(k)$. We apply the crystalline Dieudonn\'e functor ${\bf D}$:
\medskip
--  in the context of the principally quasi-polarized $2$-divisible groups endowed with endomorphisms of $\tilde u_m^*(A^\prime_{\Mr/2\Mr},\lambda_{A^\prime_{\Mr/2\Mr}},\Mb)$ and of $\Mt_m$, and 
\smallskip
--  in the context of the thickenings attached naturally to the closed embeddings ${\rm Spec}\,U_{m-1}\hookrightarrow {\rm Spec}\,R/(t^m)$ and ${\rm Spec}\,U_m\hookrightarrow {\rm Spec}\,R/(t^m)$.
\medskip
We get that the extension of $\Mc_y$ through $\tilde l(m)$ is isomorphic to the sextuple 
$$(M_R(m),\tilde F^1_R(m),\Phi_{M_R}(m),\nabla_{M_R}(m),\lambda_{M_R}(m),\Mb^{\rm opp})$$ 
under an isomorphism which is denoted by $\My(m)$ and which modulo $(t)/(t^m)$ is $1_{M_y}$. As $\Phi_{\Mr}(t_1,\ldots,t_d)\subseteq (t_1,\ldots,t_d)^2$, the isomorphism $\My(m-1)$ is uniquely determined by the property (ii). Thus $\My(m)$ lifts $\My(m-1)$ and therefore $\tilde F^1_R(m)$ is a direct summand of $M_R(m)$ that lifts $F^{1\prime}_{R}(m-1)$. 
\smallskip
We check that under $\My(m)$, the reductive subgroup scheme $\Mg_y\times_{W(k)} \Mr$ of ${\bf GL}_{M_y\otimes_{W(k)} \Mr}$ pulls back to the reductive subgroup scheme $\Mg_{R}^\prime(m)$ of ${\bf GL}_{M_R(m)}$. For this it suffices to check that under the isomorphism $\My(m)[{1\over 2}]$, the reductive subgroup scheme $\Mg_y\times_{W(k)} \Mr[{1\over 2}]$ of ${\bf GL}_{M_y\otimes_{W(k)} \Mr[{1\over 2}]}$ pulls back to the reductive subgroup scheme $\Mg_R^\prime(m)[{1\over 2}]$ of ${\bf GL}_{M_R(m)[{1\over 2}]}$. But this holds as $\Mg_y\times_{W(k)} \Mr[{1\over 2}]$ (resp. as $\Mg_R^\prime\times_R R[{1\over 2}]$) is the identity component of the subgroup scheme of ${\bf GSp}(M_y\otimes_{W(k)} \Mr[{1\over 2}],\lambda_{M_y})$ (resp. of ${\bf GSp}(M_R[{1\over 2}],\lambda_{M_R})$) that fixes all elements of $\Mb^{\rm opp}$. For $\tilde m\in\{m-1,m\}$ we will identify naturally $\Mu_y\times_{W(k)} R/(t^{\tilde m})$ with a closed subgroup scheme of $\Mg_R^\prime(\tilde m)$ and thus we will view $\Mu_y(R/(t^{\tilde m}))$ as a subgroup of $\Mg_R^\prime (R/(t^{\tilde m}))=\Mg_R^\prime(\tilde m)(R/(t^{\tilde m}))$.
\smallskip
Let $\mu_{R,m}:{\bf G}_{m,R/(t^m)}\to\Mg_R^\prime(m)$ be the reduction of the cocharacter $\mu_R$ modulo $(t^m)$. Let $\mu_{y,m}:{\bf G}_{m,R/(t^m)}\to\Mg_R^\prime(m)$ be the cocharacter obtained naturally from $\mu_{y,\Mr}$ via $\tilde l(m)$ and $\My(m)$. As over $W(k)$ the two cocharacters $\mu_{R,m}$ and $\mu_{y,m}$ of $\Mg_R^\prime(m)$ coincide, there exists an element $h_3\in{\rm Ker}(\Mg_R^\prime(R/(t^m))\to\Mg_R^\prime(R/(t)))$ such that we have an identity
$$h_3\mu_{R,m}h_3^{-1}=\mu_{y,m}$$ 
(cf. [12, Vol. II, Exp. IX, Thm. 3.6]). As $\tilde F^1_R(m)$ lifts $F^{1\prime}_R(m-1)$, we can write $h_3=h_1h_2$, where $h_1\in {\rm Ker}(\Mu_y(R/(t^m))\to \Mu_y(R/(t^{m-1})))$ and where $h_2\in {\rm Ker}(\Mg_R^\prime(R/(t^m))\to \Mg_R^\prime(R/(t)))$ normalizes $F^{1\prime}_R(m)$. We have $h_1(F^{1\prime}_R(m))=h_3(F^{1\prime}_R(m))=\tilde F^1_R(m)$.
\smallskip
As $h_1\in {\rm Ker}(\Mu_y(R/(t^m))\to \Mu_y(R/(t^{m-1})))$, as the trivial divided power structure on $J_m$ is nilpotent, and (see Subsection 6.4) as we have an identity ${\got I}_y={\rm Lie}(\Mu_y)\otimes_{W(k)} \Mr\subseteq {\rm Hom}(F^1,M_y/F^1)\otimes_{W(k)} \Mr$, we can replace $\tilde l(m)$ by another lift $l(m)$ of $l(m-1)$ such that its corresponding $h_1$ element is (i.e., becomes) the identity element. Let $u_m:{\rm Spec}\,U_m\to {\rm Spec}\,\Mr/2\Mr$ be defined by $l(m)$ modulo $2W(k)$. 
\smallskip
Thus we have $F^{1\prime}_R(m)=\tilde F^1_R(m)$. As the closed embedding ${\rm Spec}\,U_{m-1}\hookrightarrow {\rm Spec}\,U_m$ is a nilpotent thickening, from the identity $F^{1\prime}_R(m)=\tilde F^1_R(m)$ and from the Grothendieck--Messing deformation theory we get that the principally quasi-polarized $2$-divisible groups endowed with endomorphisms associated to $u_m^*(A^\prime_{\Mr/2\Mr},\lambda_{A^\prime_{\Mr/2\Mr}},\Mb)$ and $\Mt(m)$ are the same lift to ${\rm Spec}\,U_m$ of the principally quasi-polarized $2$-divisible group endowed with endomorphisms of $u_{m-1}^*(A^\prime_{\Mr/2\Mr},\lambda_{A^\prime_{\Mr/2\Mr}},\Mb)=\Mt(m-1)$ (the last identification is via the isomorphism whose existence is guaranteed by the inductive assumption, whose crystalline realization is $\My(m-1)$, and which lifts the identity automorphism of $(A,\lambda_A,\Mb)$). From this and the Serre--Tate deformation theory we get that $u_m^*(A^\prime_{\Mr/2\Mr},\lambda_{A^\prime_{\Mr/2\Mr}},\Mb)$ and $\Mt(m)$ are canonically identified via an isomorphism whose crystalline realization is $\My(m)$ and which lifts the identity automorphism of $(A,\lambda_A,\Mb)$. This ends the induction.\endproof 
\medskip
 From Theorem 6.6 we get directly that for each $m\in{\bf N}^{\ast}$, the point $m_{U_m}\in\Mm_k(U_m)$ is $m_{\Mr/2\Mr}\circ u_m\in\Mm_k(U_m)$. This implies that:
\medskip\noindent
{\bf 6.6.1. Corollary.} {\it Each morphism $m_{R/2R}:{\rm Spec}\,R/2R\to\Mm_k$ that factors through $\Mn_k$ in such a way that it lifts $y\in\Mn_k(k)$, factors also through $m_{\Mr/2\Mr}:{\rm Spec}\,\Mr/2\Mr\to\Mm_k$.} 
\medskip\noindent
{\bf 6.6.2. Fact.} {\it Let $S$ be a reduced $k$-algebra which is a local, complete, noetherian ring of residue field $k$ (thus $S$ is a reduced quotient of $k[[t_1,\ldots,t_s]]$ for some $s\in {\bf N}^{\ast}$). Then there exists a family $\got F$ of $k$-homomorphisms $S\to k[[t]]=R/2R$ such that $\cap_{h\in {\got F}} {\rm Ker}(h)=0$ (i.e., $S$ is naturally a $k$-subalgebra of $\prod_{h\in {\got F}} k[[t]]$).}
\medskip\noindent
{\it Proof:} By replacing $S$ with a quotient of $S$ by a minimal prime ideal of it, we can assume that $S$ is integral. We proceed by induction on $v=\dim(S)$. The case $v=0$ is trivial. If $v=1$, then the normalization of $S$ is isomorphic to $k[[t]]$ and thus we can take ${\got F}$ to be formed by one inclusion $S\hookrightarrow k[[t]]$. For $v\ge 2$, the passage from $v-1$ to $v$ goes as follows. Let $\Lambda_1(S)$ be the set of prime ideals of $S$ of height $1$. For ${\got p}\in \Lambda_1(S)$ let ${\got F}_{\got p}$ be the set of all $k$-homomorphisms $S\to k[[t]]$ whose kernels contain ${\got p}$. As $\dim(S/{\got p})\le v-1$, by induction we get that ${\got p}=\cap_{h\in {\got F}_{{\got p}}} {\rm Ker}(h)$. For each non-zero element $a$ of the maximal ideal ${\got m}_S$ of $S$, the number of prime ideals ${\got p}\in \Lambda_1(S)$ containing $a$ is a positive integer. As $v\ge 2$, ${\got m}_S$ is not a finite union of elements of $\Lambda_1(S)$ and thus the set $\Lambda_1(S)$ is infinite. From the last two sentences we get that $\cap_{{\got p}\in\Lambda_1(S)} {\got p}=0$. If ${\got F}:=\cup_{{\got p}\in \Lambda_1(S)} {\got F}_{{\got p}}$, then we have $\cap_{h\in {\got F}} {\rm Ker}(h)=\cap_{{\got p}\in \Lambda_1(S)} \cap_{h\in {\got F}_{{\got p}}} {\rm Ker}(h)=\cap_{{\got p}\in \Lambda_1(S)} {\got p}=0$.\endproof
\medskip\smallskip\noindent
{\bf 6.7. Proposition.} {\it The scheme $\Mn_{k(v),\rm red}$ is regular and formally smooth over $k(v)$.}
\medskip\noindent
{\it Proof:} We recall that $y\in\Mn_{W(k)}(k)$ is arbitrary. Let $\bar\Mo_y$ (resp. $\bar\Mo_y^{\rm big}$) be the completion of the local ring of the $k$-valued point of $\Mn_{k,\rm red}:=\Mn_{k(v),{\rm red}}\times_{k(v)} k$ (resp. of $\Mm_k$) defined by $y$. As $\Mn_k$ is a pro-\'etale cover of $\Mn_k/H_0$, $\bar\Mo_y$ is a local, complete, excellent $k$-algebra of dimension $d=\dim(\Mr/2\Mr)$. To prove the proposition, it suffices to show that $\bar\Mo_y$ is isomorphic to $\Mr/2\Mr$. 
\smallskip
The $k$-homomorphism $\bar m_y:\bar\Mo_y^{\rm big}\to \Mr/2\Mr$ defined naturally by $m_{\Mr/2\Mr}$ is onto, cf, Fact 6.3.1. From Corollary 6.6.1 applied to morphisms $n_{R/2R}:{\rm Spec}\,R/2R\to\Mn_k$ of $k$-schemes that factor through ${\rm Spec}\,\bar\Mo_y$ and from Fact 6.6.2 applied with $S=\bar\Mo_y$, we get that ${\rm Ker}(\bar m_y)$ has a trivial image in the $\bar\Mo_y^{\rm big}$-algebra $\bar\Mo_y$. Thus the $k$-epimorphism $\bar\Mo_y^{\rm big}\twoheadrightarrow \bar\Mo_y$ factors through $\bar m_y$ inducing a $k$-epimorphism $\bar q_y:\Mr/2\Mr\twoheadrightarrow \bar\Mo_y$. By reasons of dimensions, $\bar q_y$ is an isomorphism. Thus $\bar\Mo_y$ is isomorphic to $\Mr/2\Mr$ and in fact $m_{\Mr/2\Mr}:{\rm Spec}\,\Mr/2\Mr\to \Mm_k$ factors uniquely through a morphism 
$$n_{\Mr/2\Mr}:{\rm Spec}\,\Mr/2\Mr\to \Mn_{k,{\rm red}}\leqno (15)$$
that lifts $y\in\Mn_{W(k)}(k)=\Mn_{k,{\rm red}}(k)$. \endproof 
\bigskip\smallskip\noindent
{\bigsll {\bf 7. Proof of Theorem 1.4 (b)}}
\bigskip
The following key proposition is a natural application of [35] and Lemma 3.1.4. 
\medskip\smallskip\noindent
{\bf 7.1. Proposition.} {\it Suppose $y\in\Mn_{W(k)}(k)$ factors through the ordinary locus of $\Mn_k$ (i.e., the abelian variety $A$ is ordinary). Then $\Mn_{W(k)}$ is regular and formally smooth over $W(k)$ at $y$.}
\medskip\noindent
{\it Proof:} Let $\Mo_y$ (resp. $\Mo_y^{\rm big}$) be the completion of the local ring of the $k$-valued point of $\Mn_{W(k)}$ (resp. of $\Mm_{W(k)}$) defined by $y$. As $\Mn_{W(k)}$ is a pro-\'etale cover of $\Mn_{W(k)}/H_0$, $\Mo_y$ is a local, complete, excellent $W(k)$-algebra. To prove the proposition, it suffices to show that $\Mo_y$ is isomorphic to $\Mr$. Let $\Mo_y^{\rm n}$ be the normalization of $\Mo_y$ in the ring of fractions of $\Mo_y$. As $\Mo_y$ is reduced and excellent, there exists $v\in {\bf N}^{\ast}$ such that we have a finite product decomposition $\Mo_y^{\rm n}=\prod_{i=1}^v \Mo_{y,i}$ into local, complete, normal rings which are finite $\Mo_y$-modules. With the notations of the proof of Proposition 6.7, one has $\bar\Mo_y=(\Mo_y/2\Mo_y)_{\rm red}$ and $\bar\Mo_y^{\rm big}=\Mo_y^{\rm big}/2\Mo_y^{\rm big}$. 
\smallskip
As $A$ is ordinary, it is well known that ${\rm Spf}\,\Mo_y^{\rm big}$ (i.e., the formal deformation space of $(A,\lambda_A)$) over ${\rm Spf}\,W(k)$ has a canonical structure of a formal torus. From [35, Cor. 3.8] we get that $\Mo_y^{\rm n}$ is regular and formally smooth over $W(k)$ and thus that in fact each local $W(k)$-algebra $\Mo_{y,i}$ is isomorphic to $\Mr$. From [35, Thm. 3.7 and Cor. 3.8] we get that in fact the $W(k)$-homomorphism $\Mo_y\to \Mo_{y,i}$ is surjective and that ${\rm Spf}\,\Mo_{y,i}$ is a closed formal subscheme of ${\rm Spf}\,\Mo_y^{\rm big}$ which is the translation of a formal subtorus of ${\rm Spf}\,\Mo_y^{\rm big}$ by a torsion ${\rm Spf}\,W(k)$-valued point. As a formal torus over ${\rm Spf}\,W(k)$ has no torsion ${\rm Spf}\,W(k)$-valued  point which is not a $2$-torsion ${\rm Spf}\,W(k)$-valued point, we get that ${\rm Spf}\,\Mo_{y,i}$ is a closed formal subscheme of ${\rm Spf}\,\Mo_y^{\rm big}$ which is the translation of a formal subtorus of ${\rm Spf}\,\Mo_y^{\rm big}$ by a $2$-torsion ${\rm Spf}\,W(k)$-valued point. This implies that the number of ${\rm Spf}\,W(k)$-valued points of ${\rm Spf}\,\Mo_{y,i}$ which define $2$-torsion ${\rm Spf}\,W(k)$-valued points of ${\rm Spf}\,\Mo_y^{\rm big}$ (i.e., define quasi-canonical lifts of $y$) is exactly $2^d$ and it does not depend on $i$. Thus, as $\Mo_y[{1\over 2}]=\Mo_y^{\rm n}[{1\over 2}]$, the number of ${\rm Spf}\,W(k)$-valued points of ${\rm Spf}\,\Mo_y$ which define $2$-torsion ${\rm Spf}\,W(k)$-valued points of ${\rm Spf}\,\Mo_y^{\rm big}$ is exactly $v2^d$. If $v=1$, then the $W(k)$-homomorphism $\Mo_y\to \Mo_{y,i}=\Mo_y^{\rm n}$ is surjective as well as injective and thus it is an isomorphism. Thus if $v=1$, then $\Mo_y$ is isomorphic to $\Mr$. Therefore to complete the proof it suffices to show that $v=1$. 
\smallskip
Let $V=W(k)$ and let $z\in\Mn_{W(k)}(W(k))$ be a lift of $y$ which defines a $2$-torsion ${\rm Spf}\,W(k)$-valued point of ${\rm Spf}\,\Mo_y^{\rm big}$. Thus $z\in\Mn_{W(k)}(W(k))$ is a quasi-canonical lift of $y\in\Mn_{W(k)}(k)$; the number of such quasi-canonical lifts $z$ is exactly $v2^d$. Let $A_V=A_{W(k)}$ and $b_z:\Me_z\tilde\to \Me_z^{\rm t}$ be as in Subsection 4.3. Let $A_{W(k)}^{\rm can}$ be the canonical lift of $A$ to $W(k)$. Each ${\bf Z}_{2}$-endomorphism of $A$ (or isomorphism $A\to A^{\rm t}$) lifts to a ${\bf Z}_{2}$-endomorphism of $A_{W(k)}^{\rm can}$ (or to an isomorphism $A_{W(k)}^{\rm can}\to (A_{W(k)}^{{\rm can}})^{\rm t}$). This implies that the $2$-divisible group of $A_{W(k)}^{\rm can}$ is an $s$-th power $\Me^s_{\rm can}$ and moreover we have an isomorphism $b_{\rm can}:\Me_{\rm can}\tilde\to \Me_{\rm can}^{\rm t}$ whose crystalline realization is symmetric and whose special fibre $b_{\rm can,k}:\Me_{\rm can,k}\tilde\to \Me_{\rm can,k}^{\rm t}$ is $b_{z,k}:\Me_{z,k}\tilde\to \Me_{z,k}^{\rm t}$ (to be compared with the definitions of $\Me_z$ and $b_z$). As the $2$-divisible group of $A_{W(k)}^{\rm can}$ is a direct sum of an \'etale $2$-divisible group and a $2$-divisible group of multiplicative type, we have a direct sum decomposition $\Me_{\rm can}=\Me_{\rm can,0}\oplus \Me_{\rm can,1}$ into $2$-divisible groups which are isotropic with respect to $b_{\rm can}$, with $\Me_{\rm can,0}$ \'etale and with $\Me_{\rm can,1}$ of multiplicative type. This implies that $b_{\rm can}[2]$ is a principal quasi-polarization. 
\smallskip
Thus to $\Me_{\rm can}=\Me_{\rm can,0}\oplus \Me_{\rm can,1}$ and to $b_{\rm can}$ corresponds a direct sum decomposition of $2$-divisible groups over an algebraic closure of $B(k)$ which can be naturally identified with a direct sum decomposition of ${\bf Z}_2$-modules $\Ml_n\otimes_{{\bf Z}_{(2)}} {\bf Z}_2=L_0\oplus L_1$, with $L_0$ and $L_1$ isotropic with respect to ${\got B}_n$. Here $\Ml_n$ and ${\got B}_n$ are as in Subsection 3.1 and we have applied Proposition 3.4 (a) with $S={\bf Z}_2$. 
\smallskip
One has a canonical short exact sequence $0\to \Me_{z,1}\to \Me_z\to \Me_{z,0}\to 0$ of $2$-divisible groups, where $\Me_{z,0}$ is \'etale and $\Me_{z,1}$ is of multiplicative type. As $\Me_{z,0,k}=\Me_{\rm can,0,k}$ and $\Me_{z,1,k}=\Me_{\rm can,1,k}$, we can identify canonically $\Me_{z,0}=\Me_{\rm can,0}$ and $\Me_{z,1}=\Me_{\rm can,1}$. As $z$ defines a $2$-torsion ${\rm Spf}\,W(k)$-valued point of ${\rm Spf}\,\Mo_y^{\rm big}$, the last short exact sequence splits after a push forward through the isogeny $2:\Me_{z,1}\to\Me_{z,1}$. By considering the pull back of $\Me_z$ to the spectrum of an algebraic closure of $B(k)$, we conclude that to $\Me_z$ and to the isomorphism $b_z:\Me_z\tilde\to \Me_z^{\rm t}$ corresponds uniquely (via the natural identification of the last paragraph) a ${\bf Z}_2$-lattice $\Ml_z$ of $\Ml_n\otimes_{{\bf Z}_{(2)}} {\bf Q}_2$ which has the following three properties:
\medskip
{\bf (i)} it contains $L_1$ and it is contained in ${1\over 2}L_1\oplus L_0$;
\smallskip
{\bf (ii)} we have a short exact sequence $0\to L_1\to \Ml_z\to L_0\to 0$ of ${\bf Z}_2$-modules;
\smallskip
{\bf (iii)} we have a perfect symmetric bilinear form ${\got B}_n:\Ml_z\times\Ml_z\to {\bf Z}_2$ whose reduction modulo $2{\bf Z}_2$ is alternating (cf. Fact 4.3.1). 
\medskip
As the quasi-canonical lift $z$ of $y$ is uniquely determined by the $2$-divisible group $\Me_z^s$ of $A_{W(k)}$ (cf. Serre--Tate deformation theory) and thus by $\Me_{z,B(k)}^s$ (cf. Tate's extension theorem), we get that $z$ is also uniquely determined by the ${\bf Z}_2$-lattice $\Ml_z$ of $\Ml_n\otimes_{{\bf Z}_{(2)}} {\bf Q}_2$. Thus the number of such lattices $\Ml_z$ is at least equal to the number $v2^d=v2^{{n(n-1)}\over 2}$ of quasi-canonical lifts $z\in\Mn_{W(k)}(W(k))$ of $y$. But the number of ${\bf Z}_2$-lattices of $\Ml_n\otimes_{{\bf Z}_{(2)}} {\bf Q}_2$ that have the above three properties is exactly $2^d=2^{{n(n-1)}\over 2}$, cf. Lemma 3.1.4. Thus $v2^d\le 2^d$ and therefore $v\le 1$. As $v\in {\bf N}^{\ast}$, we get that $v=1$. 
\smallskip
Moreover, there exists a unique quasi-canonical lift $z_{\rm can}\in\Mn_{W(k)}(W(k))$ of $y$ such that $\Ml_{z_{\rm can}}=\Ml_n\otimes_{{\bf Z}_{(2)}} {\bf Z}_2$ i.e., $z_{\rm can}^*(\Ma\times_{\Mn} \Mn_{W(k)})$ is the canonical lift of $A$.\endproof 
\medskip\noindent
{\bf 7.1.1. Corollary.} {\it If $y\in\Mn_{W(k)}(k)$ factors through the ordinary locus of $\Mn_k$, then the formal deformation space ${\rm Spf}\,\Mo_y$ of $y$ in $\Mn_{W(k)}$ is a formal subtorus of the formal torus defined by the formal deformation space ${\rm Spf}\,\Mo_y^{\rm big}$ of $y$ in $\Mm_{W(k)}$.}
\medskip\noindent
{\it Proof:} We know that ${\rm Spf}\,\Mo_y={\rm Spf}\,\Mo_y^{\rm n}$ is the translation of a formal subtorus of ${\rm Spf}\,\Mo_y^{\rm big}$ by a $2$-torsion ${\rm Spf}\,W(k)$-valued point, cf. proof of Proposition 7.1. But the existence of the lift $z_{\rm can}$ of $y$ (see above) implies that ${\rm Spf}\,\Mo_y$ contains the identity section of the formal torus ${\rm Spf}\,\Mo_y^{\rm big}$ over ${\rm Spf}\,W(k)$. The corollary follows from the last two sentences.\endproof
\medskip\smallskip\noindent
{\bf 7.2. Proof of Proposition 1.4.1.} Proposition 1.4.1 (a) follows from the end of the proof of Proposition 7.1 (or from Corollary 7.1.1). We now prove Proposition 1.4.1 (b). 
\smallskip
Let $y\in\Mn_{W(k)}(k)$. We use the notations of Subsection 6.2. Let $(\bar M_y,\phi_y)$ be the reduction modulo $2W(k)$ of $(M_y,\Phi_y)$. Let $l_y:L_{(2)}^\vee\otimes_{{\bf Z}_{(2)}} W(k)\tilde\to M_y$ be the inverse of the transpose (dual) of the isomorphism $j_y$ of Subsection 5.2. We view naturally $G_{{\bf Z}_2}$ as a closed subgroup scheme of ${\bf GL}_{L_{(2)}^\vee\otimes_{{\bf Z}_{(2)}} {\bf Z}_2}$. As $G_{{\bf Z}_2}$ is split, there exists a cocharacter $\mu_{\acute et}:{\bf G}_{m,{\bf Z}_2}\to G_{{\bf Z}_2}$ such that
$$l_y{\mu_{\acute et}}_{W(k)} l_y^{-1}:{\bf G}_{m,W(k)}\to\Mg_y$$ 
is $\Mg_y(W(k))$-conjugate to the cocharacter $\mu_y:{\bf G}_{m,W(k)}\to\Mg_y$ of Proposition 6.1. 
\smallskip
By replacing $l_y$ with its left composite with an element of $\Mg_y(W(k))$, we can assume that in fact we have  $l_y{\mu_{\acute et}}_{W(k)} l_y^{-1}=\mu_y$. Thus $l_y^{-1}\Phi_yl_y$ is a $\sg$-linear endomorphism of $L_{(2)}^\vee\otimes_{{\bf Z}_{(2)}} W(k)$ of the form $g_{\acute et}(1_{L_{(2)}^\vee}\otimes\sg)\mu_{\acute et W(k)}({1\over 2})$ for some element $g_{\acute et}\in G_{{\bf Z}_{(2)}}(W(k))$. As $(L_{(2)}^\vee\otimes_{{\bf Z}_{(2)}} W(k),(1_{L_{(2)}^\vee}\otimes\sg){\mu_{\acute et}}_{W(k)}({1\over 2}))$ is an ordinary $F$-crystal over $k$, there exists an element $g_y\in\Mg_y(W(k))$ such that $(M_y,g_y\Phi_y)$ is an ordinary $F$-crystal over $k$. 
\smallskip
We consider the open subscheme $\Mz_y$ of $\Mg_{y,k}$ with the property that for an element $g\in\Mg_y(W(k))$, the $F$-crystal $(M_y,g\Phi_y)$ over $k$ is ordinary if and only if the reduction $\bar g$ of $g$ modulo $2W(k)$ is a $k$-valued points of $\Mz_y$. The existence of $\Mz_k$ is a direct consequence of the existence of Hasse--Witt invariants: for $g\in\Mg_y(W(k))$, the $F$-crystal $(M_y,g\Phi_y)$ over $k$ is ordinary if and only if the $\sigma$-linear map $\bar M_y/\bar F^1_y\to \bar M_y/\bar F^1_y$ induced naturally by $\bar g\phi_y$ is a $\sigma$-linear automorphism (equivalently, its matrix representation with respect to a fixed basis of the $k$-vector space $\bar M_y/\bar F^1_y$ has a non-zero determinant). As $\bar g_y\in\Mz_y(k)$ (cf. end of the previous paragraph), $\Mz_y$ is an open, Zariski dense subscheme of the integral scheme $\Mg_{y,k}$. 
\smallskip
For each $h\in \Mg_y(W(k))$ that lifts a $k$-valued point of $\Mp_{y,k}$, we have $\Phi_y(h):=\Phi_y\circ h\circ \Phi_y^{-1}\in \Mg_y(W(k))$. Let $\Mp_y$ be a parabolic subgroup scheme of $\Mg_y$ whose special fibre is $\Mp_{y,k}$ and through which $\mu_y$ factors. The rule that associates to a pair $(h,u)\in \Mp_y(W(k))\times \Mu_y(W(k))$ the element $hu\Phi_y(h^{-1})\in\Mg_y(W(k))$ induces over $k$ an \'etale morphism $\Mp_{y,k}\times_{{\rm Spec}\,k} \Mu_{y,k}\to \Mg_{y,k}$. Thus there exists such a pair $(h,u)$ such that $h u\Phi_y(h^{-1})$ lifts a $k$-valued point of $\Mz_y$. As the $F$-crystals $(M_y,hu\Phi_y(h^{-1})\Phi_y)$ and $(M_y,u\Phi_y)$ are isomorphic, we conclude that there exists $u\in \Mu_y(W(k))$ such that $(M_y,u\Phi_y)$ is ordinary. In other words, $\Mz_y\cap \Mu_{y,k}$ is an open, Zariski dense subscheme of the integral scheme $\Mu_{y,k}$.
\smallskip
We write $\Mu_{y,k}={\rm Spec}\,S_y$. Let $\bar u_{\rm mod}\in\Mu_y(S_y)$ be defined by the closed embedding ${\rm Spec}\,S_y=\Mu_{y,k}\hookrightarrow \Mu_y$. It makes sense to speak about the pair $(\bar M_y\otimes_k S_y,\bar u_{\rm mod}(\phi_y\otimes\Phi_{S_y}))$ being ordinary at some point of $\Mu_{y,k}$; the set of such points is the open subscheme $\Mz_y\cap \Mu_{y,k}$ of $\Mu_{y,k}$. From the end of the previous paragraph we get that $(\bar M_y\otimes_k S_y,\bar u_{\rm mod}(\phi_y\otimes\Phi_{S_y}))$ is ordinary at the generic point of $\Mu_{y,k}$.
\smallskip
As the reduction modulo $2\Mr$ of the universal element $u_{\rm univ}\in\Mu_y(\Mr)$ of Subsection 6.3 is the pull back of $\bar u_{\rm mod}$ via the dominant morphism ${\rm Spec}\,\Mr/2\Mr\to \Mu_{y,k}$, from the last sentence and from the description of $\Mc_y$ in Subsection 6.4, we get that (the truncated Barsotti--Tate group of level $1$ of) the abelian scheme $A_{\Mr/2\Mr}^\prime$ over $\Mr/2\Mr$, is generically ordinary. As we have $n_{\Mr/2\Mr}^*(\Ma_{\Mn_{k,\rm red}})=A_{\Mr/2\Mr}^\prime$ (cf. (15)), we conclude that the ordinary locus of $\Mn_k$ specializes to the $k$-valued point of $\Mn_k$ defined by $y$. As $y$ is arbitrary, Proposition 1.4.1 (b) follows.\endproof
\medskip\smallskip\noindent
{\bf 7.3. End of the proof of Theorem 1.4 (b).} From Propositions 7.1 and 1.4.1 we get that the smooth locus of $\Mn_k/H_0$ is Zariski dense in $\Mn_k/H_0$. The connected components of $\Mn_k/H_0$ are irreducible, cf. Proposition 6.7. Thus the connected components of $\Mn_k/H_0$ are generically smooth and irreducible. Moreover the reduced scheme of $\Mn_k/H_0$ is smooth over $k$ (cf. Proposition 6.7) and thus normal. Therefore from a lemma of Hironaka (see [19, Part II, Prop. (5.12.8)]) applied to all local rings of $\Mn_{W(k)}/H_0$ at $k$-valued points and to their element $2$,  we get that $\Mn_k/H_0$ is normal and thus (cf. Proposition 6.7) a smooth $k$-algebra. This implies that $\Mn_{W(k)}/H_0$ is a smooth $W(k)$-scheme. Therefore $\Mn$ is regular and formally smooth over $O_{(v)}$ i.e., Theorem 1.4 (b) holds.\endproof
\medskip\noindent
{\it Acknowledgments.} We would like to thank the Univeristy of Arizona, Binghamton University, and TIFR--Mumbai for good working conditions and P. Deligne for some comments. We would also like to thank O. Gabber for several suggestions related to Section 3. We would also like to thank the referee for many valuable comments. This research was partially supported by the NSF grant DMS \#0900967.
\bigskip
\centerline{\bigsll {\bf References}}
\bigskip
\item{[1]} M. Artin, {\sl Algebraization of formal moduli. I}, Global Analysis (Papers in Honor of K. Kodaira),  21--71, Univ. Tokyo Press, Tokyo, 1969.
\item{[2]} M. Artin, {\sl Versal deformations and algebraic stacks}, Invent. Math. {\bf 27} (1974),  165--189.
\item{[3]} P. Berthelot, {\it Cohomologie cristalline des sch\'emas de caract\'eristique $p>0$}, Lecture Notes in Math., Vol. {\bf 407}, Springer-Verlag, Berlin-New York, 1974.
\item{[4]} P. Berthelot, L. Breen, and W. Messing, {\it Th\'eorie de Dieudonn\'e cristalline II}, Lecture Notes in Math., Vol. {\bf 930}, Springer-Verlag, Berlin, 1982.\item{[5]} P. Berthelot and W. Messing, {\it Th\'eorie de Dieudonn\'e cristalline III}, The Grothendieck Festschrift, Vol. I, 173--247, Progr. in Math., Vol. {\bf 86}, Birkh\"auser Boston, Boston, MA, 1990. 
\item{[6]} A. Borel, {\it Linear algebraic groups}, Grad. Texts in Math., Vol. {\bf 126}, Springer-Verlag, New York, 1991.
\item{[7]} N. Bourbaki, {\it Lie groupes and Lie algebras. Chapters 4--6}, Elements of Mathematics (Berlin), Springer-Verlag, Berlin, 2002.
\item{[8]} B. Conrad, O. Gabber, and G. Prasad, {\it Pseudo-reductive groups},
New Mathematical Monographs, Vol. {\bf 17}, Cambridge University Press, Cambridge, 2010. 
\item{[9]} P. Deligne, {\it Travaux de Shimura}, S\'em. Bourbaki, Exp. no 389, 123--165, Lecture Notes in Math., Vol. {\bf 244}, Springer-Verlag, Berlin, 1971.
\item{[10]} P. Deligne, {\it Vari\'et\'es de Shimura: interpretation modulaire, et techniques de construction de mod\`eles canoniques}, Automorphic forms, representations and $L$-functions (Oregon State Univ., Corvallis, OR, 1977), Part 2,   247--289, Proc. Sympos. Pure Math., {\bf 33}, Amer. Math. Soc., Providence, RI, 1979.
\item{[11]} P. Deligne, {\it Hodge cycles on abelian varieties}, Hodge cycles, motives, and Shimura varieties, 9--100, 
Lecture Notes in Math., Vol. {\bf 900}, Springer-Verlag, Berlin-New York, 1982.
\item{[12]} M. Demazure, A. Grothendieck, et al. {\it Sch\'emas en groupes. Vols. II and III}, S\'eminaire de G\'eom\'etrie Alg\'ebrique du Bois Marie 1962/64 (SGA 3), Lecture Notes in Math., Vol. {\bf 152--153}, Springer-Verlag, Berlin-New York, 1970.
\item{[13]} V. G. Drinfeld, {\it Elliptic modules}, (Russian), Mat. Sb. (N.S.) {\bf 94 (136)} (1974),  594--627, 656. 
\item{[14]} G. Faltings, {\it Integral crystalline cohomology over very ramified valuation rings}, J. Amer. Math. Soc. {\bf 12} (1999), no. 1,  117--144.
\item{[15]} G. Faltings and C.-L. Chai, {\it Degeneration of abelian varieties}, Ergebnisse der Math. und ihrer Grenzgebiete (3), Vol. {\bf 22}, Springer-Verlag, Heidelberg, 1990.
\item{[16]} J.-M. Fontaine, {\it Groupes $p$-divisibles sur les corps locaux}, Ast\'erisque {\bf 47/48}, Soc. Math. de France, Paris, 1977.
\item{[17]} J.-M. Fontaine, {\it Sur certain types de repr\'esentations p-adiques du groupe de Galois d'un corps local; construction d'un anneau de Barsotti--Tate}, Ann. of Math. (2) {\bf 115} (1982), no. 3,  529--577. 
\item{[18]} A. Grothendieck, {\it \'El\'ements de g\'eom\'etrie alg\'ebrique. III}, \'Etude cohomologique des faisceaux coh\'erents. I. (French) 
Inst. Hautes \'Etudes Sci. Publ. Math. No. {\bf 11} 1961, 167 pp.
\item{[19]} A. Grothendieck et al., {\it \'El\'ements de g\'eom\'etrie alg\'ebrique. IV. \'Etude locale des sch\'emas et des morphismes de sch\'ema (Seconde Partie \'et Quatri\`eme Partie)}, Inst. Hautes \'Etudes Sci. Publ. Math., Vols. {\bf 24} (1965) and {\bf 32} (1967). 
\item{[20]} M. Harris and R. Taylor, {\it The geometry and cohomology of some simple Shimura varieties}, Annals of Mathematics Studies, Vol. {\bf 151}, Princeton Univ. Press, Princeton, NJ, 2001.
\item{[21]} S. Helgason, {\it Differential geometry, Lie groups, and symmetric spaces}, Academic Press, New-York, 1978.
\item{[22]} L. Illusie, {\it D\'eformations des groupes de Barsotti--Tate (d'apr\`es A. Grothendieck)}, Seminar on arithmetic bundles: the Mordell conjecture (Paris, 1983/84),  151--198, Ast\'erisque {\bf 127}, Soc. Math. de France, Paris, 1985.
\item{[23]} N. Katz, {\it Serre--Tate local moduli},  Algebraic surfaces (Orsay, 1976--78),   138--202, Lecture Notes in Math., Vol. {\bf 868}, Springer, Berlin-New York, 1981.
\item{[24]} M. Kisin, {\it Integral models for Shimura varieties of abelian type}, J. Amer. Math. Soc. {\bf 23} (2010),  no. 4, 967--1012.
\item{[25]} M.-A. Knus, A. Merkurjev, M. Rost, and J.-P. Tignol, {\it The book of involutions}, Amer. Math. Soc. Colloquium Publ., Vol. {\bf 44}, Amer. Math. Soc., Providence, RI, 1998. 
\item{[26]} R. E. Kottwitz, {\it Points on some Shimura varieties over finite fields}, J. Amer. Math. Soc. {\bf 5} (1992), no. 2,  373--444.
\item{[27]} R. Langlands and M. Rapoport, {\it Shimuravariet\"aten und Gerben}, J. Reine Angew. Math. {\bf 378} (1987),  113--220.
\item{[28]} H. Matsumura, {\it Commutative algebra. Second edition}, Mathematics Lecture Note Series, Vol. {\bf 56}, Benjamin/Cummings Publishing Co., Inc., Reading, Mass., 1980.
\item{[29]} W. Messing, {\it The crystals associated to Barsotti--Tate groups, with applications to abelian schemes}, Lecture Notes in Math., Vol. {\bf 264}, Springer-Verlag, Berlin-New York, 1972.
\item{[30]} J. S. Milne, {\it The points on a Shimura variety modulo a prime of good reduction}, The Zeta functions of Picard modular surfaces,  153--255, Univ. Montr\'eal, Montreal, Quebec, 1992.
\item{[31]} J. S. Milne, {\it Shimura varieties and motives}, Motives (Seattle, WA, 1991),  447--523, Proc. Sympos. Pure Math., Vol. {\bf 55}, Part 2, Amer. Math. Soc., Providence, RI, 1994.
\item{[32]} Y. Morita, {\it Reduction mod $\got{B}$ of Shimura curves}, Hokkaido Math. J. {\bf 10} (1981), no. 2,  209--238.
\item{[33]} D. Mumford, {\it Abelian varieties. Second Edition}, Tata Inst. of Fund. Research Studies in Math., No. {\bf 5}, Published for the Tata Institute of Fundamental Research, Bombay; Oxford Univ. Press, London, 1970 (reprinted 1988).
\item{[34]} D. Mumford, J. Fogarty, and F. Kirwan, {\it Geometric invariant theory. Third edition}, Ergebnisse der Math. und ihrer Grenzgebiete (2), Vol. {\bf 34}, Springer-Verlag, Berlin, 1994.
\item{[35]} R. Noot, {\it Models of Shimura varieties in mixed characteristic}, J. Algebraic Geom. {\bf 5} (1996), no. 1,  187--207.
\item{[36]} J.-P. Serre, {\it Galois Cohomology}, Springer-Verlag, Berlin, 1997.
\item{[37]} G. Shimura, {\it On analytic families of polarized abelian varieties and automorphic functions}, Ann. of Math. {\bf 78} (1963), no. 1,  149--192.
\item{[38]} J. Tits, {\it Reductive groups over local fields}, Automorphic forms, representations and $L$-functions (Oregon State Univ., Corvallis, OR, 1977), Part 1,   29--69, Proc. Sympos. Pure Math., {\bf 33}, Amer. Math. Soc., Providence, RI, 1979.
\item{[39]} A. Vasiu, {\it Integral canonical models of Shimura varieties of preabelian type}, Asian J. Math. {\bf 3} (1999), no. 2,  401--518.
\item{[40]} A. Vasiu, {\it A purity theorem for abelian schemes}, Michigan Math. J. {\bf 52} (2004), no. 1,  71--81. 
\item{[41]} A. Vasiu, {\it Integral canonical models of unitary Shimura varieties}, Asian J. Math. {\bf 12} (2008), no. 2,  151--176.
\item{[42]} A. Vasiu, {\it Good reductions of Shimura varieties of preabelian type in arbitrary mixed characteristic, Part I,}
48 pages, available at http://arxiv.org/abs/0707.1668.
\item{[43]} A. Vasiu and T. Zink, {\it Purity results for $p$-divisible groups and abelian schemes over regular bases of mixed characteristic}, Documenta Math. {\bf 15} (2010), 571--599.
\item{[44]} T. Zink, {\it Isogenieklassen von Punkten von Shimuramannigfaltigkeiten mit Werten in einem endlichen K\"orper}, Math. Nachr. {\bf 112} (1983),  103--124.
\medskip
\noindent
Adrian Vasiu

\noindent
Department of Mathematical Sciences, Binghamton University,

\noindent
Binghamton, P. O. Box 6000, New York 13902-6000, U.S.A. 

\noindent
fax: 1-607-777-2450, e-mail: adrian@math.binghamton.edu 

\end